\documentclass[10pt]{article}
\usepackage{amssymb, amsmath, cancel}
\usepackage{graphicx,psfrag,epsf,booktabs}
\usepackage{enumerate}
\usepackage{natbib}
\usepackage{url} 

\newcommand{\blind}{0}

\usepackage{amsthm,amsmath,amsfonts,amssymb}
\usepackage{dsfont,bm,multirow,stmaryrd}
\usepackage{amsthm}
\usepackage{multirow}
\usepackage{natbib}
\usepackage{bm}
\usepackage{hyperref}
\usepackage[all]{xy}
\usepackage{graphicx,caption, subcaption}
\usepackage{float}
\usepackage{color}

\numberwithin{equation}{section}

\theoremstyle{plain}
\newtheorem{theorem}{Theorem}[section]
\newtheorem{corollary}{Corollary}[section]
\newtheorem{lemma}{Lemma}[section]
\newtheorem{prop}{Proposition}[section]

\theoremstyle{remark}
\newtheorem{definition}[theorem]{Definition}

\newcommand{\bed}{\begin{definition}}
\newcommand{\eed}{\end{definition}}
\renewcommand{\b}{\mathbf b}

\newcommand{\rom}[1]{\uppercase\expandafter{\romannumeral #1\relax}}

\usepackage{bbm}

\newcommand{\bitem}{\begin{itemize}}
\newcommand{\eitem}{\end{itemize}}

\newcommand{\goto}{\rightarrow}

\newcommand{\beqn}{\begin{equation}}
\newcommand{\eeqn}{\end{equation}}
\newcommand{\balign}{\begin{align}}
\newcommand{\ealign}{\end{align}}

\newcommand{\tPi}{\tilde{\Pi}}

\newcommand{\tG}{\tilde{G}}
\newcommand{\tpi}{\tilde{\pi}}
\newcommand{\tH}{\tilde{h}}
\newcommand{\Tr}{\text{Tr}}
\newcommand{\tOmega}{\overline{\Omega}}

\newcommand{\s}{\sigma}

\newcommand{\oOmega}{\overline{\Omega}}

\newcommand{\tr}{\mathrm{tr}}
\newcommand{\E}{\mathbb{E}}
\newcommand{\V}{\text{Var}}
\newcommand{\cov}{\text{Cov}}
\usepackage{amssymb}
\newcommand{\beq}{\begin{equation}}
\newcommand{\eeq}{\end{equation}}

\newcommand{\diag}{\mathrm{diag}}

\allowdisplaybreaks
\newcommand*\xbar[1]{%
   \hbox{%
     \vbox{%
       \hrule height 0.5pt 
       \kern0.5ex
       \hbox{%
         \kern-0.1em
         \ensuremath{#1}%
         \kern-0.1em
       }%
     }%
   }%
}

\usepackage{listings}
\begin{document}
	
\def\spacingset#1{\renewcommand{\baselinestretch}%
	{#1}\small\normalsize} \spacingset{1}


\if0\blind
{
	\title{\bf Power Enhancement and Phase Transitions for Global Testing of the Mixed Membership Stochastic Block Model}
	\author{Louis Cammarata and Zheng Tracy Ke\\
		\vspace{0.4 em} 
	 Harvard University} 
 \date{}

	\maketitle
} \fi

\if1\blind
{
	\bigskip
	\bigskip
	\bigskip
	\begin{center}
		{\LARGE\bf Power Enhancement and Phase Transitions for Global Testing of the Mixed Membership Stochastic Block Model}
	\end{center}
	\medskip
} \fi

\begin{abstract}
The mixed-membership stochastic block model (MMSBM) is a common model for social networks. Given an $n$-node symmetric network generated from a $K$-community MMSBM, we would like to test $K=1$ versus $K>1$. We first study the degree-based $\chi^2$ test and the orthodox Signed Quadrilateral (oSQ) test. These two statistics estimate an order-2 polynomial and an order-4 polynomial of a ``signal'' matrix, respectively. We derive the asymptotic null distribution and power for both tests. However, for each test, there exists a parameter regime where its power is unsatisfactory. 
It motivates us to propose a power enhancement (PE) test to combine the strengths of both tests. We show that the PE test has a tractable null distribution and improves the power of both tests.  
To assess the optimality of PE, we consider a randomized setting, where the $n$ membership vectors are independently drawn from a distribution on the standard simplex. 
We show that the success of global testing is governed by
a quantity $\beta_n(K,P,h)$, which depends on the community structure matrix $P$ and the mean vector $h$ of memberships. 
For each given $(K, P, h)$, a test is called {\it optimal} if it distinguishes two hypotheses when $\beta_n(K, P,h)\to\infty$. A test is called {\it optimally adaptive} if it is optimal for all $(K, P, h)$. 
We show that the PE test is optimally adaptive, while many existing tests are only optimal for some particular $(K, P, h)$, hence, not optimally adaptive.  
\end{abstract}

\noindent%

\maketitle

\vfill

\spacingset{1.1}

\newpage

\section{Introduction} \label{sec:intro}
Statistical analysis of large social networks has received much recent attention. 
In this paper, we are interested in testing whether an undirected network has one community or multiple communities (a.k.a., global testing). This problem has several applications: it is useful in the design of stopping rules in recursive community detection \citep{SCC2021}; it can also be applied to ego-networks to measure the neighborhood diversity of individual nodes \citep{gao2017testing,SCC2021}. 
Theoretically, a study of the global testing problem (especially its lower bound) also provides valuable insights for related problems, such as community detection \citep{zhang2016minimax}, mixed membership estimation \citep{Mixed-SCORE}, and estimation of the number of communities \citep{EstK}. 

For an undirected network with $n$ nodes, the adjacency matrix $A=(A_{ij})_{1\leq i,j\leq n}$ is a symmetric matrix, where 
\begin{equation} \label{model1a} 
A_{ij} = 
\left\{
\begin{array}{ll} 
1, &\quad  \mbox{if $i\neq j$ and there is an edge between node $i$ and node $j$},  \\
0, &\quad   \mbox{otherwise}.    \\
\end{array}\right.
\end{equation} 
Assuming that the network has $K$ perceivable communities, we model $A$ with the Mixed-Membership Stochastic Block Model (MMSBM) \citep{airoldi2009mixed} as follows. The mixed-membership vector of node $i$ is a weight vector $\pi_i \in \mathbb{R}^K$ such that $\pi_i(k)\geq 0$ and $\sum_{k=1}^K\pi_i(k)=1$, with $\pi_i(k)$ denoting the `weight' that node $i$ puts on community $k$. For a symmetric nonnegative matrix $P \in \mathbb{R}^{K\times K}$ that models the community structure,  we assume that the upper triangle of $A$ contains independent Bernoulli variables, where   
\begin{equation} \label{model1b} 
\mathbb{P}(A_{ij}=1) = \pi_i'P\pi_j,   \qquad 1 \leq i < j \leq n.  
\end{equation} 

The well-known Stochastic Block Model (SBM) and Erd\"os-Renyi (ER) model are special cases of MMSBM.  When all $\pi_i$'s are degenerate (meaning that 
one entry of $\pi_i$ is $1$ and all other entries are $0$), MMSBM reduces to SBM; furthermore, if $K = 1$, then SBM reduces to ER.  

The global testing problem is formulated as testing between the two hypotheses: 
\begin{equation} \label{gt-0} 
H_0^{(n)}: K = 1 \qquad \text{versus} \qquad H_1^{(n)}: K > 1. 
\end{equation}
In MMSBM, write $\Omega=\Pi P\Pi'$, with $\Pi=[\pi_1,\pi_2,\ldots,\pi_n]'\in\mathbb{R}^{n\times K}$. We call $\Omega$ the Bernoulli probability matrix. It follows that, for some $\alpha_n>0$, 
\beq \label{gt-1A}
\mathbb{E}A=\Omega - \mathrm{diag}(\Omega), \qquad \mbox{with}\quad \Omega= \begin{cases}
\alpha_n {\bf 1}_n{\bf 1}_n', &\mbox{under }H_0^{(n)},\\
\Pi P\Pi', & \mbox{under }H_1^{(n)}.
\end{cases}
\eeq

The signals to separate two hypotheses are captured by the following $n\times n$ matrix:
\beq \label{Define-tOmega}
\widetilde{\Omega} = \Omega  - \tilde{\alpha}_n {\bf 1}_n{\bf 1}_n', \qquad\mbox{where}\quad \tilde{\alpha}_n=n^{-2}({\bf 1}_n'\Omega {\bf 1}_n). 
\eeq
The null hypothesis holds if and only if $\widetilde{\Omega}$ is a zero matrix. We will investigate testing ideas that aim to estimate a polynomial of (the entries of) $\widetilde{\Omega}$:
\begin{itemize}
\item We consider the {\it degree-based $\chi^2$} statistic, which targets on estimating ${\bf 1}_n'\widetilde{\Omega}^2{\bf 1}_n$.  
\item For each $m\geq 3$, we consider the order-$m$ {\it orthodox Signed Polygon} statistic, which targets on estimating $\tr(\widetilde{\Omega}^m)$.  
\end{itemize}
Here, ${\bf 1}_n'\widetilde{\Omega}^2{\bf 1}_n$ is a {\it second order polynomial} of $\widetilde{\Omega}$, and $\tr(\widetilde{\Omega}^m)$ is an {\it $m$-th order polynomial} of $\widetilde{\Omega}$, for $m\geq 3$. 
There is no natural testing idea based on estimating a first order polynomial of $\widetilde{\Omega}$. For example, ${\bf 1}_n'\widetilde{\Omega}{\bf 1}_n$ is always equal to zero and hence useless for global testing; $\tr(\widetilde{\Omega})$ is hard to estimate, partially because the diagonals of $A$ are always zero (i.e., self-edges are not allowed).

We study the asymptotic performances of the above tests. We also derive information theoretic lower bounds for this global testing problem. By comparing the upper/lower bounds, we discover: (i) None of the above tests can attain the lower bound across all parameter regimes;
(ii) In some parameter regimes, the degree-based $\chi^2$ test (abbreviated as the $\chi^2$ test) is optimal; and in the remaining parameter regimes, the order-4 orthodox Signed Polygon test (abbreviated as the oSQ test) is optimal. 
This motivates us to design a new test statistic to combine the strengths of the $\chi^2$ test and the oSQ test.  

We propose the  {\it Power Enhancement} (PE) test. It is inspired by a key result about the joint distribution of the $\chi^2$ test statistic and the oSQ test statistic: Under the null hypothesis, they jointly converge to a bivariate normal distribution with a covariance matrix $I_2$; especially, the two test statistics are asymptotically uncorrelated. 
The PE test statistic is defined as the sum of squares of these two test statistics, and it converges to a $\chi_2^2(0)$ distribution under the null hypothesis. 
Therefore, we can conveniently control the level of the PE test.

To assess the power and optimality of the PE test, we adopt the phase transition framework in Jin et al. \citep{SP2021}. For arbitrary parameters $(K, P)$ and distribution $F$ on the probability simplex of $\mathbb{R}^K$, writing $h=\mathbb{E}_{\pi\sim F}[\pi]$, we consider the following pair of hypotheses:  
\beq \label{pair-RMM}
\Omega= \left\{\begin{array}{lll}
\alpha_0 {\bf 1}_n{\bf 1}_n', & \mbox{where }\alpha_0=h'Ph, &\qquad \mbox{under }H_0^{(n)},\\
\Pi P\Pi', & \mbox{where }\pi_i\overset{iid}{\sim} F, &\qquad \mbox{under }H_1^{(n)}.
\end{array}\right. 
\eeq
As $n\to\infty$, we fix $K$ and allow $(P, h)$ to depend on $n$ (i.e., we consider a sequence of $(P_n, h_n)$ indexed by $n$). Our lower bound result tells when the chi-square distance between two hypotheses converges to $0$ for {\it every} $(K, P_n, h_n)$. In particular, we identify a quantity (as before, $\alpha_0=h_n'P_nh_n$)
\beq \label{SNR-combined}
\beta_n(K, P_n,h_n) \equiv \max\bigl\{ n^{3/2}\alpha_{0}^{-1}\|P_n h_n-\alpha_{0}{\bf 1}_K\|^2,\;\;  n^2\alpha_{0}^{-2}\|P_n - \alpha_{0}{\bf 1}_K{\bf 1}_K' \|^4  \bigr\}, 
\eeq
such that the chi-square distance between two hypotheses tends to 0 if $\beta_n(K, P_n,h_n)\to 0$. 
We call the parameter regimes where $\beta_n(K, P_n, h_n)\to 0$ the {\it Region of Impossibility}. In this region, the two hypotheses are asymptotically inseparable. We call the parameter regimes where $\beta_n(K, P_n, h_n)\to\infty$ the {\it Region of Possibility}. A test is called {\it optimally adaptive} if it is able to distinguish two hypotheses for {\it any} $(K, P_n, h_n)$ in the Region of Possibility. We show that the PE test is optimally adaptive.


\subsection{Related literature}
The likelihood ratio test (LRT) was studied by Mossel et al. \citep{mossel2015reconstruction} and Banerjee and Ma \citep{banerjee2017optimal} for a special case of $K=2$, $P_n=(a_n-b_n)I_2+b_n{\bf 1}_2{\bf 1}_2'$ and $h_n=(1/2, 1/2)'$. 
The LRT may be generalized to other $(K, P_n, h_n)$, but it requires prior knowledge of parameters in the alternative hypothesis. By the Neyman-Pearson lemma, the LRT has the highest power; however, it not a polynomial-time test. The tests we study here, $\chi^2$, oSQ and PE, need no prior knowledge of the parameters, and are polynomial-time tests.  

The eigenvalue-based tests were also studied before. For example,  Lei \citep{lei2016goodness} used the maximum singular value of the centered and rescaled adjacency matrix as test statistic. However, the eigenvalue-based tests are not optimally adaptive: their SNRs are linked to the second term in \eqref{SNR-combined}; hence, in the Region of Possibility, for those $(K, P_n, h_n)$ such that  the first term in \eqref{SNR-combined} $\to\infty$ but the second term $\to 0$, the eigenvalue-based tests are unable to separate two hypotheses.

Arias-Castro and Verzelen \citep{arias2014community} considered the testing of a planted clique model v.s. the ER model. The planted clique model can be viewed as a special case of MMSBM with $K=2$, $P_n=b_n{\bf 1}_2{\bf 1}_2' + (a_n-b_n)e_1e_1'$ and $h=(\epsilon_n, 1-\epsilon_n)'$, where $\epsilon_n=o(1)$ and $a_n>b_n>0$. They derived the optimal detection boundaries for many different cases of $(n, a_n, b_n, \epsilon_n)$ and proposed optimal tests. When $b_n$ is unknown and  $n\epsilon_n\gg n^{3/2}$, they used the $\chi^2$ test (called the {\it degree variance test} in their paper). Our result about the SNR of the $\chi^2$ test agrees with their result (see their Table 1, the bottom right cell) 
in this special case. 
However, there are major differences between two papers: 
First, they focused on a particular $(K, P_n, h_n)$, for which the first term in \eqref{SNR-combined} always dominates, so that the $\chi^2$ test alone is enough to achieve optimality (provided that $n\epsilon_n\gg n^{3/2}$). In contrast, we seek to find an {\it optimally adaptive} test that works for a broad collection of $(K, P_n, h_n)$, where the power enhancement idea is crucial. Second, we focus on $\epsilon_n\asymp 1$, but their main interest was in $\epsilon_n=o(1)$. In their setting, they could take advantage of sparsity by using the scan tests, which is unnecessary in our setting. 
Last, we provide the asymptotic null distribution for the $\chi^2$ test, which was not given in \citep{arias2014community}.

The cycle count statistics were also studied in recent literature \citep{banerjee2018contiguity, banerjee2018asymptotic, bubeck2016testing,Gao2017Test2,OGC,SP2021,lu2020contextual}. 
Our oSQ test is the same as the order-4 {\it signed-cycle} statistic introduced by Bubeck et al. \citep{bubeck2016testing} (also, see Banerjee \citep{banerjee2018contiguity}). Under a 2-community SBM model (and the related contextual SBM model and Gaussian covariance model), Banerjee and Ma \cite{banerjee2018asymptotic} and Lu and Sen \citep{lu2020contextual} derived  asymptotic distributions of order-$m$ cycle count statistics for a general $m$.
However, these works focused on the special case of $K=2$, $P_n=(a_n-b_n)I_2+b_n{\bf 1}_2{\bf 1}_2'$ and $h=(1/2, 1/2)'$, in which the oSQ test alone is enough to attain optimality. We seek to find an optimally adaptive test that works for rather arbitrary $(K, P_n, h_n)$, where we do need to combine oSQ with the $\chi^2$ test to achieve optimality. Moreover, these works only studied the asymptotic behavior of cycle counts, but we study the joint distribution of the cycle count and the $\chi^2$ statistic (this is one of our key results that inspires the PE test). Last, none of these works revealed the phase transition in \eqref{SNR-combined}.

Jin et al. \citep{SP2021} studied the phase transition of global testing under the Degree-Corrected Mixed Membership (DCMM) model \citep{Mixed-SCORE}, a model more general than the MMSBM considered here. They proposed the Signed Quadrilateral (SQ) test and showed that it is optimally adaptive. Although our model is a special case of DCMM with no degree heterogeneity, the phase transition and the optimal test are different. Restricting from DCMM to MMSBM, the prior knowledge of `no degree heterogeneity' brings additional signals for separating two hypotheses. For example, under DCMM, one can construct a pair of null and alternative hypotheses such that the expected degree of each node is perfectly matched under two hypotheses \citep{SP2021}, and so the $\chi^2$ statistic contains no signals and power enhancement is useless. But such `degree-matched' hypothesis pairs do not exist when we restrict to MMSBM; for MMSBM, power enhancement is crucial for achieving optimality. 
\citep{SP2021} showed that the Region of Possibility and Region of Impossibility for global testing under DCMM are determined by $
\gamma_n=|\lambda_2(\Omega)|/\sqrt{\lambda_1(\Omega)}$. 
This quantity restricted to MMSBM is different from $\beta_n$ in \eqref{SNR-combined}. 
Moreover, the SQ test in \citep{SP2021} is also different from our oSQ test, so their results about the asymptotic behavior of the SQ test cannot imply our results about the oSQ test.  

\subsection{Content}
We have made several contributions in this paper:
\begin{itemize}
\item We derive the phase transitions for global testing under MMSBM, where the Region of Impossibility and Region of Possibility are determined by the simple quantity $\beta_n(K,P_n,h_n)$.
\item We study the (degree-based) $\chi^2$ test and the oSQ test. For each test statistic, we derive its asymptotic distribution under the null hypothesis and SNR under the alternative hypothesis. We also derive the asymptotic {\it joint null} distribution of two test statistics.  
\item We propose the Power Enhancement (PE) test to combine the strengths of the oSQ test and the $\chi^2$ test while overcoming their respective limitations.
\item We show that the PE test statistic has an asymptotic null distribution of $\chi^2_2(0)$. We also show that the PE test is optimally adaptive. In comparison, several popular tests are not optimally adaptive. 
\end{itemize}
Below, in Section~\ref{sec:3tests}, we formally introduce the $\chi^2$, oSQ and PE tests, and explain how PE combines strengths of the other two tests. In Section~\ref{sec:main}, we present the main theoretical results, including asymptotic properties of three test statistics, lower bounds and phase transitions. Section~\ref{sec:identifiability} is of independent interest, where we discuss the identifiability of parameters of MMSBM, especially, the identifiability of $K$. Section~\ref{sec:simul} contains simulation studies, and Section~\ref{sec:conclude} concludes the paper.

\section{Three test statistics}  \label{sec:3tests}

Recall that $A$ is the adjacency matrix. Let $d = A {\bf 1}_n$ be the vector of degrees, and $\bar{d} = \frac{1}{n} \sum_{i = 1}^n d_i$ be the average degree.  Under the null hypothesis, $\Omega = \alpha_n{\bf 1}_n{\bf 1}_n'$, and a good estimate of $\alpha_n$ is  
\begin{equation} \label{Definehata} 
\hat{\alpha}_n = [n(n-1)]^{-1} {\bf 1}_n' A {\bf 1}_n = (n-1)^{-1} \bar{d}.  
\end{equation}
Under the null hypothesis, $(d_i-\bar{d})/\sqrt{(n-1)\hat{\alpha}_n(1-\hat{\alpha}_n)} \approx \mathcal{N}(0,1)$, for each $1\leq i\leq n$. 
Aggregating these terms for all $i$ gives rise to the \textit{degree-based $\chi^2$ test statistic} 
(also known as the {\it degree of variance statistic} \citep{arias2014community}; throughout this paper, we call it the $\chi^2$ test for short): 
\begin{equation} \label{chi1} 
X_n = [(n-1) \hat{\alpha}_n (1 - \hat{\alpha}_n)]^{-1}\sum_{i = 1}^n(d_i - \bar{d})^2.    \footnote{In the rare event that $\hat{\alpha}_n=0$ or $\hat{\alpha}_n=1$, we replace $\hat{\alpha}_n$ by $\frac{2}{n(n-1)}$ and $\frac{n(n-1)-2}{n(n-1)}$, respectively, to make the $\chi^2$ test statistic well defined. Similar operations apply to the oSQ and PE tests.}
\end{equation}
The $\chi^2$ test looks for evidence against the null hypothesis through degree heterogeneity. 
Since all nodes have the same expected degree under the null hypothesis, a significant degree heterogeneity provides a strong evidence against the null.  
By some simple calculations, we find that 
\beq \label{chi1-equivalent} 
[(n-1) \hat{\alpha}_n (1 - \hat{\alpha}_n)](X_n-n) = \sum_{i_1,i_2,i_3 \text{ (distinct)}}(A_{i_1i_2}-\hat{\alpha}_n)(A_{i_2i_3}-\hat{\alpha}_n). 
\eeq 
Recall that we have introduced $\widetilde{\Omega}$ in \eqref{Define-tOmega}. The matrix $A-\hat{\alpha}_n{\bf 1}_n{\bf 1}_n'$ is a stochastic proxy of $\widetilde{\Omega}$. Therefore, the right hand side above is approximately $\sum_{i_1,i_2,i_3\text{ (distinct)}}\widetilde{\Omega}_{i_1i_2}\widetilde{\Omega}_{i_2i_3}\approx {\bf 1}_n'\widetilde{\Omega}^2{\bf 1}_n$. This suggests that $[(n-1) \hat{\alpha}_n (1 - \hat{\alpha}_n)](X_n-n)$ is an estimate of ${\bf 1}_n'\widetilde{\Omega}^2{\bf 1}_n$, under the alternative hypothesis.

The {\it orthodox Signed Polygon} is a family of statistics that extends the {\it Signed Triangle} statistic \citep{bubeck2016testing}, where for $m = 3, 4, \ldots$,  the $m$-th order statistic in the family is defined as 
\beq  \label{oSP}
    U_n^{(m)} = \sum_{i_1, i_2, \ldots i_m \text{ (distinct)}}  (A_{i_1 i_2} - \hat{\alpha}_n) (A_{i_2i_3} - \hat{\alpha}_n) \ldots  (A_{i_mi_1} - \hat{\alpha}_n).
\eeq
Centering each $A_{ij}$ by $\hat{\alpha}_n$ is reasonable, because $\mathbb{E}[A_{ij}] = \alpha_n$ under the null hypothesis.   
It was noted in \citep{SP2021} that the Signed Polygon
statistic may experience signal cancellation if $m$ is odd, but it can avoid signal cancellation if $m$ is even.  For this reason, it is preferred to only consider the even order statistics in the family. In this paper, we focus our discussion on the orthodox Signed Quadrilateral (oSQ), which corresponds to the smallest even order $m$ (i.e. $m=4$): 
\begin{equation} \label{oSQ}
    Q_n = \sum_{i_1, i_2, i_3, i_4 \text{ (distinct)}}  (A_{i_1 i_2} - \hat{\alpha}_n) (A_{i_2i_3} - \hat{\alpha}_n) (A_{i_3 i_4} - \hat{\alpha}_n) (A_{i_4i_1} - \hat{\alpha}_n). 
\end{equation}
Again, since $A-\hat{\alpha}_n{\bf 1}_n{\bf 1}_n'$ is a stochastic proxy of $\widetilde{\Omega}$, the right hand side above is approximately equal to $\sum_{i_1, i_2, i_3, i_4 \text{ (distinct)}} \widetilde{\Omega}_{i_1i_2}\widetilde{\Omega}_{i_2i_3}\widetilde{\Omega}_{i_3i_4}\widetilde{\Omega}_{i_4i_1}\approx \tr(\widetilde{\Omega}^4)$. 
In other words, $Q_n$ is an estimate of $\tr(\widetilde{\Omega}^4)$, under the alternative hypothesis. 

These statistics are reminiscent of the classical moment statistics, as $X_n$ and $Q_n$ estimate an order-$2$ polynomial and an order-$4$ polynomial of (the entries of) $\widetilde{\Omega}$, respectively. 
We now recall some conventional insights of classical moment statistics. Suppose that we observe independent data $X_i\sim {\cal N}(\mu_i, \sigma^2)$, for $1\leq i\leq n$, and we would like to test the global null hypothesis 
\[
H_0: \;\; \mu_1=\mu_2=\cdots=\mu_n=0.
\]
Consider two moment statistics $S_1=\frac{1}{n}\sum_{i=1}^n X_i$ and $S_2=\frac{1}{n}\sum_{i=1}^n X^2_i$. It is well-known that, on the one hand, if the $\mu_i$'s have the same sign, the lower-order moment statistic $S_1$ has a better detection boundary than the higher-order moment statistic $S_2$; on the other hand, if the $\mu_i$'s have different signs, $S_1$ faces `signal cancellation', but $S_2$ has no such issue. Hence, if one only cares about the worst-case performance, using $S_2$ is enough. However, going beyond the `worst case', there are many cases where the power of $S_2$ is inferior to $S_1$,  so using $S_1$ to enhance power will be useful. 

In the network global testing we consider here, $X_n$ is analogous to $S_1$, and $Q_n$ is analogous to $S_2$. In Section~\ref{subsec:whyPE}, we will study the signal-to-noise ratios (SNRs) of $X_n$ and $Q_n$. We have the following observations: (i) The SNR of $X_n$ can be zero even when $\widetilde{\Omega}$ is a nonzero matrix, so the $\chi^2$ test faces potential signal cancellation. (ii) The SNR of $Q_n$ is always nonzero as long as $\widetilde{\Omega}$ is a nonzero matrix, but there exist cases where the SNR of $Q_n$ is strictly smaller than the SNR of $X_n$. 
Therefore, if we want a test that performs uniformly well in all cases, we should combine these two statistics to simultaneously avoid signal cancellation and enhance power.

In order for the combined test statistic to have a tractable null distribution, we must derive the asymptotic joint distribution of $X_n$ and $Q_n$. 
We show in Theorem~\ref{thm:jointNull} that with mild regularity conditions, 
\[
\begin{bmatrix}
(2n)^{-1/2}(X_n-n)\\
\bigl(2\sqrt{2}n^2\hat{\alpha}^2_n\bigr)^{-1}Q_n
\end{bmatrix}
\qquad \xrightarrow[n\to \infty]{\mathcal{L}} \qquad {\cal N}(0,I_2), \qquad \mbox{under }H_0^{(n)}. 
\]
Hence, a convenient way to combine both test statistics is to construct
\begin{equation} \label{PET}
    S_n = (2n)^{-1}(X_n-n)+(8n^4\hat{\alpha}_n^4)^{-1}Q_n^2, 
\end{equation}
which asymptotically follows the $\chi_2^2(0)$ distribution under the null hypothesis. We call $S_n$ the \textit{Power Enhancement (PE)} statistic.
We will show that the PE test combines the strengths of the $\chi^2$ test and the oSQ test while overcoming their respective limitations.

\subsection{Comparison of the signal-to-noise ratios}  \label{subsec:whyPE}

Let $\mathbb{E}_0$ and $\V_0$ denote the expectation and variance operators under the null hypothesis. Similarly, we have the notations $\mathbb{E}_1$ and $\V_1$ for the alternative hypothesis. 
The signal-to-noise ratio (SNR) of a test statistic $T_n$ is defined as 
\[
\mathrm{SNR}(T_n)\equiv \frac{|\mathbb{E}_1[T_n]-\mathbb{E}_0[T_n]|}{\sqrt{\max\{\V_0(T_n),\V_1(T_n)\}}}. 
\]
Under the alternative hypothesis, given $(\Pi, P)$, let $h = n^{-1} \sum_{i = 1}^n \pi_i$ and $\alpha_0=h'Ph$. 
In Theorem~\ref{thm:chi2-alt}, we show that the SNR of the $\chi^2$ statistic $X_n$ is captured by 
\begin{equation} \label{Definedeltan} 
\delta_n(K, P, h) = n^{3/2}\alpha_0^{-1}\|P h-\alpha_0{\bf 1}_K\|^2.  
\end{equation} 
Note that for any node $i$, the difference of expected degree under two hypothesis is $\mathbb{E}_1[d_i]-\mathbb{E}_0[d_i] =\sum_{j \neq i} (\pi_i'P \pi_j - \alpha_0)\approx n\pi_i'(Ph-\alpha_0{\bf 1}_K)$. Hence, the $\chi^2$ test finds evidence to reject the null hypothesis from node degrees. The $\chi^2$ test can successfully separate two hypotheses as long as $\delta_n\to\infty$. 
The oSQ test looks for evidence against the null hypothesis from $\tr(\widetilde{\Omega}^4)$. 
In Theorem~\ref{thm:SQ-alt}, we show that the SNR of the oSQ statistic $Q_n$ is captured by 
\begin{equation} \label{Definetaun} 
\tau_n(K, P, h) = n^2\alpha_0^{-2}\|P - \alpha_0{\bf 1}_K{\bf 1}_K' \|^4. 
\end{equation}
It can successfully separate two hypotheses as long as $\tau_n\to\infty$. 
The SNR of the PE statistic is captured by the quantity
\beq \label{Definebetan-add} 
\sqrt{\delta_n^2+\tau_n^2}\quad \asymp\quad  \beta_n(K, P, h)\equiv \max\{\delta_n, \tau_n\}.
\eeq
The PE test can successfully separate two hypotheses as long as $\beta_n\to\infty$.

The SNR of PE improves those of $\chi^2$ and oSQ. Such improvement can be significant. Consider the case of $Ph\propto {\bf 1}_K$ (e.g., when $P$ has equal diagonals and equal off-diagonals, and the communities have equal size), $\delta_n=0$; also, it can be shown that $\tau_n$ is always nonzero.  In this case, the $\chi^2$ test faces signal cancellation and loses power, but the PE test still has power. At the same time, when the community signals are very weak (e.g., $P$ is only a tiny perturbation of $\alpha_0{\bf 1}_K{\bf 1}_K'$), there exist cases where $\tau_n\to 0$ but $\delta_n\to\infty$. Then, the oSQ test is unsatisfactory, but the PE test is still satisfactory. 
It is illuminating to understand these cases from examples.  

\medskip
\noindent
{\bf Example 1} ({\it Two-community model}). 
Fix $K = 2$. For $a_n,b_n,d_n>0$ and $\epsilon_n\in (0,1)$, 
\[
P = \left[
\begin{array}{cccc} 
a_n  &&& b_n  \\ 
b_n &&& d_n  \\ 
\end{array}  
\right], \qquad\mbox{and}\qquad h = \begin{bmatrix} \epsilon_n \\ 1-\epsilon_n\end{bmatrix}. 
\] 
Write $\bar{a}_n=(a_n+d_n)/2$. Suppose that $|a_n-b_n|=O(a_n+b_n)$. By calculations in Appendix~\ref{subsec:Ex1-proof} of the supplementary material, we obtain the order of $\delta_n$ and $\tau_n$ in several cases (``S'' stands for ``symmetric'', and ``AS'' stands for ``asymmetric''): 
\begin{table}[htb]
\begin{tabular}{c|ccccc}
Case & Symmetry in $h$ & Symmetry in $P$ & SNR of $X_n$ & SNR of $Q_n$\\
\midrule
S  & $\epsilon_n=1/2$ & $d_n=a_n$ & $0$ & $n^2 \bigl[ \frac{(a_n-b_n)^2}{a_n+b_n}\bigr]^2$\\
AS1 &  $\epsilon_n\neq 1/2$ & $d_n=a_n$  & $(1-2\epsilon_n)^2n^{3/2}\frac{(a_n-b_n)^2}{a_n+b_n}$ & $n^2 \bigl[ \frac{(a_n-b_n)^2}{a_n+b_n}\bigr]^2$\\
AS2 & $\epsilon_n=1/2$ & $|d_n-a_n|\gg |\bar{a}_n-b_n|$ & $n^{3/2}\frac{(d_n-a_n)^2}{\bar{a}_n+b_n}$ & $n^2 \bigl[ \frac{(d_n-a_n)^2}{\bar{a}_n+b_n}\bigr]^2$\\
AS3 & $\epsilon_n=1/2$ & $|d_n-a_n|\gg |\bar{a}_n-b_n|$ & $n^{3/2} \frac{(d_n-a_n)^2}{\bar{a}_n+b_n}$ & $n^2 \bigl[ \frac{(\bar{a}_n-b_n)^2}{\bar{a}_n+b_n}\bigr]^2$
\end{tabular}
\end{table}

In Case (S) (it includes the 2-community SBM in \citep{mossel2015reconstruction,banerjee2017optimal} as a special case), the $\chi^2$ test loses power due to the fact that $\delta_n=0$, and the oSQ test and the PE test have full power provided that 
\[
(a_n-b_n)^2/(a_n+b_n)\;\; \gg\;\;  n^{-1}.
\] 
In Case (AS1), suppose $|1-2\epsilon_n|\geq c$ for a constant $c>0$. Then, if 
\[
n^{-3/2}\;\; \ll\;\; (a_n-b_n)^2/(a_n+b_n)\;\; \ll \;\; n^{-1},
\]
the oSQ test does not have full power but the $\chi^2$ test and  the PE test have full power.

\medskip
\noindent
{\bf Example 2} ({\it Rank-1 model}). Fix $K>1$ and consider an MMSBM with $P =   \eta \eta'$, for some nonnegative vector $\eta\in\mathbb{R}^K$ such that $\|\eta\|_\infty\leq 1$ and $\eta\not\propto {\bf 1}_K$. This is an example where $K$ is not identifiable. In Section~\ref{sec:identifiability}, we will define the \textit{intrinsic number of communities (INC)}, which is the smallest $K_0$ such that this model can be written as a $K_0$-community MMSBM. By calculations in Appendix~\ref{subsec:Ex2-proof}, the INC for this example is $2$ (regardless of $K$), so the alternative hypothesis holds. Consider a special case of $K=2$:
\[
P = \eta\eta', \qquad \eta = \frac{\sqrt{c_n}}{\sqrt{a_n^2+b_n^2}} \begin{bmatrix} a_n\\b_n \end{bmatrix}, \qquad h=\begin{bmatrix}1/2\\1/2\end{bmatrix}, 
\]
where $a_n,b_n>0$ and $c_n\in (0,1)$. 
By direct calculations in Appendix~\ref{subsec:Ex2-proof}, when $|a_n-b_n|=O(a_n+b_n)$, 
\[
\delta_n\;\; \asymp\;\; n^{3/2}c_n\frac{(a_n-b_n)^2}{(a_n^2+b_n^2)}, \qquad \tau_n\;\; \asymp\;\; n^2c_n^2 \frac{(a_n-b_n)^4}{(a_n^2+b_n^2)^2}\;\;\asymp \;\; n^{-1}\delta_n^2. 
\]
It is seen that $\tau_n\to\infty$ implies $\delta_n\to\infty$. Hence, whenever the oSQ test has full power, the $\chi^2$ test also has full power, so does the PE test. However, when $\delta_n\to\infty$ but $\delta_n\ll \sqrt{n}$, the $\chi^2$ test and the PE test both have full power, but the oSQ test does not have full power. 

\medskip

From these examples, the oSQ test outperforms the $\chi^2$ test sometimes (e.g., in Case (S) of Example 1), and the $\chi^2$ test outperforms the oSQ test sometimes (e.g., in Example 2). Power enhancement allows us to combine the strengths of both tests. Below, we further show that the PE test achieves the optimal phase transition.

\subsection{A preview of the phase transition} \label{subsec:intro-phase}

\begin{figure}[tb!] 
\centering
\includegraphics[height=2.5in]{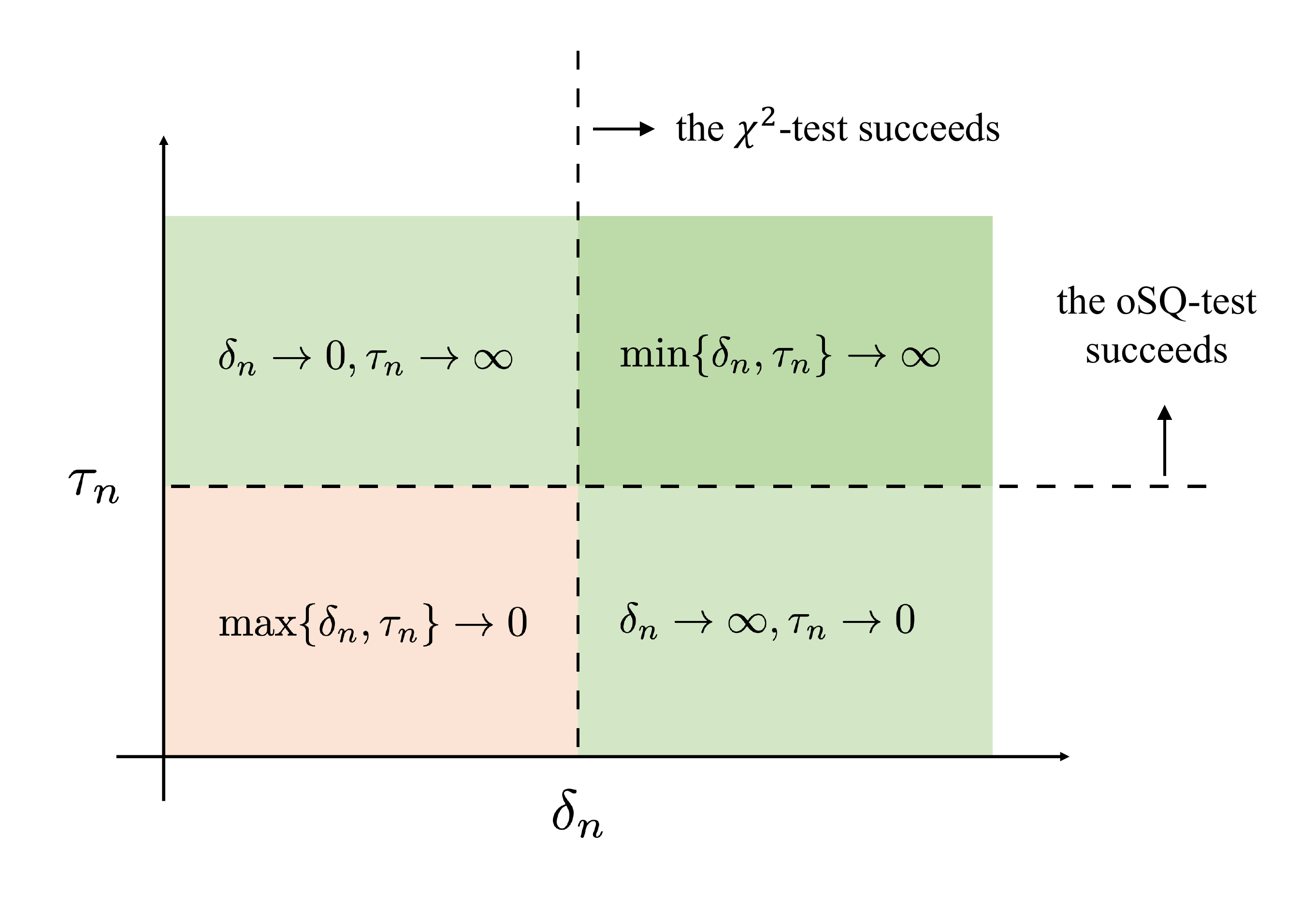}  
\caption{The phase transition of global detection for MMSBM. The light red region is the Region of Impossibility, and the three other green regions constitute the region of Possibility.}
\label{fig:phase} 
\end{figure}

The PE test can successfully separate two hypotheses if $\max\{\delta_n,\tau_n\}\to \infty$. In Theorem~\ref{thm:LB}, we provide a matching lower bound: If $\max\{\delta_n,\tau_n\}\to 0$, 
\[
\mbox{the chi-square distance between the probability densities of two hypotheses 
$\goto 0$}. 
\]
Hence, there exists no test that can asymptotically distinguish two hypotheses when $\max\{\delta_n, \tau_n\}\to0$. 
It gives rise to the following phase transition: 
Consider the two dimensional phase space for MMSBM, where the $x$-axis is $\delta_n$ which calibrates the signals in node degrees, and the $y$-axis is $\tau_n$ which calibrates the signals in cycle counts. 
The phase space is divided into two regions (see Figure~\ref{fig:phase}):  
\begin{itemize} \itemsep +5pt
\item {\it Region of Impossibility ($\max\{\delta_n, \tau_n\} \goto 0$)}.  Any alternative in this region is inseparable from a null. For any test, the sum of Type I and Type II errors tends to $1$ as $n\goto\infty$. 

\item {\it Region of Possibility ($\max\{\delta_n, \tau_n\} \goto \infty$)}.  Any alternative hypothesis in this region is separable from any null hypothesis. Specifically, the PE test statistic is able to separate two hypotheses, in the sense that for an appropriate threshold, the sum of Type I and Type II errors of the PE test tends to $0$ as $n \goto \infty$. 

\end{itemize} 
We say that a test is {\it optimally adaptive} if it can distinguish the null and alternative hypotheses in the {\it whole} Region of Possibility. The PE test is optimally adaptive, but neither the $\chi^2$ test nor the oSQ test is. The $\chi^2$ test can only distinguish two hypotheses in the sub-region of $\delta_n \goto \infty$, and the oSQ test can only distinguish two hypotheses in the sub-region of $\tau_n\goto\infty$. See Figure~\ref{fig:phase}.

\medskip
\noindent
{\bf Remark 1} {\it (Connection to minimax optimality)}: Our phase transition results are more informative than the standard {\it minimax} results. The detection boundary $\beta_n(K, P, h)$ in \eqref{Definebetan-add} is for arbitrary $(K, P, h)$, while the minimax detection boundary is the specific value of $\beta_n(K,P,h)$ at some worst-case $(K, P, h)$ in a pre-specified class. Therefore, for a test to be {\it minimax optimal}, it only requires that the SNR of this test  matches $\beta_n(K,P,h)$ at the worst-case $(K,P, h)$; however, for a test to be {\it optimally adaptive}, its SNR has to match $\beta_n(K,P,h)$ for all $(K,P,h)$. We discuss this more carefully in Section~\ref{subsec:LB}.  

\smallskip
 \noindent
{\bf Remark 2} {\it (Other moment statistics)}: In \eqref{oSP}, we have defined the {\it order-$m$ orthodox Signed Polygon statistic} $U_n^{(m)}$, for $m\geq 3$. We now define the {\it length-$m$ Signed Path} statistic $V_n^{(m)}$ by
\beq \label{Signed-Path}
V_n^{(m)}=\sum_{i_1,\ldots,i_{m+1} \text{ (distinct)}} (A_{i_1i_2}-\hat{\alpha}_n)(A_{i_2i_3}-\hat{\alpha}_n)\cdots (A_{i_mi_{m+1}}-\hat{\alpha}_n), \qquad \mbox{for }m\geq 2. 
\eeq
The two statistics distinguish two hypotheses from estimating $\tr(\widetilde{\Omega}^m)$ and ${\bf 1}_n'\widetilde{\Omega}^m{\bf 1}_n$, respectively. The oSQ statistic is $U_n^{(m)}$ with $m=4$, and the $\chi^2$ statistic is equivalent to $V_n^{(m)}$ with $m=2$ (see \eqref{chi1-equivalent}). A natural question is whether we should consider other values of $m$. In Appendix~\ref{subsec:Remark2-proof}, we show that under some regularity conditions: 
\begin{align*}
& \mathrm{SNR}(U_n^{(m)})\;\;  \asymp\;\;  n^m\alpha_0^{-m}\|P-\alpha_0{\bf 1}_K{\bf 1}_K'\|^m \;\; \asymp \;\; \tau_n^{m/4},\cr
&\mathrm{SNR}(V_n^{(m)}) \;\; \asymp\;\;  n^{\frac{m+1}{2}} \alpha_0^{-\frac{m}{2}}\|P-\alpha_0{\bf 1}_K{\bf 1}_K'\|^{m-2}\|Ph-\alpha_0{\bf 1}_K\|^2 \;\; \asymp\;\; \tau_n^{(m-2)/4}\delta_n. 
\end{align*}
Hence, either $\mathrm{SNR}(U_n^{(m)})\to\infty$ or $\mathrm{SNR}(V_n^{(m)})\to\infty$ is a stronger requirement than $\max\{\delta_n,\tau_n\}\to\infty$. In other words, we do not benefit from a better phase transition by considering other values of $m$.

\smallskip
 \noindent
{\bf Remark 3} {\it (The regime of a constant SNR)}: Our phase transition covers the case of SNR $\to 0$ and SNR $\to \infty$. It is also interesting to study the regime of a constant SNR. In the special case of symmetric 2-community SBM (see Example 1), this regime has been well understood. If $n(a_n-b_n)^2/(a_n+b_n)<1$, the two hypotheses are mutually contiguous; if $n(a_n-b_n)^2/(a_n+b_n)>1$, the two hypotheses are asymptotically singular, and the signed cycle statistic $U_n^{(m)}$, with $m\asymp \log^{1/4}(n)$, has asymptotically full power \cite{mossel2015reconstruction,banerjee2017optimal,gao2021minimax}. Beyond this special case, much less is known. 
Some partial results were obtained for general structure of $P$ \cite{banks2016information}, unequal-size communities \cite{zhang2016community}, and mixed memberships 
\cite{hopkins2017efficient}, where they were primarily interested in getting a good estimate of $\Pi$ rather than the global testing problem. 
It is largely unclear whether there exists a test with asymptotically full power in the constant SNR regime for a general MMSBM. Given the expression of $\beta_n$ in our phase transition and the full-power test for the special 2-community SBM, we conjecture that a full-power test may be constructed from a weighted combination of $\{U_n^{(m)}, V_n^{(m)}\}_{1\leq m\leq \log(n)}$, where $U_n^{(m)}$ and $V_n^{(m)}$ are the signed cycle statistics and signed path statistics as in Remark 2. 
We note that our proposed PE test is a combination of $U_n^{(4)}$ and $V_n^{(2)}$.

\smallskip
 \noindent
{\bf Remark 4} {\it (Which test to use at finite $n$ and with knowledge of parameters)}: 
Our theoretical results are for the asymptotic setting of $n\to\infty$, but our analysis does give the SNRs of these statistics for finite $n$. Recall that $\Omega=\Pi P\Pi'$ and $\widetilde{\Omega}=\Pi (P-\alpha_0{\bf 1}_K{\bf 1}_K')\Pi'$. Let $[\Omega\circ (1-\Omega)]$ be the matrix whose $(i,j)$th entry is $\Omega_{ij}(1-\Omega_{ij})$. We can obtain that
\[
\mathrm{SNR}(X_n)\; \sim\; \frac{{\bf 1}_n'\widetilde{\Omega}^2{\bf 1}_n}{\sqrt{2({\bf 1}_n'[\Omega\circ(1-\Omega)]^2 {\bf 1}_n)}}, \qquad \mathrm{SNR}(Q_n)\;\sim\; \frac{\mathrm{tr}(\widetilde{\Omega}^4)}{\sqrt{8\mathrm{tr}([\Omega\circ(1-\Omega)]^4)}}.
\]
In principle, for finite $n$, if the parameters of the alternative hypothesis are given, we can compute these precise SNRs and decide which test to use. 
However, in practice, we always recommend using the PE test. One advantage of PE is its `adaptivity': It yields a good power uniformly in many settings, so we do not worry about choosing between $\chi^2$ and oSQ.

\smallskip
 \noindent
{\bf Remark 5} {\it (Comparison with DCMM)}: 
The DCMM model \citep{Mixed-SCORE,SP2021} generalizes MMSBM by accommodating degree heterogeneity. It introduces a degree parameter $\theta_i>0$ for each node $i$ and assumes 
\[
\Omega_{ij}\equiv \mathbb{P}(A_{ij}=1) =\theta_i\theta_j (\pi_i'P\pi_j), \qquad 1\leq i<j\leq n.
\]
Although MMSBM is a sub-class of DCBM by forcing $\theta_i\equiv 1$, the global testing for these two models is quite different.  Consider Example 2 in Section~\ref{subsec:whyPE}, where $P=\eta\eta'$. By letting $\tilde{\theta}_i=\theta_i\cdot \pi_i'\eta$, we can also write $\Omega_{ij}=\tilde{\theta}_i\tilde{\theta}_j$ for all $1\leq i,j\leq n$. Then, it becomes a null model under DCBM, although it is still an alternative model under MMSBM (where the intrinsic number of communities is 2; see Section~\ref{sec:identifiability}). This example shows that restricting to a sub-class of models can change the detection boundary.

Compared with MMSBM, DCMM has many more free parameters, so `degree matching' \cite{SP2021} is possible: Given any alternative DCMM, there exists a null DCMM such that for each node, its expected degree under the null model is matched with its expected degree under the alternative model. Then, any degree-based test loses power. In contrast, such `degree matching' is impossible under MMSBM; and we can find many settings where the (degree-based) $\chi^2$ test has superior power. Hence, to achieve the optimal phase transition, it is crucial to use the $\chi^2$ statistic for `power enhancement'.

\section{Main results} \label{sec:main}
Recall that we consider the global testing problem \eqref{gt-0} under the MMSBM model, where the Bernoulli probability matrix $\Omega$ under two hypotheses are as in \eqref{gt-1A}. In Section~\ref{subsec:main-null}, we derive the null distributions of the three test statistics ($\chi^2$, oSQ and PE). In Section~\ref{subsec:main-alt}, we study the power of the three tests. We provide lower bound arguments and phase transitions in Section~\ref{subsec:LB}.

\subsection{The asymptotic null distributions} \label{subsec:main-null}
Under the null hypothesis, $\Omega=\alpha_n{\bf 1}_n{\bf 1}'_n$, with $\alpha_n\in (0,1)$ calibrating the sparsity level of the network. We estimate $\alpha_n$ by $
\hat{\alpha}_n= (n-1)^{-1} \bar{d}$, where $\bar{d}$ is the average node degree. Let $X_n$ and $Q_n$ be the $\chi^2$ statistic and the oSQ statistic in \eqref{chi1} and \eqref{oSQ}, respectively. The following theorem characterizes the joint null distribution of $(X_n, Q_n)$, which is proved in the supplementary material.

\begin{theorem}[Asymptotic joint null distribution] \label{thm:jointNull} 
Consider the global testing problem in \eqref{gt-0}, where $\Omega=\alpha_n{\bf 1}_n{\bf 1}_n'$ under the null hypothesis. Suppose that $\alpha_n\leq 1- c_0$ for a constant $c_0\in (0,1)$ and that $n\alpha_n\to\infty$ as $n\to\infty$. Then, under the null hypothesis, 
\[
\left( \frac{X_n-n}{\sqrt{2n}},\;  \frac{Q_n}{2\sqrt{2}n^2\hat{\alpha}_n^2} \right)\quad \xrightarrow[n\to \infty]{\mathcal{L}}\quad \mathcal{N}(0, I_2).  
\] 
\end{theorem} 

We immediately obtain the asymptotic distributions of the three test statistics. 

\begin{corollary} \label{cor:null}
Under the conditions of Theorem~\ref{thm:jointNull}, as $n\to\infty$, the following statements are true:
\begin{itemize}
\item The $\chi^2$ test statistic satisfies that $(X_n-n)/\sqrt{2n}\longrightarrow {\cal N}(0,1)$ in distribution.
\item The oSQ test statistic satisfies that $Q_n/(2\sqrt{2}n^2\hat{\alpha}_n^2)\longrightarrow {\cal N}(0,1)$ in distribution. 
\item The PE test statistic satisfies that $S_n\to \chi_2^2(0)$ in distribution. 
\end{itemize}
\end{corollary}

Fix any $\epsilon\in(0,1)$. The level-$\epsilon$ degee-based $\chi^2$ test rejects the null hypothesis if
\beq \label{rejectR-DC}
(X_n-n)/\sqrt{n}\quad >\quad \mbox{$(1-\epsilon)$-quantile of $ \mathcal{N}(0,1)$}. 
\eeq
The level-$\epsilon$ oSQ test rejects the null hypothesis if
\beq \label{rejectR-SQ}
Q_n/\bigl(2\sqrt{2}n^2\hat{\alpha}_n^2) \quad >\quad \mbox{$(1-\epsilon)$-quantile of $ \mathcal{N}(0,1)$}. 
\eeq
The level-$\epsilon$ PE test rejects the null hypothesis if
\beq \label{rejectR-PET}
S_n \quad >\quad \mbox{$(1-\epsilon)$-quantile of $\chi^2_2(0)$}. 
\eeq

\noindent
{\bf Remark 6}:  
We give a brief explanation of why $X_n$ and $Q_n$ are asymptotically uncorrelated. 
Let $Q_n^*$ be a proxy of $Q_n$ by replacing $\hat{\alpha}_n$ by $\alpha_n$ in \eqref{oSQ}. Moreover, from \eqref{chi1-equivalent}, we can re-write $X_n=n + \frac{1}{n\hat{\alpha}_n(1-\hat{\alpha}_n}\sum_{i,j,k (dist)}(A_{ij}-\hat{\alpha}_n)(A_{jk}-\hat{\alpha}_n)$. Replacing $\hat{\alpha}_n$ by $\alpha_n$ in this expression leads to a proxy of $X_n$, denoted by $X_n^*$.  Let $W=A-\mathbb{E}[A]$. 
Under the null hypothesis,  $A_{ij}=\alpha_n+W_{ij}$. It follows that
\[
X_n^*= n + \frac{1}{n\alpha_n(1-\alpha_n)} \sum_{i,j,k \text{ (dist)}}W_{ij}W_{jk}, \qquad Q^*_n=\sum_{i_1,i_2,i_3,i_4 \text{ (dist)}}W_{i_1i_2}W_{i_2i_3}W_{i_3i_4}W_{i_4i_1}. 
\]
A key observation is that $\mathbb{E}[(W_{ij}W_{jk})(W_{i_1i_2}W_{i_2i_3}W_{i_3i_4}W_{i_4i_1})]=0$ for all $(i, j, k, i_1, i_2, i_3, i_4)$ such that $(i,j,k)$ are distinct and $(i_1,i_2,i_3,i_4)$ are distinct. To verify this, it suffices to check all possible cases of $\{i,j,k\}\cap \{i_1,i_2,i_3,i_4\}\neq\emptyset$. For example, when $i=i_1$, $j=i_2$ and $k=i_3$, the expectation is equal to $\mathbb{E}[W^2_{ij}W^2_{jk}W_{ki_4}W_{i_4i}]=0$. Other cases are similar. Therefore, $X_n^*$ and $Q_n^*$ are uncorrelated. Since $X_n\approx X_n^*$ and $Q_n\approx Q_n^*$, we can show that $X_n$ and $Q_n$ are asymptotically uncorrelated.

\subsection{Power analysis} \label{subsec:main-alt}
Under the alternative hypothesis, $\Omega=\Pi P\Pi'$. We notice that the parameters are not identifiable. There may exist $K^*\neq K$ and $(\Pi^*, P^*)$ such that $\Omega=\Pi^*P^*(\Pi^*)'$ also holds. To address this issue, we follow Occam's razor to choose the parameters associated with the smallest possible $K$. This $K$ is called the Intrinsic Number of Communities (INC). \footnote{We will show in Sectioin~\ref{sec:identifiability} that $\mathrm{INC}=1$ if and only if $\Omega\propto {\bf 1}_n{\bf 1}_n'$, which is compatible with the null model in \eqref{gt-1A}. } In this subsection, we always assume that $(K, P, \Pi)$ are the parameters associated with INC. 
The detailed discussion of INC is deferred to Section~\ref{sec:identifiability}. 

Write
\[
h=\frac{1}{n}\sum_{i=1}^n \pi_i, \qquad G = \frac{1}{n}\sum_{i=1}^n \pi_i\pi_i', \qquad \alpha_0=\alpha_0(h,P)=h'Ph.
\]
We assume there exists a constant $C>0$ such that
\beq\label{cond1}
\frac{\max_{1\leq k\leq K}h_k}{\min_{1\leq k\leq K}h_k}\leq C, \qquad\mbox{and}\qquad \|G^{-1}\|\leq C.  
\eeq
For a constant $c\in (0,1)$, we assume
\beq \label{cond2}
\alpha_0\leq c, \qquad \mbox{and}\qquad n\alpha_0\geq c^{-1}. 
\eeq
These conditions are mild. Condition \eqref{cond1} is about balance of communities. This is easier to see in the case of no mixed membership (i.e., $\pi_i$ only takes  values in $\{e_1,\ldots,e_K\}$). In this case, $G$ is a diagonal matrix, and both $h_k$ and $G_{kk}$ are equal to the fraction of nodes in community $k$. Then, \eqref{cond1} says that the fraction of nodes in each community is bounded away from zero, which is a mild condition (e.g., \cite{zhang2016minimax} used a similar condition). Condition \eqref{cond2} is about network sparsity. Under our model, the average node degree is at the order of $n\alpha_0$; therefore, \eqref{cond2} allows the average node degree to range from $O(1)$ to $O(n)$, which covers a wide range of sparsity.

First, we study the (degree-based) $\chi^2$ test. 

\begin{theorem}[Power of the $\chi^2$ test] \label{thm:chi2-alt} 
Consider the global testing problem in \eqref{gt-0}, where \eqref{cond1}-\eqref{cond2} hold under the alternative hypothesis. Let 
\[
\delta_n=n^{3/2}\alpha_0^{-1}\|Ph-\alpha_0{\bf 1}_K\|^2. 
\]
Suppose $\delta_n\geq C$. There exist a constant $c_1>0$ such that, under the alternative hypothesis,
\begin{align*}
\mathbb{E}\biggl[\frac{X_n-n}{\sqrt{2n}}\biggr] & \geq c_1\delta_n- O\bigl(n^{-1/2}\alpha_0^{-1/2}\bigr), \cr
\mathrm{Var}\biggl( \frac{X_n-n}{\sqrt{2n}}\biggr) & =  
O\Bigl(1+n^{-1/2}\delta_n+ n^{-2}\alpha_0^{-1}\delta_n^2\log(n) \Bigr). 
\end{align*}
\end{theorem}

We are interested in the scenario of $\delta_n\to\infty$. Write $\psi^{(1)}_n=(X_n-n)/\sqrt{2n}$ for short. 
By Theorem~\ref{thm:chi2-alt} and the assumption $n\alpha_n\geq c^{-1}$, we have:
\[
\mathbb{E}[\psi^{(1)}_n]\geq C\delta_n, \qquad \mathrm{SD}(\psi^{(1)}_n)\leq \begin{cases}
C\bigl(1+\delta_n\sqrt{\log(n)/n}\bigr), & \mbox{if }\delta_n=O(\sqrt{n}),\\
C\bigl(n^{-1/4}\sqrt{\delta_n}+\delta_n\sqrt{\log(n)/n}\bigr), &\mbox{if }\delta_n\gg\sqrt{n}. 
\end{cases}
\]
It follows that the SNR is
\[
\geq 
\begin{cases}
C\min\bigl\{\delta_n, \, \sqrt{n/\log(n)}\bigr\}, & \mbox{if }\delta_n=O(\sqrt{n}),\\
C\min\bigl\{n^{1/4}\sqrt{\delta_n}, \, \sqrt{n/\log(n)}\bigr\}, &\mbox{if }\delta_n\gg\sqrt{n}. 
\end{cases}
\]
In either case, the SNR tends to $\infty$. We expect that the $\chi^2$ test successfully separates the null from the alternative. This gives the following corollary: 

\begin{corollary} \label{cor:chi2}
Consider a pair of hypotheses as in \eqref{gt-0}-\eqref{gt-1A}, where $\alpha_n\leq 1-c_0$ and $n\alpha_n\to\infty$ under the null hypothesis and the conditions \eqref{cond1}-\eqref{cond2} hold under the alternative hypothesis. Suppose $\delta_n\to\infty$ under the alternative hypothesis. 
For any fixed $\epsilon\in (0,1)$, consider the level-$\epsilon$ $\chi^2$ test in \eqref{rejectR-DC}. Then, the level of the test tends to $\epsilon$ and the power of the test tends to  $1$.   
\end{corollary}

Next, we study the oSQ test. 

\begin{theorem}[Power of the oSQ test] \label{thm:SQ-alt} 
Consider the global testing problem in \eqref{gt-0}, where \eqref{cond1}-\eqref{cond2} hold under the alternative hypothesis. Let
\[
\tau_n = n^2\alpha_0^{-2}\|P-\alpha_0{\bf 1}_K{\bf 1}_K'\|^4. 
\]
Suppose $\tau_n\geq C$. 
There exists a constant  $c_2>0$ such that, under the alternative hypothesis, 
\begin{align*}
\mathbb{E}\biggl[\frac{Q_n}{2\sqrt{2}n^2 \hat{\alpha}_n}\biggr] & \geq  c_2\tau_n-o\left(1+n^{-1/2}\tau_n^{3/4}\right), \cr
\mathrm{Var}\biggl(\frac{Q_n}{2\sqrt{2}n^2\hat{\alpha}_n^2}\biggr) &=O\Bigl(1+ n^{-1}\tau_n^{3/2}+n^{-2}\alpha_0^{-1}\tau_n^2\log(n)\Bigr).
\end{align*}
\end{theorem}

We are interested in the cases where $\tau_n\to\infty$. Write $\psi^{(2)}_n=Q_n/\bigl(2\sqrt{2}n^2 \hat{\alpha}_n\bigr)$ for short. 
By Theorem~\ref{thm:SQ-alt} and the assumption $n\alpha_n\geq c^{-1}$, we have:
\[
\mathbb{E}[\psi^{(2)}_n]\geq C\tau_n, \qquad \mathrm{SD}(\psi^{(2)}_n)\leq \begin{cases}
C\bigl(1+\tau_n\sqrt{\log(n)/n}\bigr), & \mbox{if }\tau_n=O(n^{2/3}),\\
C\bigl(n^{-1/2} \tau_n^{3/4}+\tau_n\sqrt{\log(n)/n}\bigr), &\mbox{if }\tau_n\gg n^{2/3}. 
\end{cases}
\]
It follows that the SNR is
\[
\geq 
\begin{cases}
C\min\bigl\{\tau_n,\,\sqrt{n/\log(n)}\bigr\}, & \mbox{if }\tau_n=O(n^{2/3}),\\
C\min\bigl\{n^{1/2}\tau_n^{1/4},\, \sqrt{n/\log(n)}\bigr\}, &\mbox{if }\tau_n\gg n^{2/3}. 
\end{cases}
\]
In either case, the SNR tends to $\infty$. We expect that the oSQ test can successfully separate two hypotheses, as stated in the following corollary:

\begin{corollary} \label{cor:SQ}
Consider a pair of hypotheses as in \eqref{gt-0}-\eqref{gt-1A}, where $\alpha_n\to 0$ and $n\alpha_n\to\infty$ under the null hypothesis and \eqref{cond1}-\eqref{cond2} hold under the alternative hypothesis. Suppose $\tau_n\to\infty$ under the alternative hypothesis. 
For any fixed $\epsilon\in (0,1)$, consider the level-$\epsilon$ oSQ test in \eqref{rejectR-SQ}. Then, as $n\to\infty$, the level of the test tends to $\epsilon$ and the power of the test tends to $1$.   
\end{corollary}

Last, we study the PE test,  where the test statistic is $S_n=\left(\psi_n^{(1)}\right)^2+\left(\psi_n^{(2)}\right)^2$. 

\begin{theorem}[Power of the PE test] \label{thm:PET-alt} 
Consider a pair of hypotheses as in \eqref{gt-0}-\eqref{gt-1A}, where $\alpha_n\to 0$ and $n\alpha_n\to\infty$ under the null hypothesis and \eqref{cond1}-\eqref{cond2} hold under the alternative hypothesis. Suppose $\max\{\delta_n, \tau_n\}\to\infty$ under the alternative hypothesis. Then, under the alternative hypothesis, 
\[
S_n \xrightarrow{\mathbb{P}} \infty. 
\]
Furthermore, for any fixed $\epsilon\in (0,1)$, consider the level-$\epsilon$ PE test in \eqref{rejectR-PET}. Then, as $n\to\infty$, the level of the test tends to $\epsilon$ and the power of the test tends to $1$.   
\end{theorem}

By Theorem~\ref{thm:PET-alt}, the PE test successfully distinguishes two hypotheses as long as $\delta_n^2+\tau_n^2\to\infty$ (or equivalently, $\max\{\delta_n, \tau_n\}\to\infty$).

\subsection{The lower bounds and phase transitions} \label{subsec:LB}

To obtain lower bounds, we switch to the random-membership MMSBM (this follows the convention: If we use a non-random $\Pi$, only trivial lower bounds can be obtained, which is uninteresting). Fix $K\geq 2$ and consider a sequence of $(P_n, F_n)$, indexed by $n$, where $P_n \in\mathbb{R}^{K\times K}$ is an eligible community matrix and $F_n$ is a distribution on the probability simplex of $\mathbb{R}^K$. Let $h_n=\mathbb{E}_{\pi\sim F_n}[\pi]\in\mathbb{R}^{K}$. We often drop the subscript $n$ in $(P_n, F_n, h_n)$ to simplify the notations. The (randomized) alternative hypothesis is
\beq \label{RMM-alt}
H_1^{(n)}:\qquad\Omega = \Pi P\Pi', \qquad\mbox{where}\quad \Pi =[\pi_1,\pi_2,\ldots,\pi_n]', \quad\pi_i\overset{iid}{\sim}F. 
\eeq
We pair this alternative hypothesis with the null hypothesis below:
\beq \label{RMM-null}
H_1^{(n)}:\qquad\Omega=\alpha_0{\bf 1}_n{\bf 1}_n', \qquad \mbox{where}\quad \alpha_0=h'Ph. 
\eeq
Let $f_1(A)$ and $f_0(A)$ be the probability densities associated with \eqref{RMM-alt} and \eqref{RMM-null}, respectively. The $\chi^2$-distance between two hypotheses is defined as $\int [\frac{f_1(A)}{f_0(A)}-1]^2f_0(A)dA$. Two hypotheses are asymptotically indistinguishable if the $\chi^2$-distance $\to 0$. 
The following theorem is proved in the supplementary material.

\begin{theorem}[Lower bound] \label{thm:LB}    
Consider a sequence of hypothesis pairs \eqref{RMM-alt}-\eqref{RMM-null} indexed by $n$. Let
\[
\beta_n=\beta_n(K, P_n,h_n)=\max\bigl\{ n^{3/2}\alpha_0^{-1}\|P h-\alpha_{0}{\bf 1}_K\|^2,\quad  n^2\alpha_0^{-2}\|P - \alpha_{0}{\bf 1}_K{\bf 1}_K' \|^4  \bigr\}. 
\]
If $\beta_n\to 0$, then the chi-square distance between two hypotheses converges to $0$ as $n\to\infty$.
\end{theorem}

We now combine Theorem~\ref{thm:LB} with Theorem~\ref{thm:PET-alt} and obtain the phase transitions:
\begin{itemize}
\item {\it Region of Impossibility}. When $\beta_n=\max\{\delta_n,\tau_n\}\to 0$, by Theorem~\ref{thm:LB}, the two hypotheses are asymptotically inseparable, where for any test, the sum of type I and type II errors tends to $1$ as $n\goto\infty$. 
\item {\it Region of Possibility}. When $\beta_n=\max\{\delta_n, \tau_n\}\to\infty$, by Theorem~\ref{thm:PET-alt}, the PE test can successfully separate the two hypotheses: for properly chosen $\epsilon_n\to 0$,  the sum of type I and type II errors of the level-$\epsilon_n$ PE test tends to $0$ as $n\goto\infty$.
\end{itemize}
We conclude that the PE test is optimally adaptive. As we have mentioned in Section~\ref{sec:intro}, none of the previously existing tests are optimally adaptive.

Theorem~\ref{thm:PET-alt} is more informative than the standard minimax lower bound. To prove a minimax lower bound, we only need to pick one `worst-case' configuration of $(K, P_n, h_n)$, but Theorem~\ref{thm:PET-alt} is for all configurations of $(K, P_n, h_n)$.  
The test that works well for the `worst-case' configuration may have unsatisfactory performances for other configurations. 
For example, a commonly studied configuration in the literature is $P_n=(a_n-b_n)I_K+b_n{\bf 1}_K{\bf 1}_K'$ and $h_n=(1/K){\bf 1}_K$. For this configuration, $\delta_n=0$, and the oSQ test alone is optimal.  
However, when $(P_n, h_n)$ deviate from this configuration, the PE test can outperform the oSQ test (as seen in the simulations in Section~\ref{sec:simul}).  

Theorem~\ref{thm:LB} can be used to derive the minimax lower bounds for different parameter classes. In standard minimax arguments, we adopt the original MMSBM with a non-random $\Pi$, but we will consider the worst case performance for a class of parameters. Let $\{\alpha_n\}_n$ be a positive sequence in $[0,1]$ and $\{\gamma_n\}_n$ be another positive sequence. Given $(K, P, \Pi)$, write $h=h(\Pi)=n^{-1}\sum_{i=1}^n\pi_i$ and $\alpha_0=\alpha_0(P, \Pi)=h'Ph$. We introduce the following classes of Bernoulli probability matrices:
\begin{align*}
{\cal M}_{0n}(c_0, \alpha_n) &= \bigl\{ \Omega=b\bm{1}_n\bm{1}_n': \; b\leq 1-c_0, \; b \geq \alpha_n \bigr\},\cr
{\cal M}_{1n}(t_0;  K, c, C,\alpha_n, \gamma_n)&= \left\{\begin{array}{ll}
\Omega=\Pi P\Pi': & \mbox{ $(\Pi, P)$ satisfies \eqref{cond1}-\eqref{cond2} for $c, C>0$,}\\
&\alpha_0 \equiv \alpha_0(P,\Pi) \geq \alpha_n/2, \\
& \|Ph(\Pi)-\alpha_0 {\bf 1}_K \|_F\geq t_0\|P-\alpha_0{\bf 1}_K{\bf 1}_K'\|, \\
&\mbox{and } 2\alpha_0^{-1}\|P-\alpha_0{\bf 1}_K{\bf 1}_K'\|\geq \gamma_n
\end{array}
\right\}.
\end{align*}
We abbreviate the two classes as $M_{0n}$ and $M_{1n}(t_0)$, respectively. We note that when $t_0>0$, an additional constraint on $(P, \Pi)$ is imposed. 
Define the minimax testing risk as (below, the infimum is taken over all possible tests $\psi$)
\beq \label{MinimaxRisk}
Risk_n^*(t_0)= \inf_{\psi}\biggl\{ \sup_{\Omega\in {\cal M}_{0n}}\mathbb{P}(\psi=1) + \sup_{\Omega \in {\cal M}_{1n}(t_0)} \mathbb{P}(\psi=0)\biggr\}. 
\eeq
The following theorem is proved in the supplementary material:

\begin{theorem}[Minimax lower bound] \label{thm:minimax}
Fix $K\geq 2$. Suppose $\alpha_n\goto 0$, $\gamma_n\goto 0$ and $n\alpha_n\goto\infty$. 
\begin{itemize}
\item Fix $t_0=0$. If $n^2\alpha_n^2\gamma_n^4\to 0$, then $\lim_{n\to\infty}\{ Risk_n^*(t_0)\} = 1$. 
\item Fix $0<t_0<\sqrt{\frac{(K-1)(K+3)}{16K}}$. If $n^{3/2}\alpha_n\gamma_n^2\to 0$, then $\lim_{n\to \infty}\{ Risk_n^*(t_0)\} = 1$. 
\end{itemize}
\end{theorem}

Theorem~\ref{thm:minimax} implies that the minimax lower bound changes with the parameter class. When $t_0=0$, it is a very broad class, including the symmetric cases in Example 1 of Section~\ref{subsec:intro-phase} where the $\chi^2$ test loses power. In this broad class, the minimax lower bound is governed by $n^2\alpha_n^2\gamma_n^4$ (corresponding to the previous $\tau_n$), and the oSQ test is minimax optimal.   
When $t_0>0$, we restrict to a narrower class, with those extremely symmetrical settings excluded. The minimax lower bound is governed by $n^{3/2}\alpha_n\gamma_n^2$ (corresponding to the previous $\delta_n$), and the $\chi^2$-test is minimax optimal. In comparison, for both classes, the PE test is minimax optimal.

\section{The identifiability of $K$} \label{sec:identifiability}
To our best knowledge, the identifiability of MMSBM has not yet been carefully studied. We present our results, which are of independent interest.  

\bed \label{def:eligibility}
Fix $n \geq K$. We call $\Pi \in \mathbb{R}^{n\times K}$ an eligible membership matrix if each row is a weight vector and the $K \times K$ identity matrix $I_K$ is a sub-matrix of $\Pi$. We call $P\in[0,1]^K$ eligible if it is entry-wise non-negative. 
\eed

In MMSBM, $\Omega=\Pi P\Pi'$, for some $K$ and eligible $(\Pi, P)$. However, there may exist $K^*\neq K$ and eligible $(\Pi^*, P^*)$ such that $\Omega=\Pi^*P^*(\Pi^*)'$ also holds. 
Below is an example. 


\smallskip
\noindent
{\bf Example 3}. Let $K = 4$ and $K^*=2$. Fix any mixed membership matrix $\Pi=[\pi_1,\pi_2,\ldots,\pi_n]'\in\mathbb{R}^{n\times 4}$. Introduce $P$, $P^*$ and $\Pi^*=[\pi_1^*, \pi_2^*,\ldots,\pi_n^*]'\in\mathbb{R}^{n\times 2}$ as follows:   
\[
P =  0.01 \times\begin{bmatrix}
1  &  2  &  1.8 &   3\\
2   &  4  &  3.6  &   6\\
1.8  &  3.6  &  3.24  &  5.4\\
3  &  6 &   5.4 &   9
\end{bmatrix}, 
\quad 
P^*=0.01\times\begin{bmatrix}9 &\; 3\\ 3 &\; 1\end{bmatrix},
\quad 
\pi_i^*= \begin{bmatrix} 
0.5\pi_i(2)+0.4\pi_i(3)+\pi_i(4)\\
\pi_i(1)+0.5\pi_i(2)+0.6\pi_i(3)
\end{bmatrix}. 
\]

We connect $(K,\Pi, P)$ and $(K^*, \Pi^*, P^*)$. Let $\eta^*=(0.3, 0.1)'$, $h_1=(0, 0.5, 0.4, 1)'$, $h_2=(1, 0.5, 0.6, 0)'$, and $H=[h_1, h_2]\in\mathbb{R}^{4\times 2}$. It is straightforward to see  that $P^*=\eta^*(\eta^*)'$ and $\Pi^*=\Pi H$. Furthermore, we can verify that $P=(H\eta^*)(H\eta^*)'$. It follows that
\[ 
\Pi^* P^*(\Pi^*)' =(\Pi H) [\eta^*(\eta^*)'] (H'\Pi') =\Pi P\Pi'=\Omega. 
\] 
We can view this example as an MMSBM with $K=4$ communities or an MMSBM with $K^*=2$ communities. At the same time, the rank of $\Omega$ is $r=1$. Which of $(K, K^*, r)$ shall we use as the correct definition of `number of communities'? 

\medskip

To address this issue, we follow Occam's razor to define the number of communities as the smallest $K$ that is compatible with the matrix $\Omega$.

\bed \label{def:INC}
The {\it Intrinsic Number of Communities (INC)} of 
$\Omega$, denoted as $k_\Omega$, is the smallest integer $K$ such that $\Omega = \Pi P \Pi'$ for some eligible $\Pi \in \mathbb{R}^{n\times K}$ and $P \in \mathbb{R}^{K\times K}$. 
\eed

The next proposition is proved in the supplementary material: 
\begin{prop}[Identifiability of parameters of MMSBM]  \label{prop:iden} 
Suppose $\Omega = \Pi_0 P_0 \Pi_0'$ for some eligible $(\Pi_0, P_0, K_0)$ as in Definition~\ref{def:eligibility}. Recall that the INC $k_{\Omega}$ is as in Definition~\ref{def:INC}. 
\begin{itemize}
   \item  There exists a pair of eligible $\Pi \in \mathbb{R}^{n, k_{\Omega}}$ and $P \in \mathbb{R}^{k_{\Omega}, k_{\Omega}}$ such that $\Omega = \Pi P \Pi'$.
   \item If there is another pair of eligible $\Pi^* \in \mathbb{R}^{n, k_{\Omega}}$ and $P^* \in \mathbb{R}^{k_{\Omega}, k_{\Omega}}$ such that $\Omega = \Pi^*  P^*  (\Pi^*)'$, then there must exist a permutation matrix $D \in \mathbb{R}^{k_{\Omega},k_{\Omega}}$ such that $P^* = D PD' $.
   Therefore, when $K_0=k_{\Omega}$, $P$ is identifiable up to permutation.
   \item If, in addition, $ \mathrm{rank}(P)=k_\Omega$, then both $P$ and $\Pi$ are identifiable up to permutation.
   \item It holds that $K_0\geq k_{\Omega}\geq \mathrm{rank}(P_0)$. If $P_0$ is non-singular, then $k_{\Omega}=\mathrm{rank}(P_0)=K_0$. 
   \end{itemize}
\end{prop} 

Throughout this paper, we assume that the $K$ in the alternative hypothesis is the INC defined above. 

\smallskip

The definition of INC is natural, but it is not easy to compute. We introduce an alternative formula, which connects $k_{\Omega}$ with the geometry associated with the eigenvectors of $\Omega$. This equivalent definition is much more convenient to use. 
\begin{prop}[Equivalent definition of INC]  \label{prop:INC} 
Fix $K$ and $n\geq K$. 
Suppose $\Omega = \Pi P \Pi'$ for some eligible $\Pi\in\mathbb{R}^{n\times K}$ and $P\in\mathbb{R}^{K\times K}$ as in Definition~\ref{def:eligibility}. Let $r=\mathrm{rank}(\Omega)$, and let $\xi_1,\xi_2,\ldots,\xi_r\in\mathbb{R}^n$ be the eigenvectors of $\Omega$ associated with nonzero eigenvalues. Write $\Xi = [\xi_1, \xi_2, \ldots, \xi_r]$. Let $C(\Xi)\subset\mathbb{R}^r$ be the convex hull of the $n$ rows of $\Xi$. Then, $C(\Xi)$ is a polytope and $k_{\Omega}$ is equal to the number of vertices of this polytope. 
\end{prop} 

We apply Proposition~\ref{prop:INC} to get $k_{\Omega}$ in Example 3. 
In that example, we have seen that $\Omega$ is a rank-1 matrix, $\Omega=(\Pi H\eta^*)(\Pi H\eta^*)'$, implying that $k_{\Omega}\geq 1$. Additionally, $\Xi=\xi_1\propto \Pi H\eta^*$, and $C(\Xi)$ is an interval in $\mathbb{R}$ (a simplex with 2 vertices). It follows immediately that  $k_{\Omega}=2$. 

Using Proposition~\ref{prop:INC}, we can also easily see that the definition of $k_{\Omega}$ is compatible with the form of $\Omega$ for the null hypothesis. If $\Omega=\alpha_n {\bf 1}_n{\bf 1}_n'$, then it is obvious that $k_{\Omega}=1$. If $k_{\Omega}=1$, the last bullet point of Proposition~\ref{prop:iden} implies $\Omega=\lambda_1\xi_1\xi_1'$ and $C(\Xi)\subset\mathbb{R}$. By Proposition~\ref{prop:INC}, $C(\Xi)$ has to be a singleton, i.e., all the entries of $\xi_1$ are equal; hence, $\Omega\propto{\bf 1}_n{\bf 1}_n'$.

\section{Simulations} \label{sec:simul}

We conduct numerical experiments to investigate the behavior of the degree-based $\chi^2$ test, the orthodox Signed Quadrilateral (oSQ) test and the newly proposed Power Enhancement (PE) test. 
\footnote{All the simulation code can be found at: \url{https://github.com/louiscam/SBM_phase_transition.git}}

\paragraph*{Experiment 1: The asymptotic null distributions.}
We study how well the asymptotic null distributions in Corollary~\ref{cor:null} fit the simulated data, for a moderately large $n$. In these experiments, we generate networks from the Erd\"os-R\'enyi model, where $\Omega=\alpha{\bf 1}_n{\bf 1}_n'$.

\begin{figure}[h]
\begin{subfigure}{0.33\linewidth}
\includegraphics[width=1\linewidth]{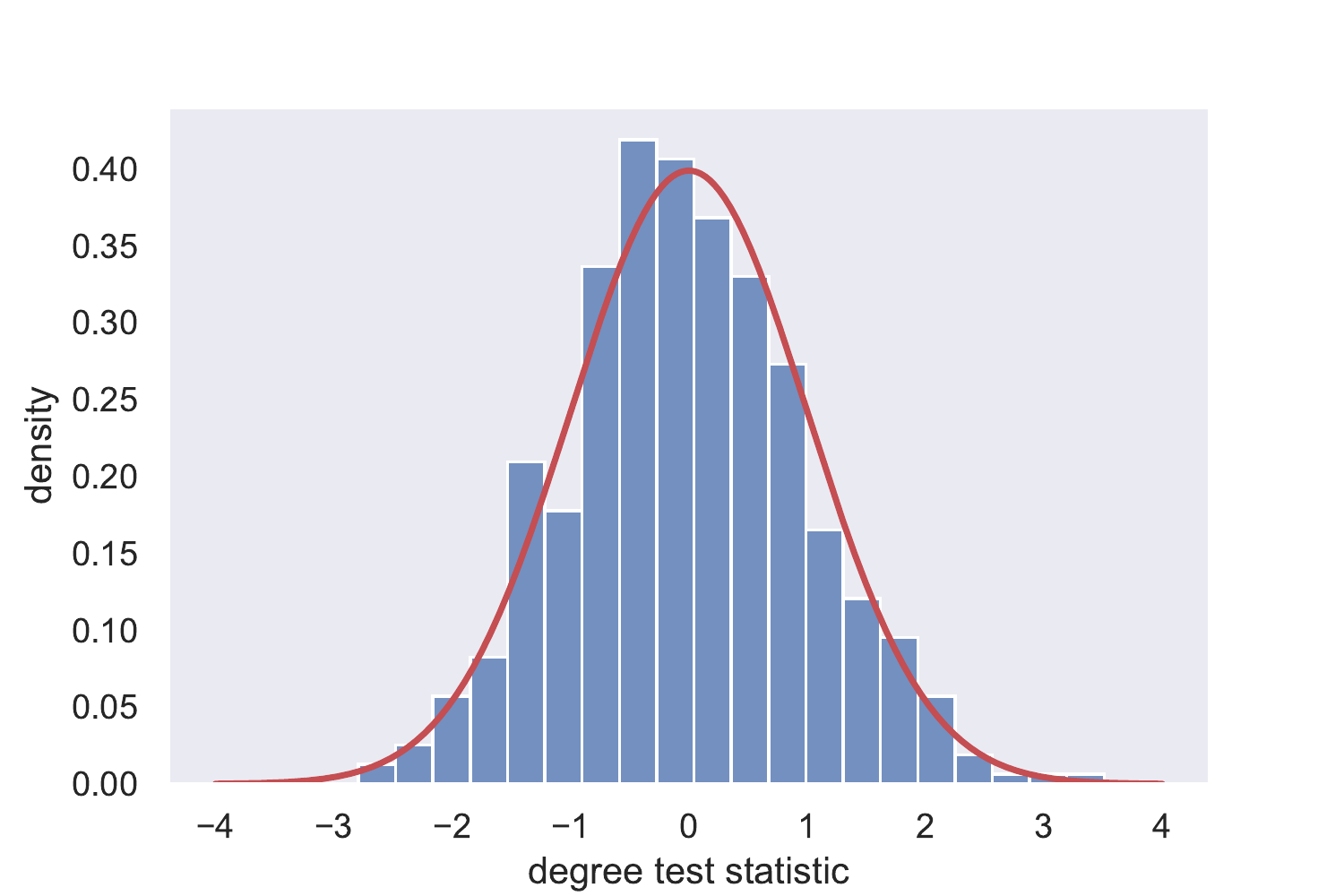} 
\caption{Degree-based $\chi^2$ test}
\label{fig:exp1_case1_deg}
\end{subfigure}
\begin{subfigure}{0.33\linewidth}
\includegraphics[width=1\linewidth]{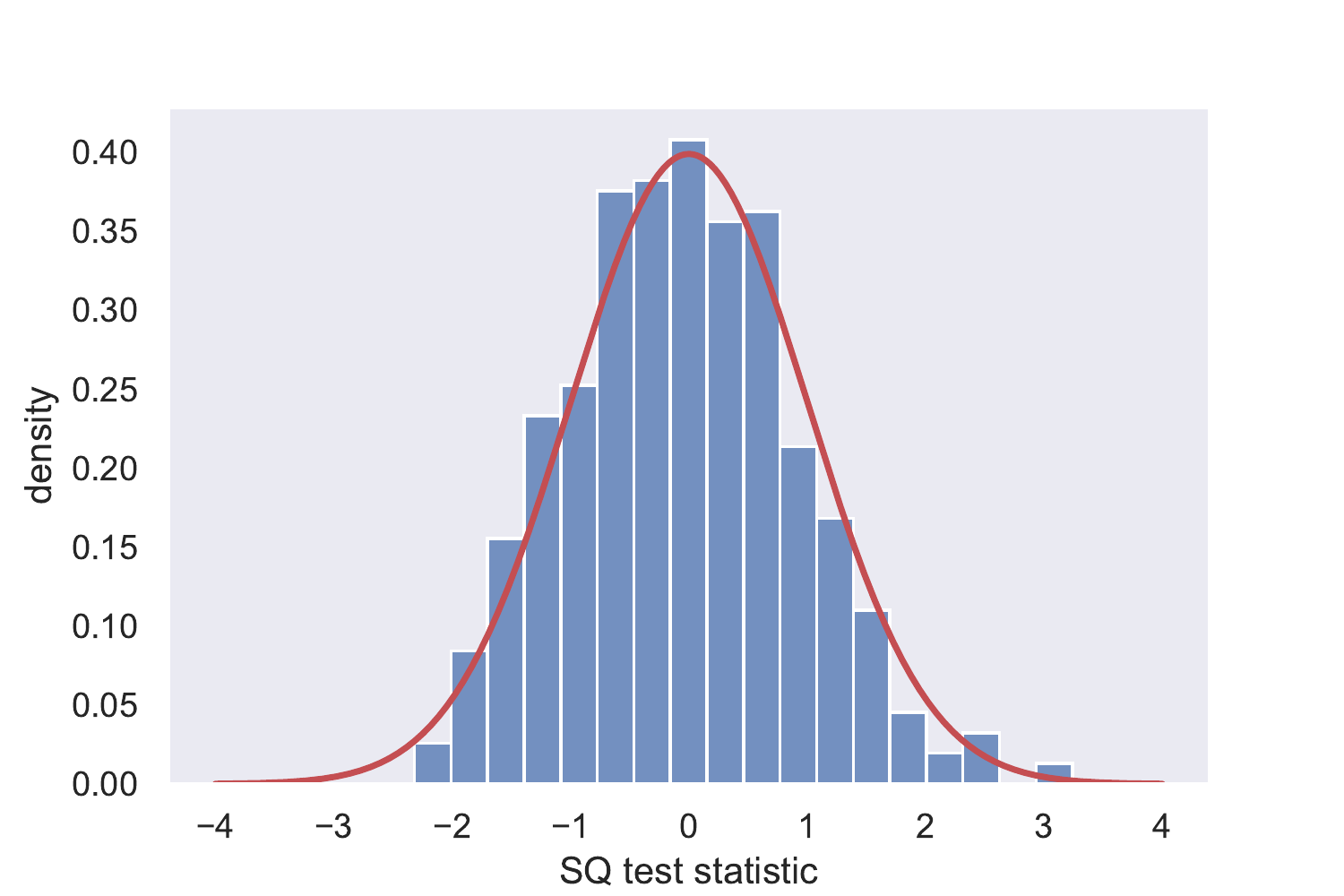}
\caption{oSQ test}
\label{fig:exp1_case1_sq}
\end{subfigure}
\begin{subfigure}{0.33\linewidth}
\includegraphics[width=1\linewidth]{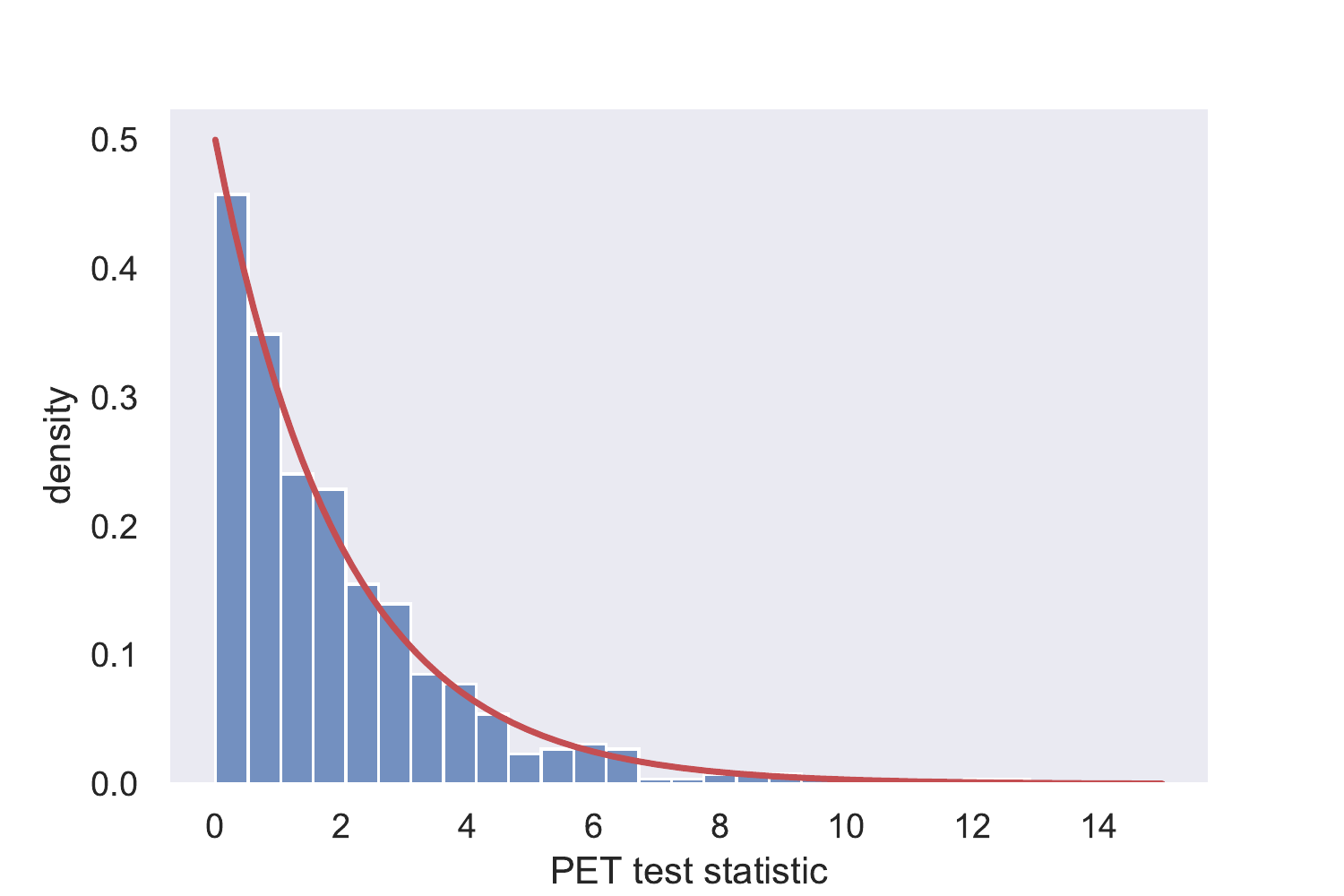}
\caption{PE test}
\label{fig:exp1_case1_pet}
\end{subfigure}
\caption{Histograms of the three test statistics under the null hypothesis ($n=200$). The red curves are the limiting null distributions in Corollary~\ref{cor:null}.}
\label{fig:Ex1-1}
\end{figure}

In Experiment 1.1, we fix $n=200$, $\alpha=0.1$ and generate 500 networks from the Erd\"os-R\'enyi model with $\Omega=\alpha{\bf 1}_n{\bf 1}_n'$. The histograms of the three test statistics are shown in Figure~\ref{fig:Ex1-1}.
They fit the limiting null distributions well. 

In Experiment 1.2, we focus on the PE test and evaluate the type I error when the target level is set at $5\%$. 
Given $(n,\alpha)$, we generate 500 networks, apply the level-$5\%$ PE test, and compute the empirical type I error. 
We let $n$ range in $\{100,200,500,1000\}$ and $\alpha$ range in $\{0.1, 0.2, 0.3, 0.4\}$. The results are shown in  Table~\ref{tb:Ex1-2}. It suggests that the type I errors are controlled satisfactorily. 
\begin{table}[hbt]
\centering
\caption{Empirical type I error of the level-$5\%$ PE test (calculated based on 500 repetitions).} \label{tb:Ex1-2}
\begin{tabular}{ccccc}
\hline
& $\alpha=0.1$ & $\alpha=0.2$ & $\alpha=0.3$ & $\alpha=0.4$\\
\hline
$n=100$ & 3.2\% & 3.4\% & 4\% & 3.2\% \\
$n=200$  & 4\% & 5.4\% & 6.2\% & 3.4\% \\
$n=500$ & 5.8\% & 3.8\% & 5.2\% & 4.2\%\\
$n=1000$ & 5\% & 6\% &5\% & 5.6\%\\
\hline
\end{tabular}
\end{table}

\paragraph*{Experiment 2: Power comparison of the three tests.}

We examine the power of the three tests and demonstrate the numerical advantage of PET over $\chi^2$ and oSQ. 
We will consider two settings, adapted from the examples in Section~\ref{subsec:whyPE}. In the first setting, $\delta_n=0$, hence, only the oSQ test has non-trivial power. In the second setting, $\tau_n\gg\delta_n$, hence, the $\chi^2$ test has a much larger SNR.

\begin{figure}[htb!]
\centering
\begin{subfigure}{0.3\linewidth}
\includegraphics[width=1\linewidth]{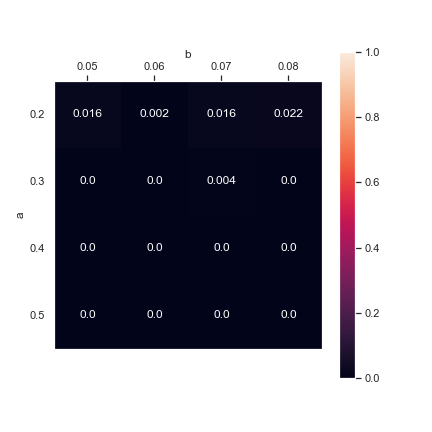} 
\caption{Degree-based $\chi^2$ test}
\label{fig:exp2_case1_deg}
\end{subfigure}
\begin{subfigure}{0.3\linewidth}
\includegraphics[width=1\linewidth]{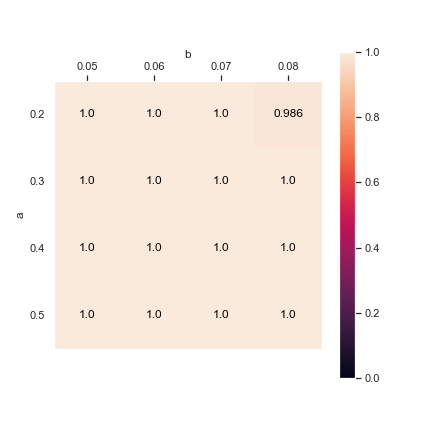}
\caption{SQ test}
\label{fig:exp2_case1_sq}
\end{subfigure}
\begin{subfigure}{0.3\linewidth}
\includegraphics[width=1\linewidth]{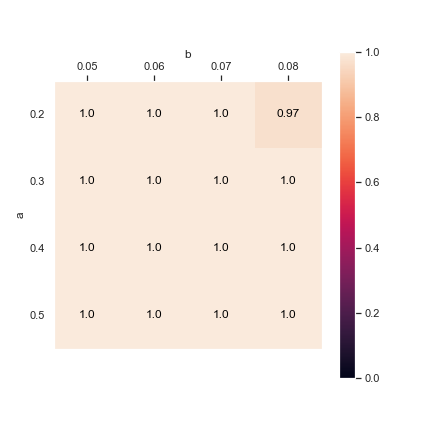}
\caption{PE test}
\label{fig:exp2_case1_pet}
\end{subfigure}
\caption{Empirical power of the three tests in Experiment 2.1 ($n=300$, $K=5$). The diagonals of $P$ are equal to $a$ and the off-diagonals are equal to $b$. Different combinations of $(a,b)$ are considered.
In this setting, $\delta_n=0$ always holds. Therefore, only the oSQ test has a non-trivial power.}
\label{fig:Ex2-1}
\end{figure}

In Experiment 2.1, we let $P = (a-b)I_K+b\bm{1}_K\bm{1}_K'$, same as in Example 1 of Section~\ref{subsec:whyPE}. We let all nodes be pure, with the same number of nodes in each community; this corresponds to $h=K^{-1}\bm{1}_K$. By direct calculations, $Ph=\alpha_0\bm{1}_K$, with $\alpha_0=h'Ph=K^{-1}[a+(K-1)b]$. In this setting, $\delta_n=0$, so that the $\chi^2$ test loses power. 
We fix $(n,K)=(300,5)$ and let $a$ range in $\{0.2, 0.3,0.4, 0.5\}$ and $b$ range in $\{0.05, 0.06, 0.07, 0.08\}$. For each $(a,b)$, we generate 500 networks from the alternative model as above and apply the three tests for a target level $5\%$; then, we report the proportion of rejections. The results are shown in Figure~\ref{fig:Ex2-1}. 
We observe that the SQ test clearly outperforms the degree-based $\chi^2$ test across these configurations, and that the PE test also benefits from this desirable behavior.


\begin{figure}[htb!]
\centering
\begin{subfigure}{0.3\linewidth}
\includegraphics[width=1\linewidth]{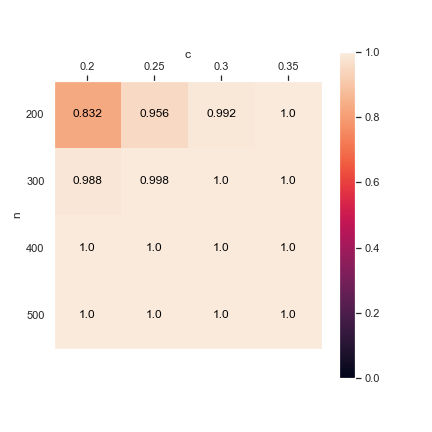} 
\caption{Degree-based $\chi^2$ test}
\label{fig:exp2_case2_deg}
\end{subfigure}
\begin{subfigure}{0.3\linewidth}
\includegraphics[width=1\linewidth]{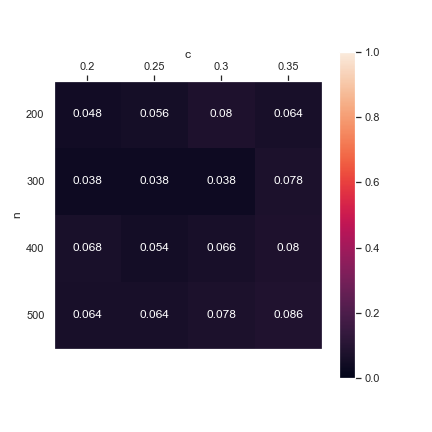}
\caption{SQ test}
\label{fig:exp2_case2_sq}
\end{subfigure}
\begin{subfigure}{0.3\linewidth}
\includegraphics[width=1\linewidth]{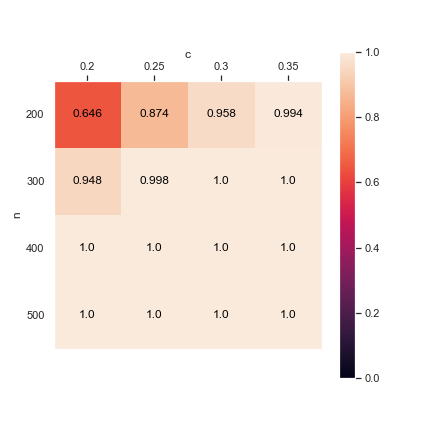}
\caption{PE test}
\label{fig:exp2_case2_pet}
\end{subfigure}
\caption{Empirical power of the three tests in Experiment 2.2 ($K=2$), where $P=c\eta\eta'$, for a vector $\eta\not\propto {\bf 1}_K$. Different combinations of $(n, c)$ are considered. The vector $\eta$ is chosen such that $\delta_n$ is always much larger than $\tau_n$. Therefore, the $\chi^2$ test has a much higher SNR.}
\label{fig:exp2_case2}
\end{figure}

In Experiment 2.2, we let $P=c\eta\eta'$, for a vector $\eta\not\propto {\bf 1}_K$ (similar to Example 2 of Section~\ref{subsec:whyPE}).  We fix $K=2$. Let all nodes be pure, with the same number of nodes in each community; hence, $h=K^{-1}\bm{1}_K$. 
We parametrize  $\eta=(a/\sqrt{a^2+b^2},b/\sqrt{a^2+b^2})'$, where $b=1$ and $a=1+n^{-1/4}$. By direct calculations, $\delta_n\sim c n/4$ and $\tau_n\sim c^2n/4$. 
We then let $n$ range in $\{200,300,400,500\}$ and $c$ range in $\{0.2,0.25,0.3,0.35\}$. For these values of $c$, the SNR of the $\chi^2$ test is considerably larger than that of the oSQ test. 
Similarly as in Experiment 2.1, we set the target level at 5\% and calculate the empirical power based on 500 repetitions. The results are shown in Figure~\ref{fig:exp2_case2}. We observe that the degree-based $\chi^2$ test clearly outperforms the SQ test across these configurations, and that the PE test also benefits from this desirable behavior.

\paragraph*{Experiment 3: Phase transitions for PE.}
We focus on the PE test and examine its power when the SNR $\beta_n=\max\{\delta_n,\tau_n\}$ gradually increases. This reveals the phase transitions associated with PE. 

\begin{figure}[tbh]
\centering
\begin{subfigure}{0.45\linewidth}
\includegraphics[width=1\linewidth]{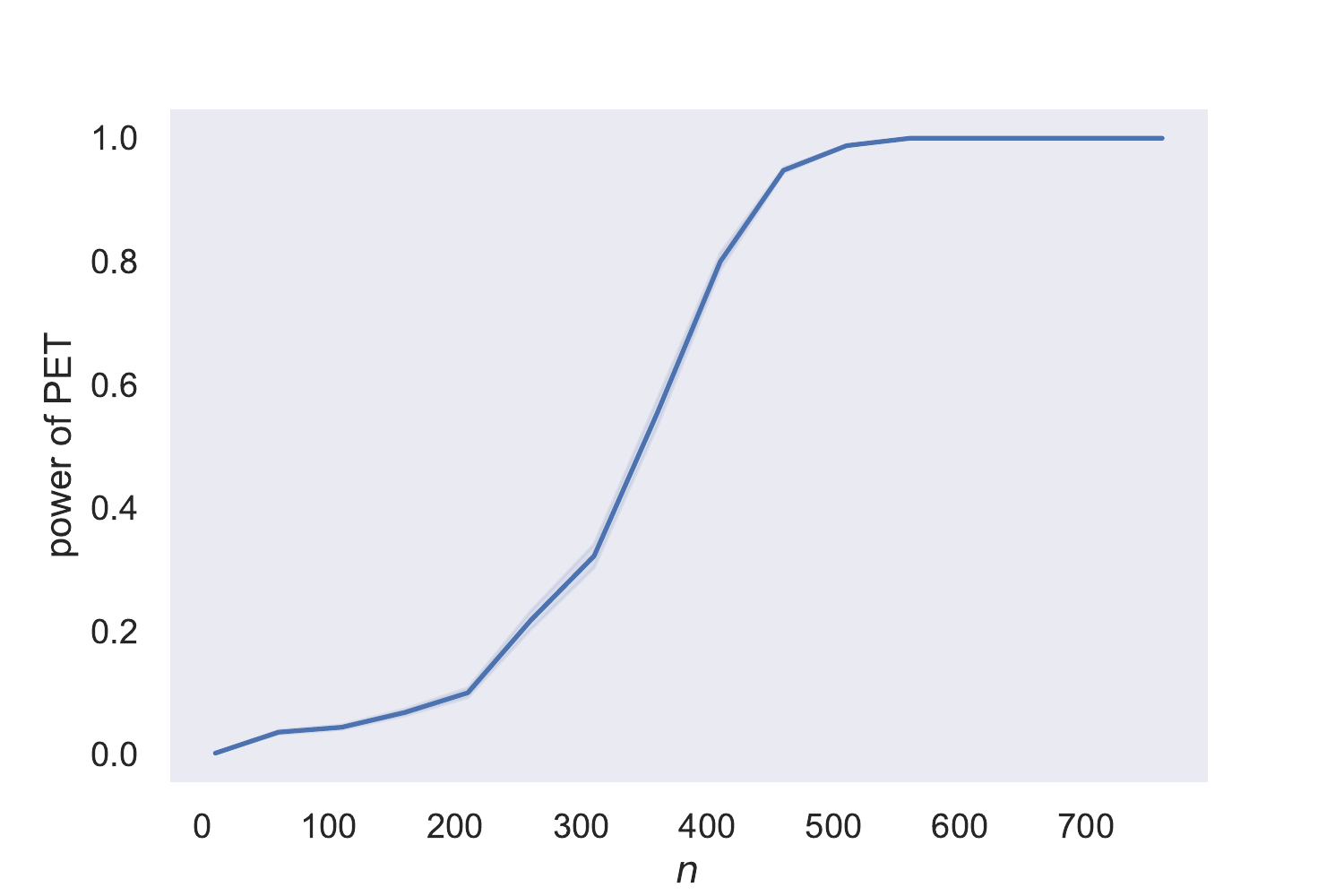} 
\caption{$b=0.1$, $n$ varies}
\label{fig:exp3_case1_n_PET}
\end{subfigure}
\begin{subfigure}{0.45\linewidth}
\includegraphics[width=1\linewidth]{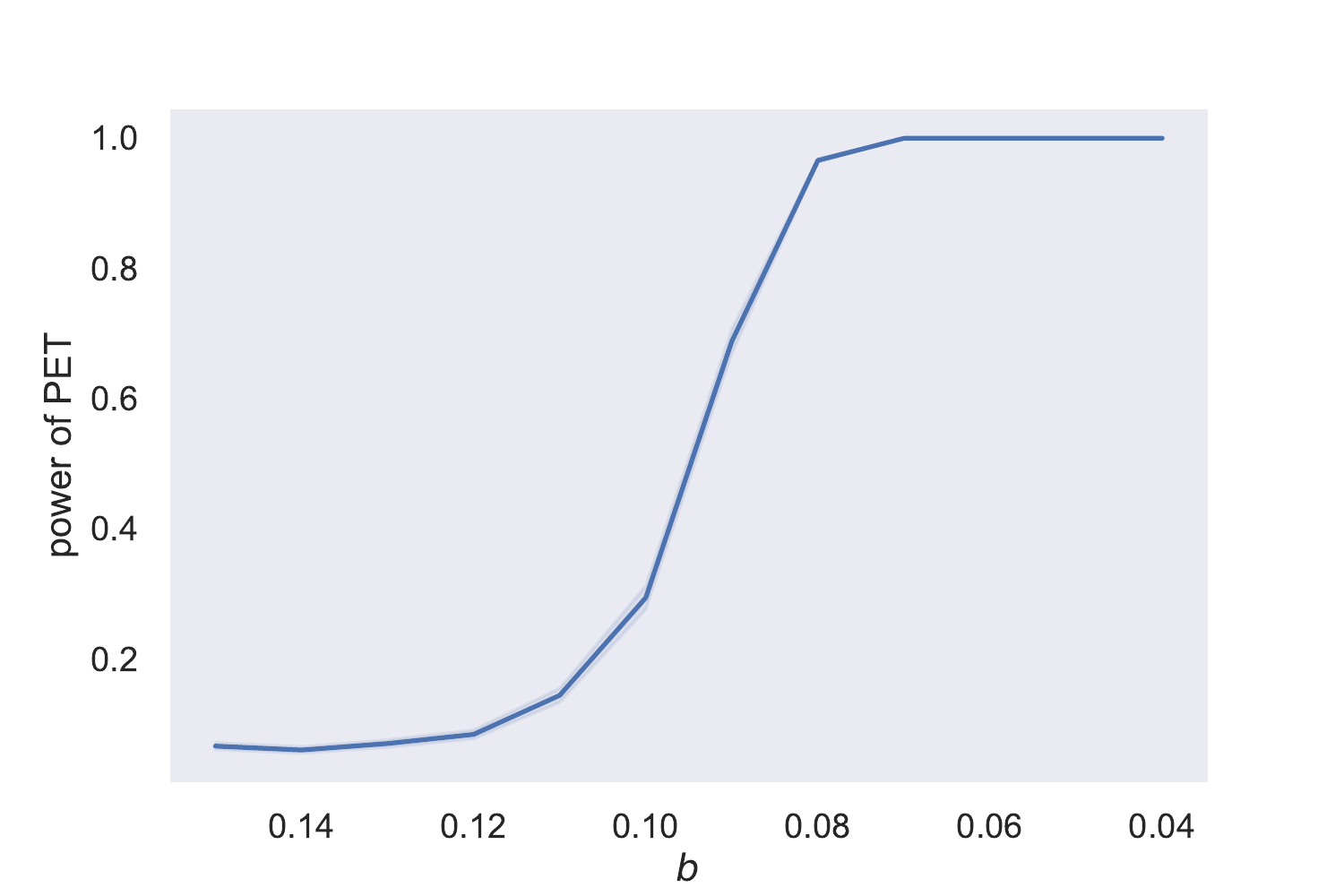}
\caption{$n=300$, $b$ varies}
\label{fig:exp3_case1_b_PET}
\end{subfigure}
\caption{Empirical power of the PE test in Experiment 3.1. In this setting, $\beta_n$ increases as $n$ increases or $b$ decreases.}
\label{fig:exp3_case1}
\end{figure}

In Experiment 3.1, we use the same model as in Experiment 2.1, where $K=5$ and $P = (a-b)I_K+b\bm{1}_K\bm{1}_K'$.  By direct calculations,  $\beta_n =n^2K^{2}[a+(K-1)b]^{-2}(a-b)^4$. We fix $a=0.2$. Then, $\beta_n$ is monotone increasing with $n$ and monotone decreasing with $b$. In Experiment 3.1(a), we fix $b=0.1$ and let $n$ vary from $10$ to $760$ with a step size of $50$. In Experiment 3.1(b), we fix $n=300$ and let $b$ vary from $0.04$ to $0.15$ with a step size of $0.01$. We report simulation results in Figure~\ref{fig:exp3_case1}, where power estimates for each configuration are obtained by averaging the number of rejections over 500 repetitions. The phase transition is visible as we move from vanishing power at low SNR to full power at high SNR.

\begin{figure}[tbh]
\centering
\begin{subfigure}{0.45\linewidth}
\includegraphics[width=1\linewidth]{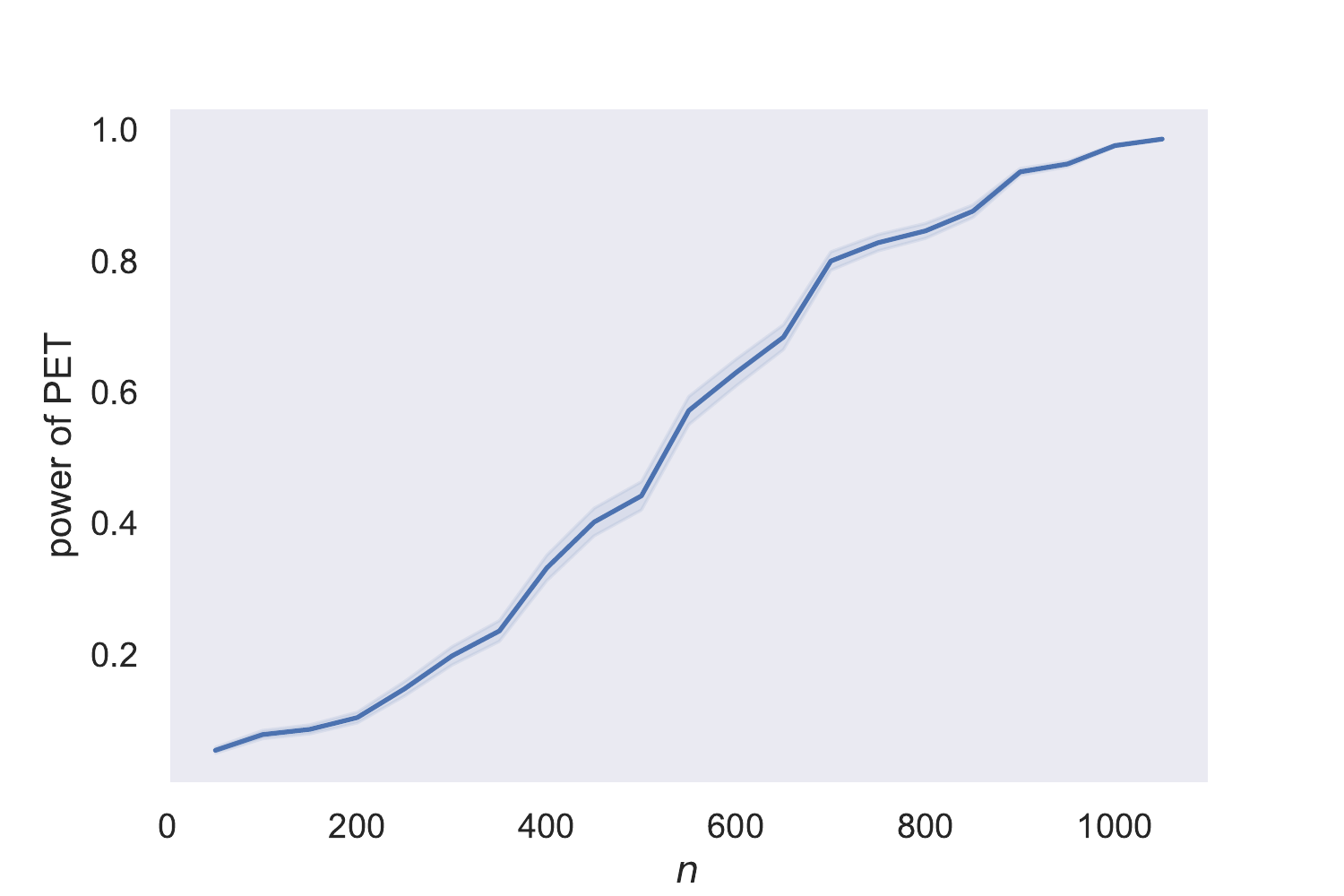} 
\caption{$c=0.06$, $n$ varies}
\label{fig:exp3_case2_n_PET}
\end{subfigure}
\begin{subfigure}{0.45\linewidth}
\includegraphics[width=1\linewidth]{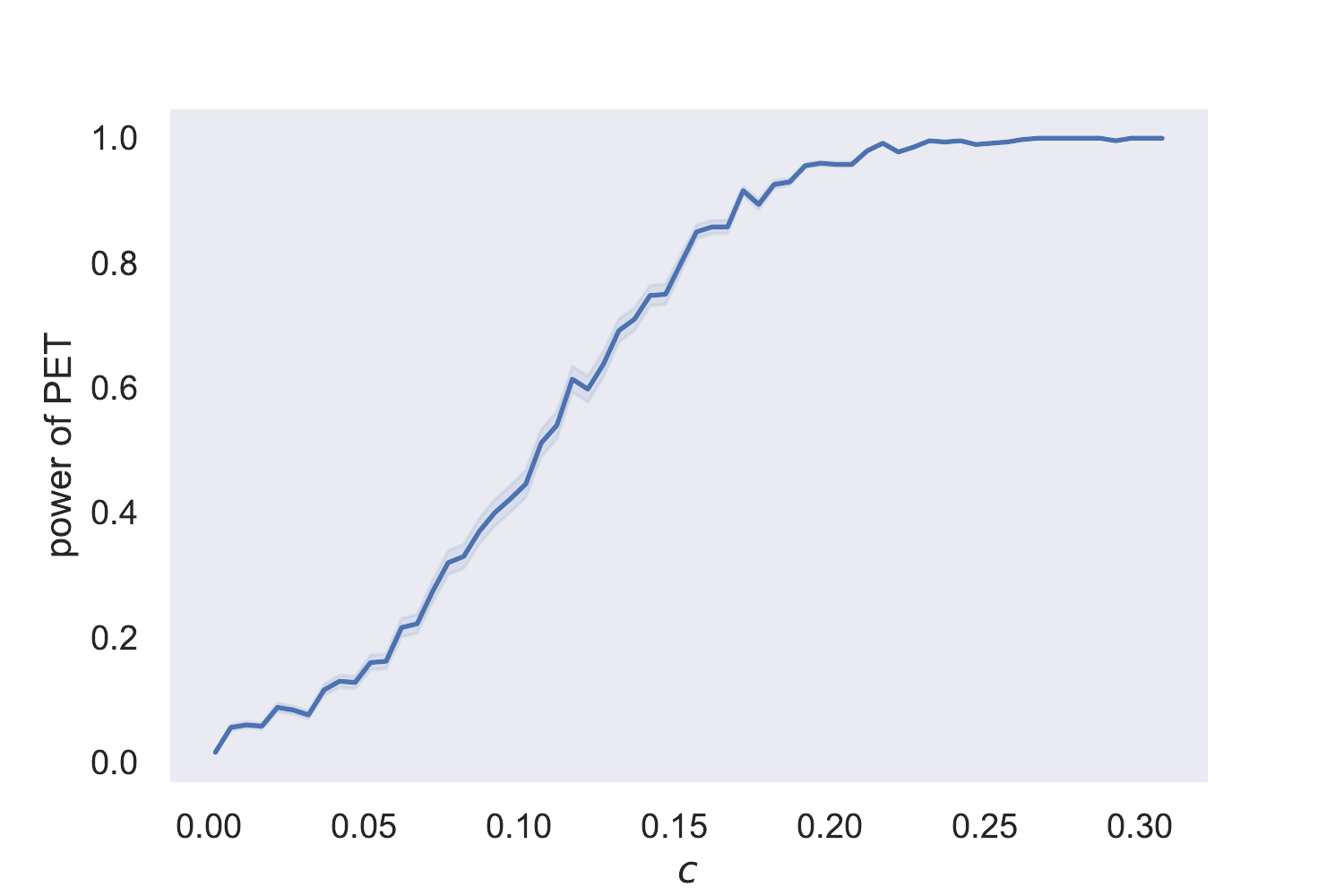}
\caption{$n=300$, $c$ varies}
\label{fig:exp3_case2_c_PET}
\end{subfigure}
\caption{Empirical power of the PE test in Experiment 3.2. In this setting, $\beta_n$ increases as $n$ increases or $c$ increases.}
\label{fig:exp3_case2}
\end{figure}

In Experiment 3.2, we use the same model as in Experiment 2.2, where $K=2$ and $P=c \eta\eta'$, with $\eta=(a, b)'/\sqrt{a^2+b^2}$. We fix $b=1$ and $a=1+n^{1/4}$. Then, $\beta_n\sim \max\{cn/4, \, c^2n/4\}$, which increases with both $n$ and $c$. 
In Experiment 3.2(a), we fix $c=0.06$ and let $n$ vary from $50$ to $1100$ with a step size of $50$. In Experiment 3.2(b), we fix $n=300$ and let $c$ vary from $0.002$ to $0.3$ with a step size of $0.005$. We report simulation results in Figure~\ref{fig:exp3_case2}. It also reveals the phase transition, from the vanishing power at low SNR to full power at high SNR.

\paragraph*{Experiment 4: Comparison with other testing ideas.} 
Other common ideas of global testing include the eigenvalue-based tests and the likelihood-ratio tests. For eigenvalue-based tests, we consider the one in Lei \cite{lei2016goodness}. The test statistic is a function of the largest and smallest eigenvalues of $A-\hat{\alpha}_n{\bf 1}_n{\bf 1}_n$.  \cite{lei2016goodness} showed that the test statistic converges to a Tracy-Widom distribution under the null hypothesis. We use this null distribution to set the rejection region. The use of likelihood-ratio tests has been limited to SBM (i.e., there is no mixed membership) and requires information on the unknown $K$ in the alternative hypothesis. We instead applied the model selection approach in Bickel and Wang \cite{wang2017likelihood}, which obtains $\hat{K}$ by successively computing the likelihood ratio between $K$ and $(K+1)$, for $K=1,2,\ldots$. We reject the null hypothesis if $\hat{K}>1$. 
In this approach, computing the likelihood ratios involves a sum over all possible community labels, and we followed \cite{wang2017likelihood} to use the EM algorithm with an initialization by spectral clustering. More details on our implementation of these methods are in our GitHub repository.

In Experiment 4.1, we study SBM with $K=2$.  We consider three models in the alternative hypothesis: (i) The symmetric SBM: $P=(a-b)I_2+b\bm{1}_2\bm{1}_2'$, with $a=0.2$ and $b=0.05$; the two communities have equal size. (ii) The asymmetric SBM: $P=(a-b)I_2+b\bm{1}_2\bm{1}_2'$, with $a=0.2$ and $b$ drawn from $\mathrm{Uniform}[0.125, 0.175]$; $\pi_i$'s are i.i.d. drawn from $\mathrm{Multinomial}(1, (0.2, 0.8)')$.   (iii) The rank-1 SBM: $P=\eta\eta'$, where $\eta=(a/\sqrt{a^2+b^2},b/\sqrt{a^2+b^2})'$, $b=1$ and $a=1+n^{-1/2}$; the two communities have equal size. Additionally, we consider the Erd\"os-R\'enyi model $\Omega=\alpha_0{\bf 1}_n{\bf 1}_n'$, with $\alpha_0=0.2$, as the null hypothesis. 
The $\chi^2$, oSQ, PE and eigenvalue tests have known null distributions, and we set the rejection region by controlling the level at 5\%. For the likelihood ratio test, as mentioned, we reject the null hypothesis if $\hat{K}>1$. For each model, we fix $n=500$, generate 100 networks, and measure the power of each test by the fraction of rejections over these 100 repetitions. In Experiment 4.2, we extend Models (i)-(iii) from SBM to MMSBM. For each model, $P$ is the same as before, except that $a=1+n^{-1/5}$ in Model (iii); $\pi_i$'s are i.i.d. generated from $\mbox{Dirichlet}(0.1,0.1)$ in Model (i),  $\mbox{Dirichlet}(0.2,0.8)$ in Model (ii), and $\mbox{Dirichlet}(0.4,0.6)$ in Model (iii). The results are in Table~\ref{tb:Ex4}.

\begin{table}[hbt]
\centering
\caption{Comparison with an eigenvalue-based test and a likelihood-ratio test. For each test, we report the empirical power over 100 repetitions. The settings are described in Experiments 4.1-4.2.} \label{tb:Ex4}
\scalebox{.8}{
\begin{tabular}{l|c|ccc|ccc}
\hline
\multirow{2}{*}{Test} & \multirow{2}{*}{Erd\"os-R\'enyi} & \multicolumn{3}{c|}{SBM} &  \multicolumn{3}{c}{MMSBM} \\
\cline{3-8}
& & Symmetric & Asymmetric & Rank-1 & Symmetric & Asymmetric & Rank-1  \\
\hline
$\chi^2$                                    & 0.04 & 0  & 0.96 & 0.95 & 0.09  & 0.87 & 0.99\\
oSQ                               & 0.05 & 1  & 0.33 & 0.04 & 1  & 0.06 & 0.03 \\
PE                                    & 0.06 & 1 &  0.92 & 0.88 &  1  & 0.76 & 0.98\\
\hline
Eigenvalue        & 0.06 & 1  & 0.56 & 0.02 & 1  & 0.10 & 0.31 \\
Likelihood     & 0.43 & 1  & 0.53 & 0.48 & 0.53  & 0.59 & 0.49 \\
\hline
\end{tabular}}
\end{table}

The `Symmetric' and `Asymmetric' models correspond to Case (S) and Case (AS1) in Example 1 of Section~\ref{subsec:whyPE}, and the `Rank-1' models correspond to Example 2. In Section~\ref{subsec:whyPE}, the SNRs of $\chi^2$, oSQ and PE have been analyzed, and their empirical powers here agree with the theoretical results. 
We now focus on comparing the eigenvalue test and the likelihood ratio test with the PE test.  
In all six models for the alternative hypothesis, the PE test outperforms the eigenvalue test and the likelihood ratio test. 
The eigenvalue test has a full power in symmetric SBM and symmetric MMSBM, but its performance is unsatisfactory in the other models. 
In fact, using the results in \cite{lei2016goodness}, we can derive that the SNR of the eigenvalue test is $\asymp \tau_n^{1/4}$; in comparison, the SNR of the PE test is $\max\{\delta_n,\tau_n\}$. Therefore, when $\delta_n\to\infty$ but $\tau_n\to 0$, the PE test has asymptotically full power but the eigenvalue test loses power (for example, in the rank-1 SBM model, $\delta_n\asymp n^{3/2}(a-b)^2\asymp n^{1/2}$, but $\tau_n^{1/4}\asymp n^{1/2}|a-b|\asymp 1\ll \delta_n$). 
The likelihood ratio test has better power than the eigenvalue test in the asymmetric and rank-1 models, but worse in the symmetric models. 
The likelihood ratio test also uniformly underperforms the PE test. For the SBM settings, the likelihood ratio test is supposed to have the best power, provided that $K$ is given and the likelihood is precisely computed. However, these requirements are practically infeasible. We had to use the model selection criteria \cite{wang2017likelihood} to avoid specifying $K$ and to compute the likelihood approximately, so its numerical performance should be inferior to the precise likelihood ratio test. For the MMSBM settings, the likelihood ratio is misspecified, which explains the unsatisfactory numerical performance. 
In terms of computing time, PE is also the fastest, especially for large $n$.

\section{Discussion} \label{sec:conclude}

We consider the global testing problem for MMSBM. First, we study the (degee-based) $\chi^2$ test and the oSQ test. These two tests existed in the literature, but their performances under MMSBM had never been studied. We derive their asymptotic null distributions and characterize their powers under the alternative. We discover that, for some parameter regimes, the $\chi^2$ test has a better performance; for some other parameter regimes, the oSQ test has a better performance. It motivates us to combine the strengths of both tests. Next, we propose the Power Enhancement (PE) test. We show that the PE test has a tractable null distribution and outperforms both the $\chi^2$ test and the oSQ test. 
Last, we study the phase transitions in global testing: We identify a quantity $\beta_n(K, P, h)$, such that the two hypotheses are asymptotically inseparable if $\beta_n(K,P,h)\to 0$, and perfectly separable by the PE test if $\beta_n(n,P,h)\to\infty$. This holds for arbitrary $(K,P,h)$ that satisfy mild regularity conditions. Therefore, the PE test is optimally adaptive.

Most existing works on global testing focused on a symmetric SBM \citep{mossel2015reconstruction, banerjee2017optimal}, which corresponds to a special choice of $(K, P, h)$ in our setting. The optimal test (e.g., the oSQ test or an eigenvalue-based test) for this special case may have unsatisfactory power for other choices of $(K, P, h)$. This motivates our study of phase transitions and optimal adaptivity, where we seek to understand the statistical limits for arbitrary $(K, P, h)$ and find a test that is both optimal and adaptive. 

In the hypothesis testing literature, it is not uncommon to combine multiple tests to attain the optimal detection boundary across the whole parameter range. For example, \citep{arias2014community} combines the $\chi^2$ test with a scan test for optimal detection of a planted clique, and \citep{JKW} combines a simple aggregation test and a sparse aggregation test for optimal global testing in a clustering model. However, these are  Bonferroni combinations, i.e., the combined test rejects the null hypothesis if any of the tests rejects. The simple Bonferroni combination does not support \textit{p}-value calculation. In contrast, our power enhancement test is based on the joint asymptotic distribution of two test statistics. As a result, the PE test has a tractable null distribution and supports the \textit{p}-value calculation.

In \citep{Yuan2018TestHyper}, the authors studied the global testing problem in hypergraphs under different sparsity settings. They introduced novel powerful statistical tests in the bounded degree regime and the dense regime. Potential avenues for future research include extending our power enhancement framework to the global testing problem in mixed-membership hypergraphs.

\bibliographystyle{Chicago} 
\bibliography{network}       

\tableofcontents

\newpage

\appendix

This supplemental material provides computations for examples and remarks, as well as proofs of theorems, corollaries and propositions. Appendix~\ref{subsec:Ex1-proof} covers the computations of $\tau_n$ and $\delta_n$ in Example 1, while Appendix~\ref{subsec:Ex2-proof} contains the calculation of the Intrinsic Number of Communities of the rank-1 model of Example 2, along with computations of $\tau_n$ and $\delta_n$ for that model. Appendix~\ref{subsec:Remark2-proof} shows the signal-to-noise ratios of the order-$m$ Signed Path and Signed Cycle statistics, for $m$ arbitrary. In Appendix~\ref{sec:null_joint_proof}, we derive the asymptotic joint null distribution of Theorem 2.1. Appendix~\ref{appendix:chi2-alt} shows the proof of Theorem 2.2, which consists in providing a lower bound for the expectation of the $\chi^2$ test statistic and an upper bound for its variance under the alternative hypothesis. Likewise, Appendix~\ref{appendix:proof_SQ_alt} derives the lower bound for the expectation of the oSQ test statistic and the upper bound for its variance under the alternative hypothesis,  presented in Theorem 2.3. Appendix~\ref{appendix:proof_cor_chi2} and Appendix~\ref{appendix:proof_cor_SQ} respectively report the proofs of Corollary 2.2 and Corollary 2.3 about the level and the power of the $\chi^2$ test and the oSQ test. The proof of Theorem 2.4 about the power and the level of the PE test is provided in Appendix~\ref{appendix:proof_PET_alt}. Appendix~\ref{appendix:LB} shows the proof of the lower bound, which corresponds to Theorem 2.5, and Appendix~\ref{appendix:proof_minimax} contains the proof of the minimax result of Theorem 2.6. Finally, Appendix~\ref{appendix:proofs_INC} shows the proof of Proposition 3.1 and Proposition 3.2  which examine the identifiability of MMSBM and give an alternative definition of the Intrinsic Number of Communities.

\section{Calculations in Example 1} \label{subsec:Ex1-proof}
Introduce
\[
y_n=1-2\epsilon_n, \qquad z_n=(d_n-a_n)/2. 
\]
Recall that $\bar{a}_n=(a_n+d_n)/2$, and $z_n=(d_n-a_n)/2$. Then, 
\[
P = (\bar{a}_n-b_n) I_2 - z_n e_1e_1' + z_n e_2e_2'+ b_n {\bf 1}_2{\bf 1}_2', \qquad h=\frac{1}{2}(1-y_n, 1+y_n)'. 
\] 
We calculate $\alpha_0$, $\|Mh\|^2$ and $\|M\|^2$ in general cases. 

First, consider $\alpha_0$. Note that $\|h\|^2=h_1^2+h_2^2=(1+y_n^2)/2$ and $h_2^2-h_1^2=y_n$. 
We have 
\begin{align} \label{Ex1-alpha0}
\alpha_0 = h'Ph &= h'\Bigl[(\bar{a}_n-b_n) I_2 - z_n e_1e_1' + z_n e_2e_2'+ b_n {\bf 1}_2{\bf 1}_2'\Bigr] h\cr
 &= (\bar{a}_n-b_n)\|h\|^2 +z_n(h_2^2-h_1^2) + b_n \cr
&= \bar{a}_n (1+y_n^2)/2 +z_ny_n + b_n(1-y_n^2)/2. 
\end{align}

Next, we calculate $\|Mh\|^2$. It follows from \eqref{Ex1-alpha0} that
\beq \label{Ex1-useful}
\alpha_0 - b_n = (\bar{a}_n-b_n) (1+y_n^2)/2 +z_ny_n.  
\eeq
We plug it into the expression of $Mh$ to get 
\begin{align*}
Mh = Ph - \alpha_0 {\bf 1}_2 & = \Bigl[ (\bar{a}_n-b_n) I_2 - z_n e_1e_1' + z_n e_2e_2'+ b_n {\bf 1}_2{\bf 1}_2' \Bigr] h - \alpha_0 {\bf 1}_2\cr
&= (\bar{a}_n-b_n)\begin{bmatrix}h_1\\h_2\end{bmatrix} - z_n \begin{bmatrix}h_1\\0\end{bmatrix} + z_n \begin{bmatrix}0\\h_2\end{bmatrix} + b_n {\bf 1}_2-\alpha_0{\bf 1}_2\cr
&= \begin{bmatrix}
(\bar{a}_n-b_n-z_n)h_1 \\
(\bar{a}_n-b_n+z_n)h_2 
\end{bmatrix} - (\alpha_0-b_n){\bf 1}_2\cr
&= \frac{1}{2}\begin{bmatrix}
(\bar{a}_n-b_n-z_n)(1-y_n) \\
(\bar{a}_n-b_n+z_n)(1+y_n) 
\end{bmatrix} - (\alpha_0-b_n){\bf 1}_2\cr
&= \frac{1}{2}\begin{bmatrix}
\bar{a}_n-b_n-z_n \\
\bar{a}_n-b_n+z_n
\end{bmatrix} +\frac{y_n}{2} \begin{bmatrix}
-(\bar{a}_n-b_n-z_n) \\
\bar{a}_n-b_n+z_n
\end{bmatrix} -(\alpha_0-b_n){\bf 1}_2\cr
&= \textcolor{black}{\frac{1}{2}\begin{bmatrix}
\bar{a}_n-b_n-z_n \\
\bar{a}_n-b_n+z_n
\end{bmatrix} +\frac{y_n}{2} \begin{bmatrix}
-(\bar{a}_n-b_n-z_n) \\
\bar{a}_n-b_n+z_n
\end{bmatrix}} \cr
&\textcolor{black}{-\frac{1}{2} \begin{bmatrix}
2z_ny_n \\
2z_ny_n
\end{bmatrix}-\frac{1}{2} \begin{bmatrix}
\bar{a}_n-b_n \\
\bar{a}_n-b_n
\end{bmatrix}-\frac{1}{2} \begin{bmatrix}
y_n^2(\bar{a}_n-b_n) \\
y_n^2(\bar{a}_n-b_n)
\end{bmatrix}}\cr
&\textcolor{black}{=\frac{z_n+y_n(\bar{a}_n-b_n)}{2}\begin{bmatrix}
-1 \\
1
\end{bmatrix}-\frac{y_nz_n+y_n^2(\bar{a}_n-b_n)}{2}\begin{bmatrix}
1 \\
1
\end{bmatrix}.}
\end{align*}
The two vectors, ${\bf 1}_2$ and $(1,-1)'$, are orthogonal to each other. It follows that
\begin{align} \label{Ex1-Mh}
\textcolor{black}{\|Mh\|^2} & \textcolor{black}{= \frac{1}{2}\Bigl[z_n+y_n(\bar{a}_n-b_n)\Bigr]^2+\frac{y_n^2}{2}\Bigl[z_n+y_n(\bar{a}_n-b_n)\Bigr]^2}\cr 
& \textcolor{black}{=\frac{1}{2}(1+y_n^2)\Bigl[z_n+y_n(\bar{a}_n-b_n)\Bigr]^2.}
\end{align}

Last, we calculate $\|M\|^2$. We have seen that
\[
M = (\bar{a}_n-b_n) I_2 - z_n e_1e_1' + z_n e_2e_2' - (\alpha_0- b_n) {\bf 1}_2{\bf 1}_2'. 
\]
Introduce $M_0= (\bar{a}_n-b_n) I_2-(\alpha_0- b_n) {\bf 1}_2{\bf 1}_2'$. Then, 
\beq \label{Ex1-temp1}
M = M_0 -z_n\diag(1, -1). 
\eeq
We compute the two eigenvalues of $M_0$. 
Write $v=(1,-1)'$. It is seen that $v$ is orthogonal to ${\bf 1}_2$; furthermore, 
\begin{align*}
& M_0v =\Bigl[(\bar{a}_n-b_n) I_2-(\alpha_0- b_n) {\bf 1}_2{\bf 1}_2'\Bigr]v = (\bar{a}_n-b_n)v \quad \propto\quad v, \cr
&M_0{\bf 1}_2 = \Bigl[(\bar{a}_n-b_n) I_2-(\alpha_0- b_n) {\bf 1}_2{\bf 1}_2'\Bigr]{\bf 1}_2 = [(\bar{a}_n-b_n) -2(\alpha_0- b_n)]{\bf 1}_2\quad\propto\quad {\bf 1}_2.
\end{align*}
It follows that ${\bf 1}_2$ and $v$ are two eigenvectors of $M^*$, with the associated eigenvalues as
\begin{align} \label{Ex1-temp2}
\lambda_1(M_0) &= (\bar{a}_n-b_n),\cr
\lambda_2(M_0) &= (\bar{a}_n-b_n) -2(\alpha_0-b_n)\cr
&= (\bar{a}_n-b_n)- [ (\bar{a}_n-b_n) (1+y_n^2) +2z_ny_n]\cr
&= -(\bar{a}_n-b_n)y_n^2 - 2z_ny_n, 
\end{align}
where we have applied \eqref{Ex1-useful} in the last equality. 
Combining \eqref{Ex1-temp1}-\eqref{Ex1-temp2}, we have
\beq \label{Ex1-Mnorm}
\|M\| \sim 
\begin{cases}
|z_n|, &\mbox{if }|z_n|\gg |\bar{a}_n-b_n|, \cr
|\bar{a}_n-b_n|, &\mbox{if } |z_n|\ll |\bar{a}_n-b_n|. 
\end{cases}
\eeq

We now combine \eqref{Ex1-alpha0}, \eqref{Ex1-Mh} and \eqref{Ex1-Mnorm}. In Case (S), $z_n=0$ and $y_n=0$. It follows that 
\[
\alpha_0=\textcolor{black}{\frac{a_n+b_n}{2}}, \qquad \|Mh\|^2=0, \qquad \|M\|^2= (\bar{a}_n-b_n)^2. 
\]
Plugging them into the definitions of $\delta_n$ and $\tau_n$ and noting that $\bar{a}_n=a_n$ in this case, we immediately get the claims for Case (S). In Case (AS1), \textcolor{black}{$\bar{a}_n=a_n$ and} $z_n=0$ but $y_n$ may be nonzero. It follows that 
\[
\textcolor{black}{\alpha_0=\frac{(1+y_n^2)a_n+(1-y_n^2)b_n}{2}, \qquad \|Mh\|^2=\frac{1}{2}(1+y_n^2)y_n^2(a_n-b_n)^2}, \qquad \|M\|^2= (a_n-b_n)^2. 
\]
\textcolor{black}{Assuming that $|a_n-b_n|=O(a_n+b_n)$, it follows that $(1+y_n^2)a_n+(1-y_n^2)b_n=(1+Cy_n^2)(a_n+b_n)$ for some constant $C>0$. We obtain}
\[\textcolor{black}{\alpha_0\asymp\frac{a_n+b_n}{2}, \qquad \|Mh\|^2\asymp\frac{1}{2}y_n^2(a_n-b_n)^2, \qquad \|M\|^2= (a_n-b_n)^2.}
\]
In Case (AS2), $y_n=0$ and $z_n\gg |\bar{a}_n-b_n|$. It follows that
\[
\textcolor{black}{\alpha_0= \frac{\bar{a}_n+b_n}{2}}, \qquad \|Mh\|^2=z_n^2/2, \qquad \|M\|^2 \sim z_n^2. 
\]
In Case (AS3), $y_n=0$ and $z_n\ll|\bar{a}_n-b_n|$. It follows that 
\[
\textcolor{black}{\alpha_0= \frac{\bar{a}_n+b_n}{2}}, \qquad \|Mh\|^2=z_n^2/2, \qquad \|M\|^2\sim (\bar{a}_n-b_n)^2. 
\]
The claims follow directly. \qed

\section{Calculations in Example 2}\label{subsec:Ex2-proof}

We start by showing that the rank-1 model of Example 2 has Intrinsic Number of Communities (INC) equal to $2$, regardless of $K$. We first recognize that the INC must be at least greater or equal to $2$. Indeed, suppose that the INC is equal to $1$, then we can find $\eta^*\in[0,1]$ such that $\Omega=(\eta^*)^2\bm{1}_n\bm{1}_n'$. From the original model formulation we had $\Omega=\Pi\eta\eta'\Pi'$, and we assumed that $\eta\not\propto\bm{1}_K$. Thus, it is impossible for $\Omega$ to have all equal entries if $\Pi$ is eligible, which contradicts the earlier fact that $\Omega=(\eta^*)^2\bm{1}_n\bm{1}_n'$, \textit{QEA}!

We now show that the INC is equal to $2$. Define
\begin{equation*}
    \eta^* = (\eta^*_1, \eta^*_2)'\in[0,1]^2, \quad \mbox{where} \quad 
    \begin{cases}
          \eta^*_1 = \max_{k\in\llbracket1,K\rrbracket}\eta_k\\
          \eta^*_2 = \min_{k\in\llbracket1,K\rrbracket}\eta_k.
    \end{cases}
\end{equation*}
We also define the matrix $H\in[0,1]^{K\times 2}$ such that
\begin{equation*}
    H = \frac{1}{\eta^*_1-\eta^*_2}
    \begin{pmatrix}
          \eta_1-\eta^*_2 & \eta^*_1-\eta_1\\
          \vdots & \vdots \\
          \eta_K-\eta^*_2 & \eta^*_1-\eta_K
    \end{pmatrix}.
\end{equation*}
It is straightforward to check that $H\eta^*=\eta$ and that $\Pi^*:=\Pi H$ is an eligible mixed membership matrix. It follows that
\begin{equation}
    \Omega = \Pi P \Pi' = \Pi \eta\eta'\Pi' = \Pi H\eta^*(\eta^*)'H'\Pi' = \Pi^*P^*(\Pi^*)',
\end{equation}
where we have defined the matrix $P^*=\eta^*(\eta^*)'\in[0,1]^{2\times 2}$. This shows that the INC of this rank-1 model is equal to $2$, regardless of $K\geq2$.

Next, we compute the Signal-to-Noise Ratios (SNR) of both tests for the rank-1 model introduced in Example 2. We start by computing the SNR of the degree test statistic, $\delta_n$. Recall that 
\begin{equation*}
    \delta_n := n^{3/2}\alpha_0^{-1}\|Ph-\alpha_0\bm{1}_K\|^2.
\end{equation*}
Direct calculations show that
\begin{equation*}
    P := \eta\eta' = \frac{c_n}{a_n^2+b_n^2}
    \begin{pmatrix}
        a_n^2 & a_nb_n\\
        a_nb_n & b_n^2
    \end{pmatrix}, \quad \mbox{and} \quad
    \alpha_0 := h'Ph = \frac{c_n(a_n+b_n)^2}{4(a_n^2+b_n^2)}.
\end{equation*}
This allows computing
\begin{equation*}
    Ph-\alpha_0\bm{1}_K = \frac{c_n(a_n+b_n)(a_n-b_n)}{4(a_n^2+b_n^2)}
    \begin{pmatrix}
          1\\
          -1
    \end{pmatrix}.
\end{equation*}
Together, the results for $\alpha_0$ and $Ph-\alpha_0\bm{1}_K $ yield the following expression of the SNR:
\begin{equation}
    \delta_n = \frac{1}{2}n^{3/2}c_n\frac{(a_n-b_n)^2}{(a_n^2+b_n^2)} \propto n^{3/2}c_n\frac{(a_n-b_n)^2}{(a_n^2+b_n^2)}.
\end{equation}

Then, we compute the SNR of the oSQ test statistic, $\tau_n$. Recall that
\begin{equation*}
    \tau_n := n^2\alpha_0^{-2}\|P-\alpha_0\bm{1}_K\bm{1}_K'\|^4.
\end{equation*}
We only need to compute $\|P-\alpha_0\bm{1}_K\bm{1}_K'\|$. Straightforward calculations reveal that
\begin{equation*}
    P-\alpha_0\bm{1}_K\bm{1}_K' = \frac{c_n(a_n-b_n)}{4(a_n^2+b_n^2)}
    \begin{pmatrix}
        3a_n+b_n & b_n-a_n\\
        b_n-a_n & 3b_n+a_n
    \end{pmatrix} =: \frac{c_n(a_n-b_n)}{4(a_n^2+b_n^2)}Q,
\end{equation*}
where we introduced the matrix $Q$ for notational convenience. The eigenvalues $\lambda_+$, $\lambda_-$ of $Q$ are the solutions to the following equation in the $x$-variable
\begin{equation*}
    x^2-\mbox{Tr}(Q)x+\mbox{det}(Q) = 0, \quad \mbox{where} \quad
    \begin{cases}
        \mbox{Tr}(Q) = 4(a_n+b_n)\\
        \mbox{det}(Q) = 2(a_n+b_n)^2+8a_nb_n.
    \end{cases}
\end{equation*}
We thus obtain that
\begin{equation*}
    \lambda_\pm = 2(a_n+b_n)\pm |a_n-b_n|, \quad \mbox{so} \quad \lambda_+ \asymp a_n+b_n,
\end{equation*}
where the last equivalence follows from our assumption that $|a_n-b_n| = O(a_n+b_n)$. It follows that
\begin{equation*}
    \|P-\alpha_0\bm{1}_K\bm{1}_K'\| \asymp \frac{c_n(a_n-b_n)(a_n+b_n)}{4(a_n^2+b_n^2)}.
\end{equation*}
As a consequence,
\begin{equation}
    \tau_n \asymp n^2c_n^2\frac{(a_n-b_n)^4}{(a_n^2+b_n^2)^2} \asymp n^{-1}\delta_n.
\end{equation}

\section{Calculations in Remark 2}\label{subsec:Remark2-proof}

\subsection{SNR of Signed Path statistics}

We consider the \textit{length-$m$ Signed Path statistic} $V^{(m)}_n$ defined as
\begin{equation*}
    V_n^{(m)} = \sum_{i_1,..., i_{m+1} \text{(distinct)}} (A_{i_1i_2}-\hat{\alpha}_n)(A_{i_2i_3}-\hat{\alpha}_n)...(A_{i_mi_{m+1}}-\hat{\alpha}_n), \quad \mbox{for $m\geq 2$},
\end{equation*}
where we recall that
\begin{equation*}
    \hat{\alpha}_n = \frac{1}{n(n-1)}\sum_{i\neq j}A_{ij}.
\end{equation*}
For simplicity, we study the corresponding ideal statistic $\bar{V}^{(m)}_n$, where we replace $\hat{\alpha}_n$ by the population null edge probability $\alpha_n$:
\begin{align*}
    \bar{V}_n^{(m)} &= \sum_{i_1,..., i_{m+1} \text{(distinct)}} (A_{i_1i_2}-\alpha_n)(A_{i_2i_3}-\alpha_n)...(A_{i_mi_{m+1}}-\alpha_n).
\end{align*}
The following lemma derives the null mean and variance as well as the alternative mean of the ideal length-$m$ Signed Path statistic. It uses the following quantities, which are defined in the main text:
\begin{equation*}
    h=\frac{1}{n}\sum_{i=1}^n\pi_i, \qquad \alpha_0=h'Ph, \qquad \mbox{and} \qquad G = \frac{1}{n}\Pi'\Pi.
\end{equation*}
In addition, we denote by $\E_1[\cdot]$ the expectation under the alternative distribution and by $\E_0[\cdot]$, $\mbox{Var}_0(\cdot)$ the expectation and variance under the null distribution, respectively.

\begin{lemma}[Moments of the ideal length-$m$ Signed Path statistic]\label{lm:sgnpath}
Suppose that conditions \eqref{cond1} and \eqref{cond2} hold. In addition, let $M=P-\alpha_0\bm{1}_K\bm{1}_K'$ and suppose that $n^{-1}\|Mh\|^{-1}\|M\|^2=o(1)$.  Then,
\begin{align*}
    \E_0\left[\bar{V}_n^{(m)}\right] = 0, \quad \mbox{Var}_0\left(\bar{V}_n^{(m)}\right) \asymp n^{m+1}\alpha_n^{m}, \quad \mbox{and} \quad \E_1\left[\bar{V}_n^{(m)}\right] \asymp n^{m+1}\|Mh\|^2\|M\|^{m-2}.
\end{align*}
\end{lemma}
\noindent \textit{Proof}\\
Under the null hypothesis, we can write
\begin{align*}
    \bar{V}_n^{(m)} &= \sum_{i_1,..., i_{m+1} \text{(distinct)}} W_{i_1i_2}W_{i_2i_3}...W_{i_mi_{m+1}},
\end{align*}
where $W_{ij}=A_{ij}-\alpha_n$ for all $i\neq j$.
It is straightforward to obtain that $ \E_0\left[\bar{V}_n^{(m)}\right] = 0$. Next, we compute the null variance of the ideal Signed Path statistic. We have, by direct calculations:
\begin{align}
    \mbox{Var}_0\left(\bar{V}_n^{(m)}\right) &= \mbox{Var}_0\left(\sum_{i_1,..., i_{m+1} \text{(distinct)}} W_{i_1i_2}W_{i_2i_3}...W_{i_mi_{m+1}}\right)\nonumber\\
    &= \E_0\left[\sum_{\substack{i_1,..., i_{m+1} \text{(distinct)}\\j_1,..., j_{m+1} \text{(distinct)}}}W_{i_1i_2}...W_{i_mi_{m+1}}W_{j_1j_2}...W_{j_mj_{m+1}}\right]\asymp n^{m+1}\alpha_n^m.
\end{align}

Under the alternative hypothesis, we choose $P$ and $h$ such that $\alpha_0:=h'Ph=\alpha_n$. This choice ensures that the network will have the same average degree under the null and alternative hypotheses, thus making the testing problem harder. As a result, we can write:
\begin{equation*}
    \bar{V}_n^{(m)} = \sum_{\substack{i_1,..., i_{m+1} \\ \text{(distinct)}}} (W_{i_1i_2}+\bar{\Omega}_{i_1i_2})(W_{i_2i_3}+\bar{\Omega}_{i_2i_3})...(W_{i_mi_{m+1}}+\bar{\Omega}_{i_mi_{m+1}}),
\end{equation*}
where $W_{ij}=A_{ij}-\Omega_{ij}$ and $\bar{\Omega}_{ij}=\Omega_{ij}-\alpha_0=\pi_i'M\pi_j$ for all $i\neq j$. It follows that
\begin{align*}
    \E_1\left[\bar{V}_n^{(m)}\right] &= \sum_{\substack{i_1,..., i_{m+1} \\ \text{(distinct)}}}\bar{\Omega}_{i_1i_2}\bar{\Omega}_{i_2i_3}...\bar{\Omega}_{i_mi_{m+1}}\\
    &= \sum_{i_1,..., i_{m+1}}\bar{\Omega}_{i_1i_2}\bar{\Omega}_{i_2i_3}...\bar{\Omega}_{i_mi_{m+1}}-\sum_{\substack{i_1,..., i_{m+1} \\ \text{(not distinct)}}}\bar{\Omega}_{i_1i_2}\bar{\Omega}_{i_2i_3}...\bar{\Omega}_{i_mi_{m+1}}\\
    &= \bm{1}_n'\bar{\Omega}^m\bm{1}_n-O(n^{m}\|M\|^m) = \bm{1}_n (\Pi M\Pi')^{m}\bm{1}_n'-O(n^{m}\|M\|^m)\\
    &=  n^{m-1}\bm{1}_n (\Pi M G M... M G M \Pi')\bm{1}_n-O(n^{m}\|M\|^m) \\
    &= n^{m+1} h'M(G M... M G)Mh-O(n^{m}\|M\|^m)
\end{align*}
Since we have assumed that $\|G\|, \|G^{-1}\| <c$ and $n^{-1}\|Mh\|^{-1}\|M\|^2=o(1)$, we obtain that 
\begin{equation}
     \E_1\left[\bar{V}_n^{(m)}\right] \asymp n^{m+1}\|Ph-\alpha_0\bm{1}_K
    \|^2\|P-\alpha_0\bm{1}_K\bm{1}_K\|^{m-2}.
\end{equation}
\qed

The results in Lemma~\ref{lm:sgnpath} allow us to compute the SNR for the length-$m$ Signed Path statistic. We derive the SNR assuming that the null variance dominates the alternative variance. Thus,
\begin{align*}
     SNR\left(\bar{V}_n^{(m)}\right) &= \frac{\left|\mathbb{E}_1\left[\bar{V}_n^{(m)}\right]-\mathbb{E}_0\left[\bar{V}_n^{(m)}\right]\right|}{\sqrt{\max\left\{\mbox{Var}_0\left(\bar{V}_n^{(m)}\right), \mbox{Var}_1\left(\bar{V}_n^{(m)}\right)\right\}}}\asymp \frac{\left|\mathbb{E}_1\left[\bar{V}_n^{(m)}\right]\right|}{\sqrt{\mbox{Var}_0\left(\bar{V}_n^{(m)}\right)}}\\
     &\asymp \frac{n^{m+1}\|Mh\|^2\|M\|^{m-2}}{n^{(m+1)/2}\alpha_0^{m/2}} = \delta_n\tau_n^{(m-2)/4}.
\end{align*}
Similar to our results in Theorem~\ref{thm:chi2-alt}, there may be instances in which the alternative variance dominates the null variance. In these cases, the SNR still depends on powers of $\delta_n$ and $\tau_n$, and the detection boundary is unchanged; details are omitted.

\subsection{SNR of Signed Cycle statistics}

We consider the \textit{length-$m$ Signed Cycle statistic} $U^{(m)}_n$ defined as
\begin{equation*}
    U_n^{(m)} = \sum_{i_1,..., i_{m} \text{(distinct)}} (A_{i_1i_2}-\hat{\alpha}_n)(A_{i_2i_3}-\hat{\alpha}_n)...(A_{i_mi_{1}}-\hat{\alpha}_n), \quad \mbox{for $m\geq 3$}.
\end{equation*}
For simplicity, we study the corresponding ideal statistic $\bar{U}^{(m)}_n$, where we replace $\hat{\alpha}_n$ by the population null edge probability $\alpha_n$:
\begin{align*}
    \bar{U}_n^{(m)} &= \sum_{i_1,..., i_{m} \text{(distinct)}} (A_{i_1i_2}-\alpha_n)(A_{i_2i_3}-\alpha_n)...(A_{i_mi_{1}}-\alpha_n), \quad \mbox{for $m\geq 3$}.
\end{align*}
\begin{lemma}[Moments of the ideal length-$m$ Signed Cycle statistic]\label{lm:sgncyc}
Suppose that conditions \eqref{cond1} and \eqref{cond2} hold. In addition, let $M=P-\alpha_0\bm{1}_K\bm{1}_K'$ and assume that $|\mbox{Tr}(MG)|\asymp \|MG\|$. Then,
\begin{align*}
    \E_0\left[\bar{U}_n^{(m)}\right] = 0, \quad \mbox{Var}_0\left(\bar{U}_n^{(m)}\right) \asymp n^m\alpha_n^m, \quad \mbox{and} \quad \left|\E_1\left[\bar{U}_n^{(m)}\right]\right| \asymp n^m\|M\|^m.
\end{align*}
\end{lemma}
\noindent \textit{Proof}\\
Under the null hypothesis, we can write
\begin{align*}
    \bar{U}_n^{(m)} &= \sum_{i_1,..., i_{m} \text{(distinct)}} W_{i_1i_2}W_{i_2i_3}...W_{i_mi_{1}},
\end{align*}
where $W_{ij}=A_{ij}-\alpha_n$ for all $i\neq j$. It is straightforward to obtain that $ \E_0\left[\bar{U}_n^{(m)}\right] = 0$. Next, we compute the null variance of the ideal Signed Cycle statistic. We have, by direct calculations:
\begin{align*}
    \mbox{Var}_0\left(\bar{U}_n^{(m)}\right) &= \mbox{Var}_0\left(\sum_{i_1,..., i_{m} \text{(distinct)}} W_{i_1i_2}W_{i_2i_3}...W_{i_mi_{1}}\right).
\end{align*}
Similar to Equation~\eqref{eq:cycdecomp}, we can decompose the sum into a sum over uncorrelated cycles. It results that
\begin{equation*}
    \mbox{Var}_0\left(\bar{U}_n^{(m)}\right) = C_m{n \choose m}\alpha_n^m(1-\alpha_n)^m \asymp n^m\alpha_n^m,
\end{equation*}
where $C_m$ is a constant that depends on $m$.

Under the alternative hypothesis, we can write
\begin{equation*}
    \bar{U}_n^{(m)} = \sum_{\substack{i_1,..., i_{m} \\ \text{(distinct)}}} (W_{i_1i_2}+\bar{\Omega}_{i_1i_2})(W_{i_2i_3}+\bar{\Omega}_{i_2i_3})...(W_{i_mi_{1}}+\bar{\Omega}_{i_mi_{1}}),
\end{equation*}
where $W_{ij}=A_{ij}-\Omega_{ij}$ and $\bar{\Omega}_{ij}=\Omega_{ij}-\alpha_0$ for all $i\neq j$. Then, direct calculations show that:
\begin{align*}
    \E_1\left[\bar{U}_n^{(m)}\right] &= \sum_{\substack{i_1,..., i_{m} \\ \text{(distinct)}}}\bar{\Omega}_{i_1i_2}\bar{\Omega}_{i_2i_3}...\bar{\Omega}_{i_mi_{1}} \\
    &= \sum_{i_1,..., i_{m}}\bar{\Omega}_{i_1i_2}\bar{\Omega}_{i_2i_3}...\bar{\Omega}_{i_mi_{1}}-\sum_{\substack{i_1,..., i_{m} \\ \text{(not distinct)}}}\bar{\Omega}_{i_1i_2}\bar{\Omega}_{i_2i_3}...\bar{\Omega}_{i_mi_{1}} \\
    &= \mbox{Tr}(\bar{\Omega}^m)-O\left(n^{m-1}\|M\|^m\right) = \mbox{Tr}\left((\Pi M\Pi')^{m}\right) -O\left(n^{m-1}\|M\|^m\right)\\
    &= n^m\mbox{Tr}\left((MG)^{m}\right)-O\left(n^{m-1}\|M\|^m\right) \asymp n^m\|MG\|^m-O\left(n^{m-1}\|M\|^m\right).
    \end{align*}
Since we have assumed that $\|G\|, \|G^{-1}\| <c$ by condition~\eqref{cond1}, we obtain that 
\begin{equation}
    \left|\E_1\left[\bar{U}_n^{(m)}\right]\right| \asymp n^m\|P-\alpha_0\bm{1}_K\bm{1}_K'\|^m.
\end{equation}
\qed

The results in Lemma~\ref{lm:sgncyc} allow us to compute the SNR for the length-$m$ Signed Cycle statistic. We derive the SNR assuming that the null variance dominates the alternative variance. Thus,
\begin{align*}
     SNR\left(\bar{U}_n^{(m)}\right) &= \frac{\left|\mathbb{E}_1\left[\bar{U}_n^{(m)}\right]-\mathbb{E}_0\left[\bar{U}_n^{(m)}\right]\right|}{\sqrt{\max\left\{\mbox{Var}_0\left(\bar{U}_n^{(m)}\right), \mbox{Var}_1\left(\bar{U}_n^{(m)}\right)\right\}}}\asymp \frac{\left|\mathbb{E}_1\left[\bar{U}_n^{(m)}\right]\right|}{\sqrt{\mbox{Var}_0\left(\bar{U}_n^{(m)}\right)}}\\
     &\asymp \frac{n^m\|M\|^m}{n^{m/2}\alpha_n^{m/2}} = \tau_n^{m/4}.
\end{align*}
Similar to our results in Theorem~\ref{thm:SQ-alt}, there may be instances in which the alternative variance dominates the null variance. In these cases, the SNR still depends on powers of $\tau_n$, and the detection boundary is unchanged; details are omitted.

\section{Proof of Theorem~\ref{thm:jointNull}}\label{sec:null_joint_proof}
Write $\varphi_n^{DC}=(X_n-n)/\sqrt{2n}$ and $\psi_n^{SQ}=Q_n/(2\sqrt{2}n^2\hat{\alpha}_n^2)$. We aim to show that $(\psi_n^{DC}, \psi_n^{SQ})$ converges to $\mathcal{N}(0, I_2)$ in distribution. By the Cram\'er-Wold theorem, it suffices to show that 
\begin{equation} \label{null-proof-goal}
u\cdot \psi_n^{DC}+v\cdot \psi_n^{SQ}\quad \xrightarrow[n\to\infty]{\mathcal{L}}\quad \mathcal{N}(0,1), \qquad\mbox{for any $u,v\in\mathbb{R}$ with $u^2+v^2=1$}. 
\end{equation}
Below, we first study the null distribution of $\psi_n^{DC}$ and $\psi_n^{SQ}$ respectively. These analyses produce useful intermediate results. We then use them to show the desirable claim in \eqref{null-proof-goal}. 

\subsection{Proof of the null distribution of $\psi_n^{DC}$} \label{subsec:nullproof-DC}
We aim to show that
\beq \label{null-proof-goal-DC}
\varphi_n^{DC}=\frac{X_n-n}{\sqrt{2n}}\quad \xrightarrow{d}\quad \mathcal{N}(0,1).  
\eeq

First, we derive an equivalent expression of $X_n$.  Let $\hat{T}_n=\sum_{i,j,k\text{ dist.}}(A_{ik}-\hat{\alpha}_n)(A_{jk}-\hat{\alpha}_n)$, where $\hat{\alpha}_n$ is the same as in the definition of $X_n$. We claim that
\beq \label{chi2-proof-1}
X_n = n + \frac{\hat{T}_n}{(n-1)\hat{\alpha}_n(1-\hat{\alpha}_n)}. 
\eeq
We now show \eqref{chi2-proof-1}. 
By definition,
\[
    X_n=\sum_{i=1}^n\frac{(d_i-\xbar{d})^2}{(n-1)\hat{\alpha}_n(1-\hat{\alpha}_n)},
\]
where
\begin{align*}
    \hat{\alpha}_n=\frac{1}{n(n-1)}\bm{1}_n'A\bm{1}_n & \text{,} & d = A\bm{1}_n & \text{,} & \xbar{d} =\frac{1}{n}\bm{1}_n'A\bm{1}_n= (n-1)\hat{\alpha}_n. 
\end{align*}
We expand $X_n$ into a sum of two terms that can be easily studied:
\begin{align*}
    X_n &= \frac{\|d\|_2^2-n\xbar{d}^2}{(n-1)\hat{\alpha}_n(1-\hat{\alpha}_n)} = \frac{\bm{1}_n'A^2\bm{1}_n}{(n-1)\hat{\alpha}_n(1-\hat{\alpha}_n)}-\frac{n(n-1)\hat{\alpha}_n}{1-\hat{\alpha}_n}.
\end{align*}
We can compute $\bm{1}_n'A^2\bm{1}_n$ as follows:
\begin{align*}
    \bm{1}_n'A^2\bm{1}_n &= \sum_{i,j}(A^2)_{ij} =\bm{1}_n'A\bm{1}_n+\sum_{i,j,k\text{ dist.}}A_{ik}A_{jk} = n(n-1)\hat{\alpha}_n+\sum_{i,j,k\text{ dist.}}A_{ik}A_{jk}.
\end{align*}
Hence we further reexpress $X_n$ as 
\begin{align*} 
    X_n &= \frac{\sum_{i,j,k\text{ dist.}}A_{ik}A_{jk}}{(n-1)\hat{\alpha}_n(1-\hat{\alpha}_n)}+\frac{n-n(n-1)\hat{\alpha}_n}{1-\hat{\alpha}_n}. 
\end{align*}
Recalling that $\hat{T}_n=\sum_{i,j,k\text{ dist.}}(A_{ik}-\hat{\alpha}_n)(A_{jk}-\hat{\alpha}_n)$, we have
\begin{align*}
    \sum_{i,j,k\text{ dist.}}A_{ik}A_{jk} &= \hat{T}_n+2(n-2)\hat{\alpha}_n\bm{1}_n'A\bm{1}_n-n(n-1)(n-2)\hat{\alpha}^2_n\\
    &= \hat{T}_n+n(n-1)(n-2)\hat{\alpha}_n^2.
\end{align*}
It follows that
\begin{align*}
    X_n-n&= \frac{\hat{T}_n+n(n-1)(n-2)\hat{\alpha}_n^2}{(n-1)\hat{\alpha}_n(1-\hat{\alpha}_n)}+\frac{n-n(n-1)\hat{\alpha}_n}{1-\hat{\alpha}_n}-n = \frac{\hat{T}_n}{(n-1)\hat{\alpha}_n(1-\hat{\alpha}_n)}.
\end{align*}
This proves \eqref{chi2-proof-1}. 

Next, we introduce an ideal counterpart to $\hat{T}_n$, $T_n=\sum_{i,j,k\text{ dist.}}(A_{ik}-\alpha_n)(A_{jk}-\alpha_n)$. Direct computations show that
\begin{align*}
    \mathbb{E}[T_n]=0 &\text{,} &\text{Var}(T_n)&= 2n(n-1)(n-2)\alpha_n^2(1-\alpha_n)^2.
\end{align*}
Thus
\begin{equation*}
    \text{Var}\left(\frac{T_n}{(n-1)\alpha_n(1-\alpha_n)}\right)=\frac{2n(n-2)}{n-1}.
\end{equation*}
Combining it with \eqref{chi2-proof-1}, we obtain
\begin{equation*}
    \frac{X_n-n}{\sqrt{2n}}=\left(\frac{\alpha_n(1-\alpha_n)}{\hat{\alpha}_n(1-\hat{\alpha}_n}\right)\left(\frac{\hat{T}_n}{T_n}\right)\left(\frac{n-2}{n-1}\right)^{1/2}\left(\frac{\frac{T_n}{(n-1)\alpha_n(1-\alpha_n)}}{\sqrt{\frac{2n(n-2)}{(n-1)}}}\right).
\end{equation*}
Define
\begin{align*}
    U_n = \frac{\alpha_n(1-\alpha_n)}{\hat{\alpha}_n(1-\hat{\alpha}_n)} & \text{,} & V_n = \frac{\hat{T}_n}{T_n} & \text{,} & Z_n=\frac{\frac{T_n}{(n-1)\alpha_n(1-\alpha_n)}}{\sqrt{\frac{2n(n-2)}{(n-1)}}}.
\end{align*}
We have the following decomposition: 
\begin{equation} \label{chi2-proof-2}
    \frac{X_n-n}{\sqrt{2n}}=\left(\frac{n-2}{n-1}\right)^{1/2}U_nV_nZ_n. 
\end{equation}
Below, we study $U_n$, $V_n$, and $Z_n$, separately. 

Consider $U_n$. Note that
\begin{equation*}
    \hat{\alpha}_n=\frac{1}{n(n-1)}\bm{1}_n'A\bm{1}_n=\frac{2}{n(n-1)}\sum_{i<j}A_{ij},
\end{equation*}
where $(A_{ij})_{i<j}$ are i.i.d. Bernoulli random variables with mean $\alpha_n$. By the Weak Law of Large Numbers we obtain that
\begin{equation} \label{chi2-proof-3}
    \frac{\hat{\alpha}_n}{\alpha_n}=\frac{2}{n(n-1)}\sum_{i<j}\frac{A_{ij}}{\alpha_n}\xrightarrow{\mathbb{P}}1,
\end{equation}
from which we conclude that $U_n\xrightarrow{\mathbb{P}}1$.

Consider $V_n$. Note that
\begin{align*}
    \hat{T}_n-T_n &= \sum_{i,j,k\text{ dist.}}(A_{ik}-\hat{\alpha}_n)(A_{jk}-\hat{\alpha}_n)-\sum_{i,j,k\text{ dist.}}(A_{ik}-\alpha_n)(A_{jk}-\alpha_n)\\
    &=\sum_{i,j,k\text{ dist.}} (\alpha_n-\hat{\alpha}_n)(A_{ik}+A_{jk}-\alpha_n-\hat{\alpha}_n)\\
    &= (\alpha_n-\hat{\alpha}_n)\left[2\left(\sum_{i,j,k\text{ dist.}}A_{ik}\right)-n(n-1)(n-2)(\alpha_n+\hat{\alpha}_n)\right]\\
    &= (\alpha_n-\hat{\alpha}_n)\left[2(n-2)\bm{1}_n'A\bm{1}_n-n(n-1)(n-2)(\alpha_n+\hat{\alpha}_n)\right]\\
    &=(\alpha_n-\hat{\alpha}_n)\left[2n(n-1)(n-2)\hat{\alpha}_n-n(n-1)(n-2)(\alpha_n+\hat{\alpha}_n)\right]\\
    &= -n(n-1)(n-2)(\alpha_n-\hat{\alpha}_n)^2. 
\end{align*}
It follows that
\begin{align*}
    \left|\frac{\hat{T}_n-T_n}{T_n}\right| &= \left|\frac{n(n-1)(n-2)(\alpha_n-\hat{\alpha}_n)^2}{T_n}\right|\\
    &= \sqrt{\frac{2(n-2)}{n(n-1)}}\left(\sqrt{\frac{n(n-1)}{2}}\frac{\hat{\alpha}_n-\alpha_n}{\sqrt{\alpha_n(1-\alpha_n)}}\right)^2\left|\frac{\sqrt{\frac{2n(n-2)}{n-1}}}{\frac{T_n}{(n-1)\alpha_n(1-\alpha_n)}}\right|\\
    &=\sqrt{\frac{n-2}{2(n-1)}}\left(\sqrt{\frac{n(n-1)}{2}}\frac{\hat{\alpha}_n-\alpha_n}{\sqrt{\alpha_n(1-\alpha_n)}}\right)^2\frac{1}{\sqrt{n}|Z_n|}.
\end{align*}
Note that $\hat{\alpha}_n=\frac{2}{n(n-1)}\sum_{i<j}A(i,j)$, where $A_{ij}$ are i.i.d. Bernoulli random variables with mean $\alpha_n$.  By the Central Limit Theorem, 
\begin{equation*}
    \sqrt{\frac{n(n-1)}{2}}\frac{\hat{\alpha}_n-\alpha_n}{\sqrt{\alpha_n(1-\alpha_n)}} \xrightarrow[n\to\infty]{\mathcal{L}} \mathcal{N}(0,1). 
\end{equation*}
We will show later that $Z_n\xrightarrow{\mathcal{L}}\mathcal{N}(0,1)$. It follows that  $(\sqrt{n}|Z_n|)^{-1}\xrightarrow{\mathbb{P}}0$ (by Slutsky's theorem) and we conclude by Slutsky's theorem again that 
\begin{equation} \label{chi2-proof-4}
    \left|\frac{\hat{T}_n-T_n}{T_n}\right|\xrightarrow{\mathbb{P}}0,
\end{equation}
which shows that $V_n\xrightarrow{\mathbb{P}}1$.

Consider $Z_n$.  
We define
\begin{equation*}
    I_m=\{(i,j,k)\in\llbracket 1,m \rrbracket^3 \text{ s.t. $i,j,k$ are distinct}\},
\end{equation*}
and the following quantities for $m\in\llbracket 1,n \rrbracket$
\begin{align*} 
    T_{n,m}&=\sum_{(i,j,k)\in I_m}W_{ik}W_{jk}, \qquad \text{ and } \qquad T_{n,0}=0,\cr
    Z_{n,m}&=\sqrt{\frac{n-1}{2n(n-2)}}\frac{T_{n,m}}{(n-1)\alpha_n(1-\alpha_n)}, \qquad \text{ and } \qquad Z_{n,0}=0.
\end{align*}
Consider the filtration $\{\mathcal{F}_{n,m}\}_{1\leq m\leq n}$ with $\mathcal{F}_{n,m}=\sigma\{W_{ij}, (i,j)\in\llbracket1,m\rrbracket^2\}$ for all $m\in\llbracket1,n\rrbracket$, $\mathcal{F}_{n,0}=\{\Omega,\emptyset\}$ (where $\Omega$ denotes the sample space). It is straightforward to see that for all $0\leq m\leq n$, $Z_{n,m}$ is $\mathcal{F}_{n,m}$-measurable, $\mathbb{E}[|Z_{n,m}|]<\infty$ and $\mathbb{E}\left[T_{n,m+1}|\mathcal{F}_{n,m}\right]=T_{n,m}$. This shows that $\{Z_{n,m}\}_{1\leq m\leq n}$ is a martingale with respect to $\{\mathcal{F}_{n,m}\}_{1\leq m\leq n}$. Define the martingale difference sequence, for all $m=1,...,n$
\begin{equation*}
    X_{n,m}=Z_{n,m}-Z_{n,m-1}.
\end{equation*}
With these notations we have $Z_n \equiv Z_{n,n}=\sum_{m=1}^nX_{n,m}$. Provided the following two conditions are met
\begin{align}
    \text{(a) }& \sum_{m=1}^n\mathbb{E}[X_{n,m}^2|\mathcal{F}_{n,m-1}]\xrightarrow{\mathbb{P}}1 \label{chi2-proof-5},\\
    \text{(b) }& \forall\epsilon>0, \sum_{m=1}^n\mathbb{E}[X_{n,m}^2\bm{1}\{|X_{n,m}>\epsilon|\}|\mathcal{F}_{n,m-1}]\xrightarrow{\mathbb{P}}0, \label{chi2-proof-6}
\end{align}
we conclude using the Martingale Central Limit Theorem that $Z_n\xrightarrow{\mathcal{L}}\mathcal{N}(0,1)$.

So far, we have shown that $Z_n\xrightarrow{\mathcal{L}}\mathcal{N}(0,1)$, $U_n\xrightarrow{\mathbb{P}}1$ and $V_n\xrightarrow{\mathbb{P}}1$. We plug them into \eqref{chi2-proof-2}. Then, \eqref{null-proof-goal-DC} follows immediately from Slutsky's theorem. 

The only remaining steps are to verify that \eqref{chi2-proof-5} and \eqref{chi2-proof-6} are indeed satisfied.

\smallskip\noindent 
{\bf Proof of Equation \eqref{chi2-proof-5}:}
It suffices to show that 
\begin{equation}     \label{a0}
    \E\left[\sum_{m=1}^n\mathbb{E}[X_{n,m}^2|\mathcal{F}_{n,m-1}]\right]=1,
\end{equation}
and
\begin{equation}     \label{b0}
    \V\left(\sum_{m=1}^n\mathbb{E}[X_{n,m}^2|\mathcal{F}_{n,m-1}]\right) \xrightarrow{n\to\infty}0. 
\end{equation}

First, we prove Equation~\eqref{a0}. For notational convenience we write
\begin{equation*}
    C_n:=(n-1)\alpha_n(1-\alpha_n)\sqrt{\frac{2n(n-2)}{n-1}}.
\end{equation*}
It follows that for all $n\in\mathbb{N}^*$ and $m\in\llbracket1,n\rrbracket$
\begin{equation*}
    C_nX_{n,m}=C_n(Z_{n,m}-Z_{n,m-1})=T_{n,m}-T_{n,m-1} = \sum_{(i,j,k)\in I_m\setminus I_{m-1}}W_{ik}W_{jk}. 
\end{equation*}
Triplets in $I_m\setminus I_{m-1}$ are such that one of the nodes is $m$: either one of the wingnodes $\{i,j\}$, or the centernode $k$. Hence,
\begin{equation} \label{chi2-proof-add}
    C_nX_{n,m}=2\sum_{\substack{1\leq j,k\leq m-1\\j\neq k}}W_{mk}W_{jk}+\sum_{\substack{1\leq i,j\leq m-1\\i\neq j}}W_{im}W_{jm}. 
\end{equation}
As a result (in the following, summations are all up to $m-1$)
\begin{equation*}
    C_n^2X_{n,m}^2 = 4\sum_{\substack{k\neq j\\i\neq l}}W_{mk}W_{jk}W_{mi}W_{il}+4\sum_{\substack{k\neq j\\i\neq l}}W_{mk}W_{jk}W_{im}W_{lm}+\sum_{\substack{i\neq j\\k\neq l}}W_{im}W_{jm}W_{km}W_{lm}.
\end{equation*}
It follows that
\begin{align} \label{chi2-proof-7}
    \E[C_n^2X_{n,m}^2|&\mathcal{F}_{n,m-1}] = 4\sum_{k\neq j;\; i\neq l}W_{jk}W_{il}\E\left[W_{mk}W_{mi}\right]
    + 4\sum_{k\neq j;\; i\neq l}W_{jk}\E\left[W_{im}W_{km}W_{lm}\right]\cr
    &+ \sum_{i\neq j;\; k\neq l}\E\left[W_{im}W_{jm}W_{km}W_{lm}\right]\cr
    =& 4\alpha_n(1-\alpha_n)\sum_{i}\sum_{j\neq i,l\neq i}W_{ij}W_{il}+2(m-1)(m-2)\alpha_n^2(1-\alpha_n)^2\cr
    =& 4\alpha_n(1-\alpha_n)\sum_{(i,j,l)\in I_{m-1}}W_{ij}W_{il} +4\alpha_n(1-\alpha_n)\sum_{i\neq j}W_{ij}^2\cr
    &+2(m-1)(m-2)\alpha_n^2(1-\alpha_n)^2\cr
    =& 4\alpha_n(1-\alpha_n)\Bigl(T_{n,m-1}+\sum_{i\neq j}W_{ij}^2\Bigr)+2(m-1)(m-2)\alpha_n^2(1-\alpha_n)^2. 
\end{align}
Let $\bm{1}_{n,m}\in\mathbb{R}^n$ be a vector whose $m$ first entries are 1, and whose remaining entries are 0. Define
\begin{equation*}
    \hat{\alpha}_{n,m}:=\frac{\bm{1}_{n,m}'A\bm{1}_{n,m}}{m(m-1)}.
\end{equation*}
By direct calculations, 
\begin{align*}
    \sum_{i\neq j}W_{ij}^2 &= \sum_{i\neq j}(A_{ij}-\alpha_n)^2=\sum_{i\neq j} \left[ A_{ij}(1-2\alpha_n)+\alpha_n^2\right] \\
    &= (m-1)(m-2)\alpha_n^2+(1-2\alpha_n)\sum_{i,j}A_{ij}\\
    &=(m-1)(m-2)\alpha_n^2+(1-2\alpha_n)(m-1)(m-2)\hat{\alpha}_{n,m-1}. 
\end{align*}
We plug the above equation into \eqref{chi2-proof-7} to get
\begin{align} \label{chi2-proof-8}
     \E[C_n^2X_{n,m}^2|\mathcal{F}_{n,m-1}] &= 4\alpha_n(1-\alpha_n)T_{n,m-1}+2(m-1)(m-2)\alpha_n^2(1-\alpha_n)^2\cr
     &+ 4(m-1)(m-2)\alpha_n^3(1-\alpha_n)\cr
     &+4(m-1)(m-2)\alpha_n\hat{\alpha}_{n,m-1}(1-\alpha_n)(1-2\alpha_n). 
\end{align}
It follows that 
\begin{align*}
    C_n^2\sum_{m=1}^n\E[X_{n,m}^2|\mathcal{F}_{n,m-1}]&= 4\alpha_n(1-\alpha_n)\sum_{m=1}^nT_{n,m-1}\\
     &+\left[ 2\alpha_n^2(1-\alpha_n)^2 +  4\alpha_n^3(1-\alpha_n)\right] \sum_{m=1}^n(m-1)(m-2)\\
     &+4\alpha_n(1-\alpha_n)(1-2\alpha_n)\sum_{m=1}^n(m-1)(m-2)\hat{\alpha}_{n,m-1}. 
\end{align*}
Recall that $\E[T_{n,m-1}]=\sum_{(i,j,k)\in I_{m-1}}\E[W_{ik}W_{jk}]=\sum_{(i,k)\in I_{m-1}}\E[W_{ik}^2]=\frac{(m-1)(m-2)}{2}\alpha_n(1-\alpha_n)$. Additionally,  $\E[\hat{\alpha}_{n,m-1}]=\alpha_n$.
We thus have
\begin{align*}
    C_n^2\E\left[\sum_{m=1}^n\E[X_{n,m}^2|\mathcal{F}_{n,m-1}]\right]&= 2\alpha_n^2(1-\alpha_n)^2\sum_{m=1}^n(m-1)(m-2)\\
     &+ \left[ 2\alpha_n^2(1-\alpha_n)^2 +  4\alpha_n^3(1-\alpha_n)\right]\sum_{m=1}^n(m-1)(m-2)\\
     &+4\alpha_n(1-\alpha_n)(1-2\alpha_n)\sum_{m=1}^n\alpha_n (m-1)(m-2) \\
     &= 6\alpha_n^2(1-\alpha_n)^2\sum_{m=1}^n(m-1)(m-2)\\
     &= 2\alpha_n^2(1-\alpha_n)^2n(n-1)(n-2)=C_n^2. 
\end{align*}
This proves \eqref{a0}. 

Second, we prove Equation \eqref{b0}. 
In the second line of \eqref{chi2-proof-7}, we have seen that
\begin{align*}
    C_n^2\E[X_{n,m}^2|\mathcal{F}_{n,m-1}] & = 4\alpha_n(1-\alpha_n)\sum_{k}\sum_{\substack{1\leq i\neq j\leq m-1\\i\neq k,j\neq k}}W_{ki}W_{kj}+2(m-1)(m-2)\alpha_n^2(1-\alpha_n)^2\cr
    & = 8\alpha_n(1-\alpha_n)\sum_{k}\sum_{\substack{1\leq i< j\leq m-1\\i\neq k,j\neq k}}W_{ki}W_{kj}+2(m-1)(m-2)\alpha_n^2(1-\alpha_n)^2.
\end{align*}
As a result,
\begin{align*}
    \V\left(C_n^2\sum_{m=1}^n\E[X_{n,m}^2|\mathcal{F}_{n,m-1}]\right) \leq  64\alpha_n^2\V\left(\sum_{m=1}^n\sum_{k}\sum_{\substack{1\leq i< j\leq m-1\\i\neq k,j\neq k}}W_{ki}W_{kj}\right). 
\end{align*}
Recall that in the previous sums, summation over the indices $i,j,k$ ranges from $1$ to $m-1$. We rearrange the terms of the sums in order to facilitate the computation of the variance. Instead of summing over the order $m$, then over centernodes $k$ ranging from $1$ to $m-1$, and finally over wingnodes $i,j$ also ranging from $1$ to $m-1$, we now sum over centernodes $k$ ranging from $1$ to $n-1$, wingnodes ranging from $1$ to $n-1$, and finally over orders $m>\max{(i,j,k)}$.
\begin{align*}
    & \V\left(C_n^2\sum_{m=1}^n\E[X_{n,m}^2|\mathcal{F}_{n,m-1}]\right) \leq 64\alpha_n^2\V\left(\sum_{k=1}^{n-1}\sum_{\substack{1\leq i<j\leq n-1\\i\neq k,j\neq k}}\sum_{m>\max(i,j,k)}W_{ki}W_{kj}\right)\\
    &\leq 64\alpha_n^2n^2\V\left(\sum_{k=1}^{n-1}\sum_{\substack{1\leq i<j\leq n-1\\i\neq k,j\neq k}}W_{ki}W_{kj}\right)
    =64\alpha_n^2n^2\sum_{k=1}^{n-1}\V\left(\sum_{\substack{1\leq i<j\leq n-1\\i\neq k,j\neq k}}W_{ki}W_{kj}\right),
\end{align*}
where the last equality comes from the fact that in the above sum, terms corresponding to different values of the index $k$ are uncorrelated. As a result
\begin{equation} \label{chi2-proof-9}
    \V\left(C_n^2\sum_{m=1}^n\E[X_{n,m}^2|\mathcal{F}_{n,m-1}]\right)\leq64\alpha_n^2n^2\sum_{k=1}^{n-1}\sum_{\substack{1\leq i<j\leq n-1\\1\leq u<v\leq n-1\\i,j,u,v\neq k}}\cov(W_{ki}W_{kj},W_{ku}W_{kv}).
\end{equation}
We examine the possible cases for $\cov(W_{ki}W_{kj},W_{ku}W_{kv})$.
\begin{itemize}
    \item Case 1: $(i,j)=(u,v)$, then $\cov(W_{ki}W_{kj},W_{ku}W_{kv})=\V(W_{ki}W_{kj})=\alpha_n^2(1-\alpha_n)^2$. 
    \item Case 2: $i=u,j\neq v$ or $i\neq u, j=v$, then $\cov(W_{ki}W_{kj},W_{ku}W_{kv})=0$.
    \item All other cases: $\cov(W_{ki}W_{kj},W_{ku}W_{kv})=0$.
\end{itemize}
It follows that
\begin{align*}
    \V\left(C_n^2\sum_{m=1}^n\E[X_{n,m}^2|\mathcal{F}_{n,m-1}]\right)&\leq64\alpha_n^2n^2\sum_{k=1}^{n-1}\sum_{\substack{1\leq i<j\leq n-1\\i,j\neq k}}\V(W_{ki}W_{kj})\\
    &= 64\alpha_n^2n^2\sum_{k=1}^{n-1}\sum_{\substack{1\leq i<j\leq n-1\\i,j\neq k}}\alpha_n^2(1-\alpha_n)^2 \leq 32\alpha_n^4n^5.
\end{align*}
Hence
\begin{align*}
     \V\left(\sum_{m=1}^n\E[X_{n,m}^2|\mathcal{F}_{n,m-1}]\right)&\leq\frac{32\alpha_n^4n^5}{C_n^4}= \frac{1}{n}\left(\frac{n^4}{(n-1)^2(n-2)^2}\right)\left(\frac{8}{(1-\alpha_n)^4}\right)\xrightarrow[n\to\infty]{} 0.
\end{align*}
This proves \eqref{b0}. 

\smallskip\noindent
{\bf Proof of Equation \eqref{chi2-proof-6}:} Notice that by the Cauchy-Schwarz and Markov inequalities we obtain the following upper bound
\begin{align*}
    \left|\sum_{m=1}^n\mathbb{E}[X_{n,m}^2\bm{1}\{|X_{n,m}>\epsilon|\}|\mathcal{F}_{n,m-1}]\right| &\leq \sum_{m=1}^n\sqrt{\mathbb{E}[X_{n,m}^4|\mathcal{F}_{n,m-1}]}\sqrt{\mathbb{P}(|X_{n,m}|>\epsilon|\mathcal{F}_{n,m-1})}\\
    &\leq \frac{1}{\epsilon^2}\sum_{m=1}^n\mathbb{E}[X_{n,m}^4|\mathcal{F}_{n,m-1}].
\end{align*}
Thus it suffices to show that $\sum_{m=1}^n\mathbb{E}[X_{n,m}^4|\mathcal{F}_{n,m-1}]\xrightarrow{\mathbb{P}}0$. Since these random variables are all non-negative, we will equivalently show that
\begin{equation} \label{chi2-proof-10}
    \sum_{m=1}^n\mathbb{E}[X_{n,m}^4] = \mathbb{E}\left[\sum_{m=1}^n\mathbb{E}[X_{n,m}^4|\mathcal{F}_{n,m-1}]\right] \xrightarrow{n\to\infty} 0.
\end{equation}

We now show \eqref{chi2-proof-10}. 
Recall that (see \eqref{chi2-proof-add})
\begin{align*}
    C_nX_{n,m}&=2\sum_{1\leq i<j\leq m-1}W_{im}W_{jm}+2\sum_{1\leq i<j\leq m-1}W_{ij}W_{jm}\\
    &=2\sum_{1\leq i<j\leq m-1}W_{jm}(W_{ij}+W_{im}).
\end{align*}
Then (with summations ranging from $1$ to $m-1$)
\begin{equation*} 
    C_n^4X_{n,m}^4 = 16\sum_{\substack{i<j\\u<v\\k<l\\r<s}}W_{jm}(W_{ij}+W_{im})W_{vm}(W_{uv}+W_{um})W_{lm}(W_{kl}+W_{km})W_{sm}(W_{rs}+W_{rm}).
\end{equation*}
Taking expectations, we consider 4 types of cases in which the expectation is non-zero:
\begin{itemize}
    \item Case 1: $i=u=k=r$ and $j=v=l=s$ (1 instance),
    \item Case 2: $i=k,u=r$ with $i\neq u$ and $j=l,v=s$  with $j\neq v$ (3 instances),
    \item Case 3: $i=u=k=r$ and $j=l,v=s$ with $j\neq v$ (3 instances),
    \item Case 4: $i=k,u=r$ with $i\neq u$ and $j=v=l=s$ (3 instances),
    \item  Other cases: $\E[W_{jm}(W_{ij}+W_{im})W_{vm}(W_{uv}+W_{um})W_{lm}(W_{kl}+W_{km})W_{sm}(W_{rs}+W_{rm}]=0$.
\end{itemize}
It follows that 
\begin{align*}
    \E[C_n^4X_{n,m}^4]=& 16\biggl[\sum_{i<j}\E[W_{jm}^4]\E[(W_{ij}+W_{im})^4]\cr
    & +3\sum_{\substack{i<j,u<v\\i\neq u, j\neq v}}\E[W_{jm}^2]\,\E[(W_{ij}+W_{im})^2]\,\E[W_{vm}^2]\,\E[(W_{ij}+W_{im})^2]\\
    &+3\sum_{\substack{i<j,v\\j\neq v}}\E[W_{jm}^2]\,\E[W_{vm}^2]\,\E[(W_{ij}+W_{im})^2(W_{iv}+W_{im})^2]\\
    &+3\sum_{\substack{i,u<j\\i\neq u}}\E[(W_{ij}+W_{im})^2]\,\E[(W_{uj}+W_{um})^2]\, \E[W_{jm}^4]\biggr].
\end{align*}
We provide upper bounds for the above expectations. Indeed for all $(a,b)\in\llbracket1,n\rrbracket^2$
\[
    \E[W_{ab}^4] = (1-\alpha_n)^4\alpha_n+\alpha_n^4(1-\alpha_n)=\alpha_n(1-\alpha_n)(\alpha_n^3+(1-\alpha_n)^3).
\]
It is then straightforward to show, taking $c_*>0$ to be a high enough constant, that
\begin{align*}
    & \E[W_{jm}^4]\leq c_*\alpha_n,\qquad 
    \E[(W_{ij}+W_{im})^4] \leq  c_*\alpha_n, \qquad
    \E[W_{jm}^2]^2\leq c_*\alpha_n^2,\cr
    & \E[(W_{ij}+W_{im})^2]^2 \leq c_*\alpha_n^2, \qquad
    \E[(W_{ij}+W_{im})^2(W_{iv}+W_{im})^2] \leq c_*\alpha_n. 
\end{align*}
It follows that
\begin{align} \label{chi2-proof-11}
    \E[C_n^4X_{n,m}^4] &\leq 16\left(\sum_{i<j}c_*^2\alpha_n^2+3\sum_{\substack{i<j,u<v\\i\neq u, j\neq v}}c_*^2\alpha_n^4+3\sum_{\substack{i<j,v\\j\neq v}}c_*^2\alpha_n^3+3\sum_{\substack{i,u<j\\i\neq u}}c_*^2\alpha_n^3\right)\cr
    &\leq 16c_*^2n^2\alpha_n^2\left(1+3n^2\alpha_n^2+6n\alpha_n\right)\cr
    &=O(n^4\alpha_n^4),
\end{align}
where in the last line we have used the assumption of  $n\alpha_n\to\infty$ to identify the dominating term. 
Note that $C_n$ is at the order of $n\sqrt{n}\alpha_n$. 
We thus obtain
\begin{align*}
    \E\left[\sum_{m=1}^n\E[X_{n,m}^4|\mathcal{F}_{n,m-1}]\right] &  = n\cdot O\left(\frac{n^4\alpha_n^4}{(n\sqrt{n}\alpha_n)^4}\right) = O\left(n^{-1}\right). 
\end{align*}
This proves \eqref{chi2-proof-10}.

\subsection{Proof of the null distribution of $\psi_n^{SQ}$} \label{subsec:nullproof-SQ}
We aim to show that
\beq \label{null-proof-goal-SQ}
\varphi_n^{SQ}=\frac{Q_n}{2\sqrt{2}n^2\hat{\alpha}_n^2}\quad \xrightarrow[n\to\infty]{\mathcal{L}}\quad \mathcal{N}(0,1).  
\eeq

Let $\hat{\delta}_n=\alpha_n-\hat{\alpha}_n$. We then have $A_{ij}-\hat{\alpha}_n=W_{ij}+\hat{\delta}_n$. It follows that
$Q_n = \sum_{(i_1,i_2,i_3,i_4) \text{ dist.}} (W_{i_1i_2}+\hat{\delta}_n)(W_{i_2i_3}+\hat{\delta}_n)(W_{i_3i_4}+\hat{\delta}_n)(W_{i_4i_1}+\hat{\delta}_n)$. 
We introduce an ideal version of $Q_n$,
\[
     \tilde{Q}_n = \sum_{(i_1,i_2,i_3,i_4) \text{ dist.}}W_{i_1i_2}W_{i_2i_3}W_{i_3i_4}W_{i_4i_1},
\]
and re-express
\begin{equation} \label{sq-proof-decompose}
    \psi_n^{SQ} = \frac{Q_n}{2\sqrt{2}n^2\hat{\alpha}_n^2} =\frac{Q_n-\tilde{Q}_n}{2\sqrt{2}n^2\alpha_n^2}\left(\frac{\alpha_n}{\hat{\alpha}_n}\right)^2+\frac{\tilde{Q}_n}{2\sqrt{2}n^2\alpha_n^2}\left(\frac{\alpha_n}{\hat{\alpha}_n}\right)^2. 
\end{equation}
If we can show that
\begin{align}
    \text{(a) }& \frac{Q_n-\tilde{Q}_n}{2\sqrt{2}n^2\alpha_n^2}\xrightarrow{\mathbb{P}}0, \label{eq:diff}\\
    \text{(b) }& \frac{\tilde{Q}_n}{2\sqrt{2}n^2\alpha_n^2}\goto_d\mathcal{N}(0,1), \label{eq:ideal}
\end{align}
then \eqref{null-proof-goal-SQ} follows from Slutsky's theorem and the fact that $\hat{\alpha}_n/\alpha_n\xrightarrow{\mathbb{P}}1$. 

What remains is to prove \eqref{eq:diff} and \eqref{eq:ideal}.

\smallskip\noindent
{\bf Proof of Equation \eqref{eq:diff}:} 
Expanding $Q_n$, we obtain:
\begin{align*}
    Q_n-\tilde{Q}_n=&n(n-1)(n-2)(n-3)\hat{\delta}_n^4 +4(n-2)(n-3)\hat{\delta}_n^3\sum_{i\neq j}W_{ij}\\
    &+4(n-3)\hat{\delta}_n^2\sum_{i,j,k \text{ dist.}}W_{ij}W_{jk}+2\hat{\delta}_n^2\sum_{i,j,k,l \text{ dist.}} W_{ij}W_{kl}\\
    &+4\hat{\delta}_n\sum_{i,j,k,l \text{ dist.}} W_{ij}W_{jk}W_{kl}. 
\end{align*}
It follows that
\begin{align} \label{sq-proof-1}
    \biggl|\frac{Q_n-\tilde{Q}_n}{n^2\alpha_n^2}\biggr|\leq&n^2\frac{\hat{\delta}_n^4}{\alpha_n^2}+4\frac{|\hat{\delta}_n|^3}{\alpha_n^2}\biggl|\sum_{i\neq j}W_{ij}\biggr| + \frac{4|\hat{\delta}_n|}{n^2\alpha_n^2}\biggl|\sum_{i,j,k,l \text{ dist.}} W_{ij}W_{jk}W_{kl}\biggr|\cr
    & +\frac{\hat{\delta}_n^2}{\alpha_n^2}\left(\frac{4}{n}\biggl|\sum_{i,j,k \text{ dist.}}W_{ij}W_{jk}\biggr|+\frac{2}{n^2}\biggl|\sum_{i,j,k,l \text{ dist.}} W_{ij}W_{kl}\biggr|\right). 
\end{align}
We will bound each of the terms on the right hand side of \eqref{sq-proof-1}.

Consider the first term in \eqref{sq-proof-1}. Note that
\begin{equation*}
    n^2\frac{\hat{\delta}_n^4}{\alpha_n^2}=4\frac{(1-\alpha_n)^2}{(n-1)^2} \left(\sqrt{\frac{n(n-1)}{2}}\frac{\hat{\alpha}_n-\alpha_n}{\sqrt{\alpha_n(1-\alpha_n)}}\right)^4. 
\end{equation*}
By Central Limit Theorem, $
    \sqrt{\frac{n(n-1)}{2}}\frac{\hat{\alpha}_n-\alpha_n}{\sqrt{\alpha_n(1-\alpha_n)}} \goto \mathcal{N}(0,1)$.
It follows from Slutsky's theorem that
\begin{equation} \label{eq:diff1}
    n^2 \hat{\delta}_n^4/\alpha_n^2\xrightarrow{\mathbb{P}}0.
\end{equation}

Consider the second term in \eqref{sq-proof-1}. Since $\hat{\delta}_n=\hat{\alpha}_n-\alpha_n$, using the definition of $\hat{\alpha}_n$, we immediately have $\sum_{i<j}W_{ij}=\frac{n(n-1)}{2}\hat{\delta}_n$. As a result, 
\begin{equation}\label{eq:diff2}
    \frac{|\hat{\delta}_n^3|}{\alpha_n^2}\biggl|\sum_{i\neq j}W_{ij}\biggr|=n(n-1)\frac{\hat{\delta}_n^4}{\alpha_n^2}\leq n^2\frac{\hat{\delta}_n^4}{\alpha_n^2}\xrightarrow{\mathbb{P}}0.
\end{equation}

Consider the fourth term in \eqref{sq-proof-1}. First, let $A_{n} = \frac{1}{n^3\alpha_n}\sum_{i,j,k \text{ dist.}}W_{ij}W_{ik}$. 
Applying Chebyshev's inequality, we have that for any $\lambda>0$,
\begin{align*}
    \mathbb{P}(|A_{n}|>\lambda)&\leq \frac{\E[A_{n}^2]}{\lambda^2}\leq\frac{6^2}{n^6\alpha_n^2\lambda^2}\sum_{\substack{i<j<k \\ u<v<w}}\E[W_{ij}W_{jk}W_{uv}W_{vw}]\\
    &= \frac{36}{n^6\alpha_n^2\lambda^2}\sum_{i<j<k}\E[W_{ij}^2W_{jk}^2] \leq \frac{36}{n^3\lambda^2} \xrightarrow{n\to\infty}0, 
\end{align*}
which shows that $A_{n}\xrightarrow{\mathbb{P}}0$. Furthermore, 
\begin{align*}
    \frac{\hat{\delta}_n^2}{n\alpha_n^2}\biggl|\sum_{i,j,k \text{ dist.}}W_{ij}W_{jk}\biggr| &=2(1-\alpha_n)\Bigl(\frac{n}{n-1}\Bigr)\biggl[\sqrt{\frac{n(n-1)}{2}}\frac{\hat{\alpha}_n-\alpha_n}{\sqrt{\alpha_n(1-\alpha_n)}}\biggr]^2|A_{n}|.  
\end{align*}
By Slutsky's theorem, we have
\begin{equation}\label{eq:diff3}
    \frac{\hat{\delta}_n^2}{n\alpha_n^2}\biggl|\sum_{i,j,k \text{ dist.}}W_{ij}W_{jk}\biggr|\xrightarrow{\mathbb{P}}0.
\end{equation}
Second, let $B_{n}=\frac{1}{\alpha_nn^4}\sum_{i,j,k,l \text{ dist.}} W_{ij}W_{kl}$. 
We apply Chebyshev's inequality: For any $\lambda>0$,
\begin{align*}
    \mathbb{P}(|B_{n}|>\lambda) &\leq \frac{\E[B_{n}^2]}{\lambda^2}\leq \frac{1}{\alpha_n^2n^8\lambda^2}\sum_{\substack{i,j,k,l \text{ dist.}\\s,t,u,v \text{ dist.}}} \E[W_{ij}W_{kl}W_{st}W_{uv}]\\
    &= \frac{24^2}{\alpha_n^2n^8\lambda^2}\sum_{\substack{i<j<k<l \\s<t<u<v}} \E[W_{ij}W_{kl}W_{st}W_{uv}]= \frac{24^2}{\alpha_n^2n^8\lambda^2}\sum_{i<j<k<l} \E[W_{ij}^2W_{kl}^2]\\
    &=\frac{24^2}{\alpha_n^2n^8\lambda^2}\sum_{i<j<k<l} \E[W_{ij}^2]\E[W_{kl}^2]\leq \frac{24^2}{n^4\lambda^2} \xrightarrow{n\to\infty}0, 
\end{align*}
which shows that $B_{n}\xrightarrow{\mathbb{P}}0$. Furthermore, 
\begin{align*}
    \frac{\hat{\delta}_n^2}{n^2\alpha_n^2}\biggl|\sum_{i,j,k,l \text{ dist.}} W_{ij}W_{kl}\biggr| &=\frac{2(1-\alpha_n)n}{n-1}\biggl[\sqrt{\frac{n(n-1)}{2}}\frac{\hat{\alpha}_n-\alpha_n}{\sqrt{\alpha_n(1-\alpha_n)}}\biggr]^2|B_{n}|. 
\end{align*}
We conclude by Slutsky's theorem that
\begin{equation}\label{eq:diff4}
    \frac{\hat{\delta}_n^2}{n^2\alpha_n^2}\biggl|\sum_{i,j,k,l \text{ dist.}} W_{ij}W_{kl}\biggr|\xrightarrow{\mathbb{P}}0. 
\end{equation}

Consider the third term in \eqref{sq-proof-1}. Write $D_n=\frac{1}{\alpha_n^{3/2}n^3}\sum_{i,j,k,l \text{ dist.}} W_{ij}W_{jk}W_{kl}$. 
By Chebyshev's inequality, for any $\lambda>0$, 
\begin{align*}
    \mathbb{P}(|D_{n}|>\lambda) &\leq \frac{\E[D_n^2]}{\lambda^2}=\frac{1}{\alpha_n^3n^6\lambda^2}\E\Biggl[\sum_{\substack{i,j,k,l \text{ dist.}\\u,v,w,z \text{ dist.}}}W_{ij}W_{jk}W_{kl}W_{uv}W_{vw}W_{wz}\Biggr]\\
    &= \frac{2}{\alpha_n^3n^6\lambda^2}\E\Biggl[\sum_{i,j,k,l \text{ dist.}}W_{ij}^2W_{jk}^2W_{kl}^2\Biggr] \leq \frac{2}{n^2\lambda^2}\xrightarrow[n\to\infty]{} 0,
\end{align*}
which implies that $D_{n}\xrightarrow{\mathbb{P}}0$. Furthermore, 
\begin{align*}
    \frac{\hat{\delta}_n}{n^2\alpha_n^2}\biggl|\sum_{i,j,k,l \text{ dist.}} W_{ij}W_{jk}W_{kl}\biggr| &=  \sqrt{\frac{2n(1-\alpha_n)}{n-1}}\biggl[\sqrt{\frac{n(n-1)}{2}}\frac{\hat{\alpha}_n-\alpha_n}{\sqrt{\alpha_n(1-\alpha_n)}}\biggr]|D_{n}|. 
\end{align*}
We conclude by Slutsky's theorem that
\begin{equation}\label{eq:diff5}
    \frac{\hat{\delta}_n}{n^2\alpha_n^2}\left|\sum_{i,j,k,l \text{ dist.}} W_{ij}W_{jk}W_{kl}\right|\xrightarrow{\mathbb{P}}0. 
\end{equation}

We plug \eqref{eq:diff1}-\eqref{eq:diff5} into \eqref{sq-proof-1} to get \eqref{eq:diff}.


\medskip\noindent
{\bf Proof of Equation \eqref{eq:ideal}:} 
We introduce some notation to simplify the computations. Given 4 distinct nodes, there are 3 different possible cycles, denoted as 
\[
CC(i_1,i_2,i_3,i_4)=\{(i_1,i_2,i_3,i_4),(i_1,i_2,i_4,i_3),(i_1,i_3,i_2,i_4)\}.
\]
Moreover, for $B\subset\{1,2,...,n\}^4$, let $CC(B)=\cup_{(i_1,i_2,i_3,i_4)\in B}CC(i_1,i_2,i_3,i_4)$. For $1\leq m\leq n$, let $I_m$ be the collection of $(i_1,i_2,i_3,i_4)$ such that $1\leq i_1 < i_2 < i_3 < i_4 \leq m$. We thus have
\begin{equation}\label{eq:cycdecomp}
    \tilde{Q}_n = 8\sum_{CC(I_n)}W_{i_1i_2}W_{i_2i_3}W_{i_3i_4}W_{i_4i_1}. 
\end{equation}
It is straightforward to see that $\E[\tilde{Q}_n]=0$. In addition, notice that the terms in the sum are uncorrelated, since they all correspond to different cycles: to obtain a non-zero correlation between $W_{i_1i_2}W_{i_2i_3}W_{i_3i_4}W_{i_4i_1}$ and $W_{i_1'i_2'}W_{i_2'i_3'}W_{i_3'i_4'}W_{i_4'i_1'}$, we would need to uniquely match each factor in $W_{i_1i_2}W_{i_2i_3}W_{i_3i_4}W_{i_4i_1}$ with a factor in $W_{i_1'i_2'}W_{i_2'i_3'}W_{i_3'i_4'}W_{i_4'i_1'}$, which is equivalent to overlaying the two cycles $[i_1i_2i_3i_4]$ and $[i'_1i_2'i_3'i_4']$. Let's compute the variance
\begin{align*}
    \V(\tilde{Q}_n) &= 64\V\left(\sum_{CC(I_n)}W_{i_1i_2}W_{i_2i_3}W_{i_3i_4}W_{i_4i_1}\right)\\
    &= 64\alpha_n^4(1-\alpha_n)^4\times3{n \choose 4}= 8\alpha_n^4(1-\alpha_n)^4n(n-1)(n-2)(n-3). 
\end{align*}
Let $Z_n:=2\sqrt{2n(n-1)(n-2)(n-3)}\alpha_n^2(1-\alpha_n)^2$. It is easy to see that 
$n^2\alpha_n^2/Z_n\xrightarrow{n\to\infty}1$. By Slutsky's theorem, to show \eqref{eq:ideal}, it suffices to show that 
\begin{equation}\label{eq:clt}
    \frac{\tilde{Q}_n}{Z_n}\xrightarrow[n\to\infty]{\mathcal{L}}\mathcal{N}(0,1). 
\end{equation}

We now prove \eqref{eq:clt}. 
For each $1\leq m\leq n$, we define
\[
X_{n,m}=\frac{\tilde{Q}_{n,m}-\tilde{Q}_{n,m-1}}{Z_n}, \qquad\mbox{where}\quad \tilde{Q}_{n,m}=\sum_{CC(I_m)}W_{i_1i_2}W_{i_2i_3}W_{i_3i_4}W_{i_4i_1}.
\]
By default, we let $\tilde{Q}_{n,0}=1$. Recall that we previously defined the filtration $\{\mathcal{F}_{n,m}:0\leq m\leq n\}$ such that $\mathcal{F}_{n,m}=\sigma\{W_{ij}:(i,j)\in\llbracket1,m \rrbracket^2\}$ for $m\geq 1$ and  $\mathcal{F}_{n,0}=\{\Omega,\emptyset\}$ (where $\Omega$ denotes the sample space).
It is easy to see that $\E[|\tilde{Q}_{n,m}|]<\infty$. Hence, $\tilde{Q}_{n,m}$ is $\mathcal{F}_{n,m}$-measurable. It is also straightforward to show that $\E[\tilde{Q}_{n,m+1}|\mathcal{F}_{n,m}]=\tilde{Q}_{n,m}$. Therefore, the sequence $\{Q_{n,m}:m\in\llbracket1,n \rrbracket\}$ is a martingale with respect to $\{\mathcal{F}_{n,m}:m\in\llbracket1,n \rrbracket\}$. It follows that the sequence $\{X_{n,m}:m\in\llbracket1,n \rrbracket\}$ is a martingale difference sequence. 
Note that 
\[
\tilde{Q}_n/Z_n=\tilde{Q}_{n,n}/Z_n=\sum_{m=1}^n X_{n,m}. 
\]
By the martingale Central Limit Theorem, to show \eqref{eq:clt}, it suffices to show: 
\begin{align}
    \text{(b1) }& \sum_{m=1}^n\mathbb{E}[X_{n,m}^2|\mathcal{F}_{n,m-1}]\xrightarrow{\mathbb{P}}1 \label{eq:clt1},\\
    \text{(b2) }& \forall\epsilon>0, \sum_{m=1}^n\mathbb{E}[X_{n,m}^2\bm{1}\{|X_{n,m}>\epsilon|\}|\mathcal{F}_{n,m-1}]\xrightarrow{\mathbb{P}}0. \label{eq:clt2}
\end{align}
Below, we show \eqref{eq:clt1} and \eqref{eq:clt2} separately.

In the first part, we prove \eqref{eq:clt1}. It suffices to show:
\beq \label{eq:clt11}
\E\left[\sum_{m=1}^n\E[X_{n,m}^2|\mathcal{F}_{n,m-1}]\right]=1,
\eeq 
and 
\beq  \label{eq:clt12}
\text{Var}\left(\sum_{m=1}^n\E[X_{n,m}^2|\mathcal{F}_{n,m-1}]\right)\xrightarrow{n\to\infty}0. 
\eeq

Consider \eqref{eq:clt11} first. Recall that by definition
\begin{equation*}
    X_{n,m}=\frac{\tilde{Q}_{n,m}-\tilde{Q}_{n,m-1}}{Z_n}=\frac{8}{Z_n}\sum_{CC(I_m)\setminus CC(I_{m-1})}W_{i_1i_2}W_{i_2i_3}W_{i_3i_4}W_{i_4i_1}. 
\end{equation*}
An alternative way to enumerate all cycles in $CC(I_m)\setminus CC(I_{m-1})$ is to  first select a set of two indices $\{i,j\}$ (we take, \textit{wlog}, $i<j$) from $\{1,...,m-1\}$ and use them as the neighboring nodes of $m$ in the cycle. Then select $k\in\{1,...,m-1\}\setminus\{i,j\}$ as the last node of the cycle.
\begin{equation*}
    X_{n,m}=\frac{8}{Z_n}\sum_{1\leq i<j\leq m-1}W_{mi}W_{mj}Y_{m-1,ij}, 
\qquad
\mbox{where}\quad 
    Y_{m-1,ij} = \sum_{\substack{1\leq k\leq m-1\\k\notin\{i,j\}}}W_{ki}W_{kj}. 
\end{equation*}
It follows that
\begin{align*}
    \E[X_{n,m}^2|\mathcal{F}_{n,m-1}]&=\frac{64}{Z_n^2}\sum_{\substack{1\leq i<j\leq m-1\\1\leq u<v\leq m-1}}\E[W_{mi}W_{mj}Y_{m-1,ij}W_{mu}W_{mv}Y_{m-1,uv}|\mathcal{F}_{n,m-1}]\\
    &=\frac{64}{Z_n^2}\sum_{\substack{1\leq i<j\leq m-1\\1\leq u<v\leq m-1}}Y_{m-1,ij}Y_{m-1,uv}\E[W_{mi}W_{mj}W_{mu}W_{mv}]\\
    &=\frac{64}{Z_n^2}\sum_{1\leq i<j\leq m-1}Y_{m-1,ij}^2\E[W_{mi}^2W_{mj}^2] = \frac{64\alpha_n^2(1-\alpha_n)^2}{Z_n^2}\sum_{1\leq i<j\leq m-1}Y_{m-1,ij}^2.
\end{align*}
Hence
\begin{align*}
    \E\left[\sum_{m=1}^n\E[X_{n,m}^2|\mathcal{F}_{n,m-1}]\right] &= \frac{64\alpha_n^2(1-\alpha_n)^2}{Z_n^2}\sum_{m=1}^n\sum_{1\leq i<j\leq m-1}\E[Y_{m-1,ij}^2],
\end{align*}
where
\begin{align*}
    \E[Y_{m-1,ij}^2]&
    =\sum_{\substack{1\leq k,l\leq m-1\\k,l\notin\{i,j\}}}\E\left[W_{ki}W_{kj}W_{li}W_{lj}\right]=\sum_{\substack{1\leq k\leq m-1\\k\notin\{i,j\}}}\E\left[W_{ki}^2W_{kj}^2\right]\\
    &= (m-3)\alpha_n^2(1-\alpha_n)^2. 
\end{align*}
It follows that
\begin{align*}
    \E\left[\sum_{m=1}^n\E[X_{n,m}^2|\mathcal{F}_{n,m-1}]\right] &= \frac{64\alpha_n^2(1-\alpha_n)^2}{Z_n^2}\sum_{m=1}^n\frac{(m-1)(m-2)(m-3)}{2}\alpha_n^2(1-\alpha_n)^2=1. 
\end{align*}
This proves \eqref{eq:clt11}.

Consider \eqref{eq:clt12} next. We decompose $\sum_{m=1}^n\E[X_{n,m}^2|\mathcal{F}_{n,m-1}]$ into a sum of two components, then calculate its variance. Note that
\begin{equation*}
    Y_{m-1,ij}^2 = \Biggl(\sum_{\substack{1\leq k\leq m-1\\k\notin\{i,j\}}}W_{ki}W_{kj}\Biggr)^2 = \sum_{\substack{1\leq k\leq m-1\\k\notin\{i,j\}}}W_{ki}^2W_{kj}^2+2\sum_{\substack{1\leq k< l\leq m-1\\k,l\notin\{i,j\}}}W_{ki}W_{kj}W_{li}W_{lj}. 
\end{equation*}
Hence
\begin{align*}
    \sum_{m=1}^n\E[X_{n,m}^2|\mathcal{F}_{n,m-1}] &= \frac{64\alpha_n^2(1-\alpha_n)^2}{Z_n^2}\sum_{m=1}^n\sum_{1\leq i<j\leq m-1}Y_{m-1,ij}^2= \frac{16n^4\alpha_n^4}{Z_n^2}(I_a+I_b),
\end{align*}
where we denote
\begin{align*}
    I_a &= \frac{4(1-\alpha_n)^2}{n^4\alpha_n^2}\sum_{m=1}^n\sum_{1\leq i<j\leq m-1}\sum_{\substack{1\leq k\leq m-1\\k\notin\{i,j\}}}W_{ki}^2W_{kj}^2,\\
    I_b &= \frac{8(1-\alpha_n)^2}{n^4\alpha_n^2}\sum_{m=1}^n\sum_{1\leq i<j\leq m-1}\sum_{\substack{1\leq k< l\leq m-1\\k,l\notin\{i,j\}}}W_{ki}W_{kj}W_{li}W_{lj}. 
\end{align*}
Using the Cauchy-Schwarz inequality we obtain
\begin{align*}
    \text{Var}\left(\sum_{m=1}^n\E[X_{n,m}^2|\mathcal{F}_{n,m-1}]\right)&=\frac{256n^8\alpha_n^8}{Z_n^4}(\text{Var}(I_a)+\text{Var}(I_b)+2\text{Cov}(I_a,I_b))\\
    &\leq\frac{256n^8\alpha_n^8}{Z_n^4}( \sqrt{\text{Var}(I_a)}+\sqrt{\text{Var}(I_b)})^2.
\end{align*}
Hence, it suffices to show that $\text{Var}(I_a)\xrightarrow{n\to\infty}0$ and $\text{Var}(I_b)\xrightarrow{n\to\infty}0$ separately. For $\text{Var}(I_a)$, we first rearrange the sums in the expression of $I_a$
\begin{align*}
    I_a &= \frac{4(1-\alpha_n)^2}{n^4\alpha_n^2}\sum_{k=1}^{n-1}\sum_{\substack{1\leq i<j\leq n-1\\i,j\neq k}}\sum_{m>\max\{i,j,k\}}W_{ki}^2W_{kj}^2\\
    &=\frac{4(1-\alpha_n)^2}{n^4\alpha_n^2}\sum_{k=1}^{n-1}\sum_{\substack{1\leq i<j\leq n-1\\i,j\neq k}}(n-\max\{i,j,k\}+1)W_{ki}^2W_{kj}^2.
\end{align*}
Note that the terms of the first sum over $k=1,...,n$ are pairwise independent, which will facilitate variance computations. Hence
\begin{align*}
    \text{Var}(I_a) &= \frac{16(1-\alpha_n)^4}{n^8\alpha_n^4}\sum_{k=1}^{n-1} \text{Var}\left(\sum_{\substack{1\leq i<j\leq n-1\\i,j\neq k}}(n-\max\{i,j,k\}+1)W_{ki}^2W_{kj}^2\right)\\
    &\leq\frac{16(1-\alpha_n)^4}{n^6\alpha_n^4}\sum_{k=1}^{n-1}\sum_{\substack{1\leq i<j\leq n-1\\i,j\neq k}}\sum_{\substack{1\leq u<v\leq n-1\\u,v\neq k}}\text{Cov}(W_{ki}^2W_{kj}^2,W_{ku}^2W_{kv}^2).
\end{align*}
We can consider four cases for $\text{Cov}(W_{ki}^2W_{kj}^2,W_{ku}^2W_{kv}^2)$:
\begin{enumerate}
    \item $(i,j)=(u,v)$, then $\text{Var}(W_{ki}^2W_{kj}^2)\leq\E[W_{ki}^4W_{kj}^4]=\E[W_{ki}^4]^2\leq c\alpha_n^2$,
    \item $i=u, j\neq v$, then $\text{Cov}(W_{ki}^2W_{kj}^2,W_{ki}^2W_{kv}^2)\leq \E[W_{ki}^4W_{kj}^2W_{kv}^2]=\E[W_{ki}^4]\E[W_{kj}^2]^2\leq c\alpha_n^3$,
    \item The previous bound will also hold for the case $i\neq u, j= v$, the case $i=v$, and the case $j=u$,
    \item For any other case, $\text{Cov}(W_{ki}^2W_{kj}^2,W_{ki}^2W_{kv}^2)=0$.
\end{enumerate}
Here, $c>0$ is a high enough constant. It follows that 
\begin{align*}
    \text{Var}(I_a)&=\frac{16(1-\alpha_n)^4}{n^8\alpha_n^4}\sum_{k=1}^{n-1}\sum_{\substack{1\leq i<j\leq n-1\\i,j\neq k}}\left\{\text{Var}(W_{ki}^2W_{kj}^2) + \sum_{\substack{v=i+1\\v\notin\{k,j\}}}^{n-1}\text{Cov}(W_{ki}^2W_{kj}^2,W_{ki}^2W_{kv}^2) \right.\\
    &+\left. \sum_{\substack{u=1\\u\notin\{k,i\}}}^{j-1}\text{Cov}(W_{ki}^2W_{kj}^2,W_{ku}^2W_{kj}^2) + \sum_{\substack{u=1\\u\neq k}}^{i-1}\text{Cov}(W_{ki}^2W_{kj}^2,W_{ku}^2W_{ki}^2)
    + \sum_{\substack{v=j+1\\v\neq k}}^{n-1}\text{Cov}(W_{ki}^2W_{kj}^2,W_{kj}^2W_{kv}^2)\right\}\\
    &\leq \frac{16c(1-\alpha_n)^4}{n^8\alpha_n^4}\sum_{k=1}^{n-1}\sum_{\substack{1\leq i<j\leq n-1\\i,j\neq k}}\left\{\alpha_n^2+4n\alpha_n^3\right\} \leq \frac{8c}{n^3(n\alpha_n)}\left(4+\frac{1}{n\alpha_n}\right)\xrightarrow{n\to\infty}0. 
\end{align*}
Let's now show that $\text{Var}(I_b)\xrightarrow{n\to\infty}0$. Recall that
\begin{align*}
     I_b &= \frac{8(1-\alpha_n)^2}{n^4\alpha_n^2}\sum_{m=1}^n\sum_{1\leq i<j\leq m-1}\sum_{\substack{1\leq k< l\leq m-1\\k,l\notin\{i,j\}}}W_{ki}W_{kj}W_{li}W_{lj}\\
     &=\frac{2(1-\alpha_n)^2}{n^4\alpha_n^2}\sum_{m=1}^n\sum_{\substack{1\leq i,j,k,l\leq m-1\\i,j,k,l \text{ dist.}}}W_{ki}W_{kj}W_{li}W_{lj}\\
     &= \frac{2(1-\alpha_n)^2}{n^4\alpha_n^2}\sum_{\substack{1\leq i,j,k,l\leq n-1\\i,j,k,l \text{ dist.}}}\sum_{m>\max\{i,j,k,l\}}W_{ki}W_{kj}W_{li}W_{lj}\\
     &= \frac{2(1-\alpha_n)^2}{n^4\alpha_n^2}\sum_{\substack{1\leq i,j,k,l\leq n-1\\i,j,k,l \text{ dist.}}}(n+1-\max\{i,j,k,l\})W_{ki}W_{kj}W_{li}W_{lj}. 
\end{align*}
Therefore, 
\begin{align*}
    \text{Var}(I_b) &= \frac{4(1-\alpha_n)^4}{n^8\alpha_n^4}\text{Var}\left(\sum_{\substack{1\leq i,j,k,l\leq n-1\\i,j,k,l \text{ dist.}}}(n+1-\max\{i,j,k,l\})W_{ik}W_{kj}W_{jl}W_{li}\right)\\
     &= \frac{4(1-\alpha_n)^4}{n^8\alpha_n^4}\text{Var}\left(8\sum_{CC(I_{n-1})}(n+1-\max\{i,j,k,l\})W_{ik}W_{kj}W_{jl}W_{li}\right)\\
    &= \frac{32(1-\alpha_n)^4}{n^8\alpha_n^4}\sum_{\substack{1\leq i,j,k,l\leq n-1\\i,j,k,l \text{ dist.}}}(n+1-\max\{i,j,k,l\})^2\text{Var}(W_{ik}W_{kj}W_{jl}W_{li})\\
    &\leq \frac{32(1-\alpha_n)^4}{n^6\alpha_n^4}\sum_{\substack{1\leq i,j,k,l\leq n-1\\i,j,k,l \text{ dist.}}}\alpha_n^4(1-\alpha_n)^4 \leq \frac{32}{n^2}\xrightarrow{n\to\infty}0.
\end{align*}
This gives $\text{Var}(I_b)\xrightarrow{n\to\infty}0$. Recall that we had:
\begin{equation*}
    \text{Var}\left(\sum_{m=1}^n\E[X_{n,m}^2|\mathcal{F}_{n,m-1}]\right) \leq\frac{256n^8\alpha_n^8}{Z_n^4}\left( \sqrt{\text{Var}(I_a)}+\sqrt{\text{Var}(I_b)}\right)^2. 
\end{equation*}
Since $\frac{256n^8\alpha_n^8}{Z_n^4}\xrightarrow{n\to\infty}4$, we obtain \eqref{eq:clt12}. In combination with \eqref{eq:clt11}, this proves \eqref{eq:clt1}.

\bigskip

In the second part, we prove \eqref{eq:clt2}. We have, using the Cauchy-Schwarz and Markov inequalities
\begin{align*}
    \sum_{m=1}^n\mathbb{E}[X_{n,m}^2\bm{1}\{|X_{n,m}>\epsilon|\}|\mathcal{F}_{n,m-1}] &\leq \sum_{m=1}^n\sqrt{\E[X_{n,m}^4|\mathcal{F}_{n,m-1}]}\sqrt{\mathbb{P}(|X_{n,m}|\geq \epsilon|\mathcal{F}_{n,m-1}])}\\
    &\leq \frac{1}{\epsilon^2}\sum_{m=1}^n\E[X_{n,m}^4|\mathcal{F}_{n,m-1}]. 
\end{align*}
Hence it suffices to show that
\begin{equation} \label{eq:clt21}
    \E\left[\sum_{m=1}^n\E[X_{n,m}^4|\mathcal{F}_{n,m-1}]\right]=\sum_{m=1}^n\E[X_{n,m}^4]\xrightarrow{n\to\infty}0. 
\end{equation}
Recall that for all $n\in\mathbb{N}^*$, for all $m\in\llbracket1,n\rrbracket$
\begin{equation*}
    X_{n,m} = \frac{2}{Z_n}\sum_{1\leq i<j\leq m-1}W_{mi}W_{mj}Y_{m-1,ij} \quad \text{ with } \quad Y_{m-1,ij}=\sum_{\substack{1\leq k\leq m-1\\k\notin\{i,j\}}}W_{ki}W_{kj}. 
\end{equation*}
It follows that 
\begin{align*}
    & \E[X_{n,m}^4|\mathcal{F}_{n,m-1}]\cr
    =& \frac{16}{Z_n^4}\sum_{\substack{i<j,u<v\\k<l,r<s}}Y_{m-1,ij}Y_{m-1,uv}Y_{m-1,kl}Y_{m-1,rs} \times \E[W_{mi}W_{mj}W_{mu}W_{mv}W_{mk}W_{ml}W_{mr}W_{ms}]\\
    =& \frac{16}{Z_n^4} \left\{ \sum_{i<j} Y_{m-1,ij}^4\E[W_{mi}^4W_{mj}^4] + 3\sum_{i}\sum_{\substack{j,v\\j,v>i \text{ and }j\neq v}}Y_{m-1,ij}^2Y_{m-1,iv}^2\E[W_{mi}^4W_{mj}^2W_{mv}^2] \right.\\
    &\left.+3\sum_{j}\sum_{\substack{i,u\\i,u<j \text{ and }i\neq u}}Y_{m-1,ij}^2Y_{m-1,uj}^2\E[W_{mj}^4W_{mi}^2W_{mu}^2]+9\sum_{i<j,u<v}Y_{m-1,ij}^2Y_{m-1,uv}^2\E[W_{mi}^2W_{mj}^2W_{mu}^2W_{mv}^2]\right\}\\
    \leq&  \frac{16}{Z_n^4} \left\{ \sum_{i<j} Y_{m-1,ij}^4c\alpha_n^2 + 3\sum_{i}\sum_{\substack{j,v\\j,v>i \text{ and }j\neq v}}Y_{m-1,ij}^2Y_{m-1,iv}^2c\alpha_n^3\right.\\
    &\left.+3\sum_{j}\sum_{\substack{i,u\\i,u<j \text{ and }i\neq u}}Y_{m-1,ij}^2Y_{m-1,uj}^2c\alpha_n^3+9\sum_{i<j,u<v}Y_{m-1,ij}^2Y_{m-1,uv}^2c\alpha_n^4\right\},
\end{align*}
where $c>0$ is a high enough constant. Hence,
\begin{align*}
     \E[X_{n,m}^4] &\leq \frac{16c}{Z_n^4}\left\{\alpha_n^2\sum_{i<j} \E[Y_{m-1,ij}^4] + 3\alpha_n^3\sum_{i}\sum_{\substack{j,v\\j,v>i \text{ and }j\neq v}}\E[Y_{m-1,ij}^2Y_{m-1,iv}^2]\right.\\
     & \left. +3\alpha_n^3\sum_{j}\sum_{\substack{i,u\\i,u<j \text{ and }i\neq u}}\E[Y_{m-1,ij}^2Y_{m-1,uj}^2]+9\alpha_n^4\sum_{i<j,u<v}\E[Y_{m-1,ij}^2]\E[Y_{m-1,uv}^2]\right\}.
\end{align*}
We will now compute upper bounds on $\E[Y_{m-1,ij}^4]$, $\E[Y_{m-1,ij}^2Y_{m-1,iv}^2]$ and $\E[Y_{m-1,ij}^2]$. We have
\begin{align*}
    \E[Y_{m-1,ij}^4] &= \E\left[\sum_{k,l,u,v\notin\{i,j\}}W_{ki}W_{kj}W_{li}W_{lj}W_{ui}W_{uj}W_{vi}W_{vj}\right]\\
    &= 3\sum_{k,u\notin\{i,j\}}\E[W_{ki}^2W_{kj}^2W_{ui}^2W_{uj}^2]\\
    &= 3\left(\sum_{k\notin\{i,j\}}\E[W_{ki}^4W_{kj}^4]+\sum_{k\neq u; \; k,u\notin\{i,j\}}\E[W_{ki}^2W_{kj}^2W_{ui}^2W_{uj}^2]\right)\\
    &\leq 12m\alpha_n^2+3m^2\alpha_n^4\leq c_1(m\alpha_n^2+m^2\alpha_n^4), 
\end{align*}
where $c_1>0$ is a constant. Similarly
\begin{align*}
    \E[Y_{m-1,ij}^2]&=\E\left[\sum_{k,l\notin\{i,j\}}W_{ki}W_{kj}W_{li}W_{lj}\right] = \sum_{k\notin\{i,j\}}\E[W_{ki}^2W_{kj}^2]\leq m\alpha_n^2,
\end{align*}
and
\begin{align*}
    \E&[Y_{m-1,ij}^2Y_{m-1,iv}^2] =\E\left[\sum_{k,l,r,s: k,l\notin\{i,j\}, r,s\notin\{i,v\}}W_{ki}W_{kj}W_{li}W_{lj}W_{ri}W_{rv}W_{si}W_{sv}\right]\\
    &= \E\left[\sum_{k,r: k\notin\{i,j\}, r\notin\{i,v\}}W_{ki}^2W_{kj}^2W_{ri}^2W_{rv}^2\right]= \sum_{k,r: k\notin\{i,j\},r\notin\{i,v\}}\E[W_{ki}^2W_{kj}^2W_{ri}^2W_{rv}^2]\\
    &= \sum_{k\notin\{i,j,v\}}\E[W_{ki}^4W_{kj}^2W_{kv}^2]+\sum_{k\neq r; k\notin\{i,j\}, r\notin\{i,v\}}\E[W_{ki}^2W_{kj}^2W_{ri}^2W_{rv}^2]\\
    &\leq 2m\alpha_n^3+m^2\alpha_n^4 \leq c_2m^2\alpha_n^3,
\end{align*}
for $n$ big enough (since $\alpha_n\xrightarrow{n\to\infty}0$), where $c_2>0$ is a constant. It follows that, for some constant $\gamma>\max\{1,c,c_1,c_2\}$, we have
\begin{align*}
    \E[X_{n,m}^4] &\leq \frac{16c}{Z_n^4}\left\{\alpha_n^2\sum_{i<j} \E[Y_{m-1,ij}^4] + 3\alpha_n^3\sum_{i}\sum_{\substack{j,v\\j,v>i \text{ and }j\neq v}}\E[Y_{m-1,ij}^2Y_{m-1,iv}^2]\right.\\
     & \left. +3\alpha_n^3\sum_{j}\sum_{\substack{i,u\\i,u<j \text{ and }i\neq u}}\E[Y_{m-1,ij}^2Y_{m-1,uj}^2]+9\alpha_n^4\sum_{i<j,u<v}\E[Y_{m-1,ij}^2]\E[Y_{m-1,uv}^2]\right\}\\
     &\leq \frac{16\gamma^2}{Z_n^4}(m^3\alpha_n^4+m^4\alpha_n^6+6m^5\alpha_n^6+9\alpha_n^8m^6)\\
     &\leq \frac{16\gamma^2}{Z_n^4}(n^3\alpha_n^4+n^4\alpha_n^6+6n^5\alpha_n^6+9\alpha_n^8n^6). 
\end{align*}
As a result, 
\begin{align*}
    \sum_{m=1}^n \E[X_{n,m}^4] &\leq \frac{16\gamma^2}{n^2(n-1)^2(n-2)^2(n-3)^2\alpha_n^8(1-\alpha_n)^8}(n^4\alpha_n^4+n^5\alpha_n^6+6n^6\alpha_n^6+9\alpha_n^8n^7)\\
    &= \left(\frac{144\gamma^2n^6}{(n-1)^2(n-2)^2(n-3)^2(1-\alpha_n)^8}\right)\left(\frac{1}{n^4\alpha_n^4}+\frac{1}{n^3\alpha_n^2}+\frac{1}{n^2\alpha_n^2}+\frac{1}{n}\right)\xrightarrow{n\to\infty}0.
\end{align*}
This gives \eqref{eq:clt21}. Then, \eqref{eq:clt2} follows immediately. 

\subsection{Proof of the joint null distribution}
We now show the desirable claim \eqref{null-proof-goal}. We shall use the previously defined notations: 
\begin{align*}
T_n&=\sum_{i,j,k\text{ dist.}}(A_{ik}-\alpha_n)(A_{jk}-\alpha_n), \qquad \hat{T}_n=\sum_{i,j,k\text{ dist.}}(A_{ik}-\hat{\alpha}_n)(A_{jk}-\hat{\alpha}_n),\cr
\tilde{Q}_n&=\sum_{(i_1,i_2,i_3,i_4) \text{ dist.}}W_{i_1i_2}W_{i_2i_3}W_{i_3i_4}W_{i_4i_1}.
\end{align*}
We have seen the decomposition of $\psi_n^{DC}$ in \eqref{chi2-proof-1} and the decomposition of $\psi_n^{SQ}$ in \eqref{sq-proof-decompose}. We plug them into the definition of $S_n$ to get: 
\begin{align} \label{proof-joinNull-decompose}
    S_n & = u\left[ \frac{\hat{T}_n}{(n-1)\hat{\alpha}_n(1-\hat{\alpha}_n)}\right] + v\left[  \frac{Q_n-\tilde{Q}_n}{2\sqrt{2}n^2\alpha_n^2}\left(\frac{\alpha_n}{\hat{\alpha}_n}\right)^2+\frac{\tilde{Q}_n}{2\sqrt{2}n^2\alpha_n^2}\left(\frac{\alpha_n}{\hat{\alpha}_n}\right)^2 \right]\cr
    &= \epsilon_n +u\frac{\frac{T_n}{(n-1)\alpha_n(1-\alpha_n)}}{\sqrt{\frac{2n(n-2)}{(n-1)}}} + v\frac{\tilde{Q}_n}{2\sqrt{2}n^2\alpha_n^2},
\end{align}
where 
\[
    \epsilon_n= u\frac{\frac{T_n}{(n-1)\alpha_n(1-\alpha_n)}}{\sqrt{\frac{2n(n-2)}{(n-1)}}}\left[\frac{\sqrt{n-1}\alpha_n(1-\alpha_n)\hat{T}_n}{\sqrt{n-2}\hat{\alpha}_n(1-\hat{\alpha}_n)T_n}-1\right]+ v\left[\frac{\alpha^2_n}{\hat{\alpha}^2_n}\frac{(Q_n-\tilde{Q}_n)}{2\sqrt{2}n^2\alpha_n^2}+\Bigl(\frac{\alpha^2_n}{\hat{\alpha}^2_n}-1\Bigr)\frac{\tilde{Q}_n}{2\sqrt{2}n^2\alpha_n^2}\right]. 
\]
In Sections~\ref{subsec:nullproof-DC}-\ref{subsec:nullproof-SQ}, we have shown that
\beq \label{nullproof-alreadyshown}
\frac{\hat{\alpha}_n}{\alpha_n}\xrightarrow{\mathbb{P}} 1, \qquad \frac{\hat{T}_n}{T_n}\xrightarrow{\mathbb{P}}1, \qquad \frac{\frac{T_n}{(n-1)\alpha_n(1-\alpha_n)}}{\sqrt{\frac{2n(n-2)}{(n-1)}}}\xrightarrow{d}\mathcal{N}(0,1), \qquad \frac{Q_n-\tilde{Q}_n}{2\sqrt{2}n^2\alpha_n^2}\xrightarrow{d}\mathcal{N}(0,1). 
\eeq
It follows immediately that  $\epsilon_n\xrightarrow{\mathbb{P}}0$. 
By Slutsky's theorem, it suffices to show that
\begin{equation} \label{jointNull-proof-key}
C_n \overset{\Delta}{=} u\frac{\frac{T_n}{(n-1)\alpha_n(1-\alpha_n)}}{\sqrt{\frac{2n(n-2)}{(n-1)}}}+v\frac{\tilde{Q}_n}{2\sqrt{2}n^2\alpha_n^2}\quad  \xrightarrow[n\to\infty]{\mathcal{L}} \quad  \mathcal{N}(0,1). 
\end{equation}

Below, we show \eqref{jointNull-proof-key}. In Section~\ref{subsec:nullproof-DC}, we have defined $I_m$ as the collection of all distinct $\{(i,j,k)$ such that $1\leq i,j,k\leq m$; in Section~\ref{subsec:nullproof-SQ}, we have defined $CC(I_m)$. 
For each $1\leq m\leq n$, let  
\[
    T_{n,m}=\sum_{(j_1,j_2,j_3)\in I_m}W_{j_1j_3}W_{j_2j_3}, \qquad \tilde{Q}_{n,m}=\sum_{CC(I_m)}W_{i_1i_2}W_{i_2i_3}W_{i_3i_4}W_{i_4i_1}, 
\]
where $T_{n,0}=\tilde{Q}_{n,0}=0$ by default. Introduce
\begin{equation*}
    C_{n,m}=u\frac{\frac{T_{n,m}}{(n-1)\alpha_n(1-\alpha_n)}}{\sqrt{\frac{2n(n-2)}{(n-1)}}}+v\frac{\tilde{Q}_{n,m}}{2\sqrt{2}n^2\alpha_n^2}, \qquad\mbox{for all }0\leq m\leq n. 
\end{equation*}
We have seen that $\{T_{n,m}\}_{0\leq m\leq n}$ and $\{\tilde{Q}_{n,m}\}_{0\leq m\leq n}$ are both martingales with respect to the filtration $\{\mathcal{F}_{n,m}\}_{0\leq m\leq n}$ defined before. It is easy to see that $\{C_{n,m}\}_{0\leq m\leq n}$ is also a martingale. Write
\[
C_n  = \sum_{m=1}^n D_{n,m}, \qquad\mbox{where}\quad D_{n,m}\equiv C_{n,m}-C_{n,m-1}. 
\]
To show $C_n\xrightarrow{d}\mathcal{N}(0,1)$, we apply the martingale Central Limit Theorem. It suffices to show:
\begin{align}
    \text{(a) }& \sum_{m=1}^n\mathbb{E}[D_{n,m}^2|\mathcal{F}_{n,m-1}]\xrightarrow{\mathbb{P}}1, \label{jointCLT-1}\\
    \text{(b) }& \forall\epsilon>0, \sum_{m=1}^n\mathbb{E}[D_{n,m}^2\bm{1}\{|D_{n,m}>\epsilon|\}|\mathcal{F}_{n,m-1}]\xrightarrow{\mathbb{P}}0. \label{jointCLT-2}
\end{align}

It remains to show \eqref{jointCLT-1}-\eqref{jointCLT-2}. Consider \eqref{jointCLT-2}. Write
\[
D^{(1)}_{n,m}= \frac{\frac{T_{n,m}-T_{n,m-1}}{(n-1)\alpha_n(1-\alpha_n)}}{\sqrt{\frac{2n(n-2)}{(n-1)}}}, \qquad\mbox{and}\qquad D^{(1)}_{n,m}=\frac{\tilde{Q}_{n,m}-\tilde{Q}_{n,m-1}}{2\sqrt{2}n^2\alpha_n^2}. 
\]
Then, $D_{n,m}=uD^{(1)}_{n,m}+vD^{(2)}_{n,m}$. It follows that $D^4_{n,m}\leq 8u^4(D^{(1)}_{n,m})^4+8v^4(D^{(2)}_{n,m})^4$. As a result, for any $\epsilon>0$, 
by the Cauchy-Schwarz inequality and the Markov inequality, we have
\begin{align*}
\biggl(\sum_{m=1}^n\mathbb{E}[D_{n,m}^2\bm{1}\{|D_{n,m}&>\epsilon|\}|\mathcal{F}_{n,m-1}]\biggr)^2 \leq\biggl(\sum_{m=1}^n\mathbb{E}[D_{n,m}^4|\mathcal{F}_{n,m-1}]\biggr)\cdot \mathbb{P}\bigl(|D_{n,m}|>\epsilon|\mathcal{F}_{n,m-1}\bigr)\cr
&\leq \sum_{m=1}^n\mathbb{E}[D_{n,m}^4|\mathcal{F}_{n,m-1}]\cr
    &\leq 8u^4\sum_{m=1}^n\mathbb{E}[(D^{(1)}_{n,m})^4|\mathcal{F}_{n,m-1}]+8v^4\sum_{m=1}^n \mathbb{E}[(D^{(2)}_{n,m})^4|\mathcal{F}_{n,m-1}]. 
\end{align*}
With significant efforts, we have shown $\sum_{m=1}^n\mathbb{E}[(D^{(1)}_{n,m})^4|\mathcal{F}_{n,m-1}]\xrightarrow{\mathbb{P}} 0$ in Section~\ref{subsec:nullproof-DC}, and we have shown $\sum_{m=1}^n\mathbb{E}[(D^{(2)}_{n,m})^4|\mathcal{F}_{n,m-1}]\xrightarrow{\mathbb{P}}$ in Section~\ref{subsec:nullproof-SQ}. Plugging them into the above inequality, we immediately obtain \eqref{jointCLT-2}. 

Consider \eqref{jointCLT-1}. Write
\begin{align*}
& A_n=\sum_{m=1}^n\mathbb{E}[(D^{(1)}_{n,m})^2|\mathcal{F}_{n,m-1}], \qquad B_n=\sum_{m=1}^n\mathbb{E}[(D^{(2)}_{n,m})^2|\mathcal{F}_{n,m-1}],\cr
& M_n=\sum_{m=1}^n\mathbb{E}[(D^{(1)}_{n,m})D^{(2)}_{n,m}|\mathcal{F}_{n,m-1}]. 
\end{align*}
Then, 
\[
\sum_{m=1}^n\mathbb{E}[D_{n,m}^2|\mathcal{F}_{n,m-1}] = u^2 A_n+v^2B_n+2uvM_n,
\]
In Sections~\ref{subsec:nullproof-DC}-\ref{subsec:nullproof-SQ}, we have shown that $A_n\xrightarrow{\mathbb{P}}1$ and $B_n\xrightarrow{\mathbb{P}}1$. We claim that 
\begin{equation} \label{nullproof-M}
M_n\xrightarrow{\mathbb{P}}0. 
\end{equation}
Then, it follows that $\sum_{m=1}^n\mathbb{E}[D_{n,m}^2|\mathcal{F}_{n,m-1}] \xrightarrow{\mathbb{P}} u^2\cdot 1+v^2\cdot 1+2uv\cdot 0 = 1$. This gives \eqref{jointCLT-1}.

It remains to show \eqref{nullproof-M}. Using the expressions of $D^{(1)}_{n,m}$ and $D^{(2)}_{n,m}$, we have
\[
M_n = \frac{\tilde{M}_n}{n^2\alpha_n^3(1-\alpha_n)\sqrt{n(n-1)(n-2)}},  
\]
where $\tilde{M}_n=\sum_{m=1}^n\mathbb{E}[(T_{n,m}-T_{n,m-1})(\tilde{Q}_{n,m}-\tilde{Q}_{n,m-1})|\mathcal{F}_{n,m-1}]$. We plug in the definitions of $T_{n,m}$ and $\tilde{Q}_{n,m}$ to get 
\[
\tilde{M}_n= \sum_{m=1}^n\left(\sum_{\substack{(j_1,j_2,j_3)\in\\ I_m\setminus I_{m-1}}}\sum_{\substack{(i_1,i_2,i_3,i_4)\in\\CC(I_m)\setminus CC(I_{m-1})}} \mathbb{E}\left[W_{j_1j_3}W_{j_2j_3}\cdot W_{i_1i_2}W_{i_2i_3}W_{i_3i_4}W_{i_4i_1}\middle|\mathcal{F}_{n,m-1}\right]\right). 
\]
Let's see when $\mathbb{E}[W_{j_1j_3}W_{j_2j_3} W_{i_1i_2}W_{i_2i_3}W_{i_3i_4}W_{i_4i_1}|\mathcal{F}_{n,m-1}]\neq 0$. Since $(i_1,i_2,i_3,i_4)\in CC(I_m)\setminus CC(I_{m-1})$, exactly one of the four indices must be $m$. We assume $i_1=m$ without loss of generality. Since $(j_1,j_2,j_3)\in I_m\setminus I_{m-1}$, exactly one of the three indices must be $m$. Without loss of generality, we assume either $j_1=m$ or $j_3=m$. 
If $j_1=m$ (and recall that we have assumed $i_1=m$), then  
\begin{align*}
& \mathbb{E}[W_{j_1j_3}W_{j_2j_3} W_{i_1i_2}W_{i_2i_3}W_{i_3i_4}W_{i_4i_1}|\mathcal{F}_{n,m-1}]\cr
=\quad &W_{j_2j_3}W_{i_2i_3}W_{i_3i_4} \cdot \mathbb{E}[W_{mj_3} W_{mi_2}W_{i_4m}|\mathcal{F}_{n,m-1}].
\end{align*}
It is nonzero only if $j_3=i_2=i_4$. However, this is impossible, because $i_2$ and $i_4$ need to be distinct. 
%
If $j_3=m$ (and recall that we have assumed $i_1=m$), we have 
\begin{align*}
& \mathbb{E}[W_{j_1j_3}W_{j_2j_3} W_{i_1i_2}W_{i_2i_3}W_{i_3i_4}W_{i_4i_1}|\mathcal{F}_{n,m-1}]\cr
=\quad &W_{i_2i_3}W_{i_3i_4}\cdot \mathbb{E}[W_{j_1m}W_{j_2m}W_{mi_2}W_{i_4m} |\mathcal{F}_{n,m-1}].
\end{align*}
Note that $j_1\neq j_2$ and $i_2\neq i_4$. For the above to be nonzero, we must have $\{i_2,i_4\}=\{j_1,j_2\}$. It follows that 
\begin{align} \label{nullproof-M2}
\tilde{M}_n &= 8\sum_{m=1}^{n}\sum_{\substack{1\leq i_2,i_3,i_4\leq m-1\\\text{(distinct)}}}W_{i_2i_3}W_{i_3i_4}\cdot \mathbb{E}[W^2_{mi_2}W^2_{mi_4} |\mathcal{F}_{n,m-1}]\cr
&=8\alpha_n^2(1-\alpha_n)^2\sum_{m=1}^n \sum_{(i_2,i_3,i_4)\in I_{m-1}}W_{i_2i_3}W_{i_3i_4}\cr
&=8 \alpha_n^2(1-\alpha_n)^2\sum_{(i_2,i_3,i_4)\in I_{n-1}}(n-\max\{i_2,i_3,i_4\})W_{i_2i_3}W_{i_3i_4}. 
\end{align}
As a result,  
\begin{align*}
\mathbb{E}[M_n^2] &=  \frac{\mathbb{E}[\tilde{M}^{2}_n]}{n^5(n-1)(n-2)\alpha_n^6(1-\alpha_n)^2}\cr
&= \frac{64 \alpha_n^4(1-\alpha_n)^4}{n^5(n-1)(n-2)\alpha_n^6(1-\alpha_n)^2}\times \mathbb{E}\left[\Biggl( \sum_{(i_2,i_3,i_4)\in I_{n-1}}(n-\max\{i_2,i_3,i_4\})W_{i_2i_3}W_{i_3i_4}\Biggr)^2\right]\cr
&\leq \frac{C}{n^7\alpha_n^2} \sum_{i_2,i_3,i_4}n^2\cdot \mathbb{E}[W^2_{i_2i_3}W^2_{i_3i_4}]\cr
&\leq \frac{C}{n^7\alpha_n^2}\times n^5\alpha_n^2 \quad \xrightarrow{n\to\infty}\quad 0. 
\end{align*}
Then, \eqref{nullproof-M} follows directly. This completes the proof of Theorem~\ref{thm:jointNull}. \qed

\section{Proof of Theorem~\ref{thm:chi2-alt}}
\label{appendix:chi2-alt}

Define
\[
U_n = \frac{\hat{\alpha}_n(1-\hat{\alpha}_n)}{\alpha_0(1-\alpha_0)}-1, \qquad\mbox{and}\qquad Z_n^*= \frac{\sum_{i=1}^n(d_i-\bar{d})^2}{(n-1)\alpha_0(1-\alpha_0)}-n. 
\]
By definition, $X_n=(1+U_n)^{-1}(n+Z_n^*)$. It follows that
\beq \label{proof-chi2alt-0}
\psi_n^{DC}\equiv \frac{X_n -n}{\sqrt{n}} =\frac{1}{\sqrt{n}(1+U_n)}(Z_n^* - nU_n).  
\eeq
The asymptotic behavior of $\psi_n^{DC}$ is mainly determined by $Z_n^*$. Below, we first calculate the mean and variance of $Z_n^*$; then, we use these results to study the mean and variance of $\psi_n^{DC}$.

\smallskip
\noindent
{\bf The mean and variance of $Z_n^*$.}  
We introduce a matrix
\[
\tilde{\Omega} = \Omega  - \alpha_0\bm{1}_n\bm{1}_n', \qquad\mbox{where}\quad \alpha_0=h'Ph. 
\]
Then, $A_{ij}=W_{ij}+\tilde{\Omega}_{ij}+\alpha_0$, for all $i\neq j$.  
Write $\tilde{\Omega}^*=\tilde{\Omega}-\diag(\tilde{\Omega})$. It follows that
\begin{align} \label{chisquare-alt-1}
\sum_{i=1}^n\left( d_i-\bar{d}\right)^2 =\;\; & \sum_{i=1}^n \left( \sum_{j: j\neq i}(W_{ij}+\tilde{\Omega}_{ij}+\alpha_0) - \frac{1}{n}\sum_{(k,\ell): k\neq \ell} (W_{k\ell}+\tilde{\Omega}_{k\ell}+\alpha_0 )\right)^2\cr
=\;\; & \sum_{i=1}^n \left( e_i'W\bm{1}_n +e_i'\tilde{\Omega}^*\bm{1}_n  - \frac{1}{n}\bm{1}_n'W\bm{1}_n - \frac{1}{n}\bm{1}_n'\tilde{\Omega}^*\bm{1}_n  \right)^2\cr
=\;\; & \sum_{i=1}^n \left( e_i'\tilde{\Omega}^*\bm{1}_n-\frac{1}{n}\bm{1}_n'\tilde{\Omega}^*\bm{1}_n\right)^2 + 2\sum_{i=1}^n \left( e_i'\tilde{\Omega}^*\bm{1}_n-\frac{1}{n}\bm{1}_n'\tilde{\Omega}^*\bm{1}_n\right)(e_i'W\bm{1}_n) \cr
&- 2 \sum_{i=1}^n\left( e_i'\tilde{\Omega}^*\bm{1}_n-\frac{1}{n}\bm{1}_n'\tilde{\Omega}^*\bm{1}_n\right)\left(\frac{1}{n}\bm{1}_n'W\bm{1}_n\right) + \sum_{i=1}^n (e_i'W\bm{1}_n)^2 \cr
& +\frac{1}{n}(\bm{1}_n'W\bm{1}_n)^2 - 2\sum_{i=1}^n (e_i'W\bm{1}_n)\left(\frac{1}{n}\bm{1}_n'W\bm{1}_n\right) \cr
=\;\; & \sum_{i=1}^n \left( e_i'\tilde{\Omega}^*\bm{1}_n-\frac{1}{n}\bm{1}_n'\tilde{\Omega}^*\bm{1}_n\right)^2 + 2\sum_{i=1}^n \left( e_i'\tilde{\Omega}^*\bm{1}_n-\frac{1}{n}\bm{1}_n'\tilde{\Omega}^*\bm{1}_n\right)(e_i'W\bm{1}_n) \cr
&+ \sum_{i=1}^n (e_i'W\bm{1}_n)^2 - \frac{1}{n}(\bm{1}_n'W\bm{1}_n)^2. 
\end{align}
We further combine the last two terms of \eqref{chisquare-alt-1}:
\begin{align*}
& \sum_{i=1}^n (e_i'W\bm{1}_n)^2 - \frac{1}{n}(\bm{1}_n'W\bm{1}_n)^2\cr
=\;\; & \sum_{i=1}^n\left(\sum_{j: j\neq i}W_{ij}\right)^2 - \frac{1}{n}\left( \sum_{i\neq j}W_{ij}\right)^2\cr
=\;\; & \sum_{i\neq j}W_{ij}^2 + \sum_{i,j,k \text{ dist}} W_{ij}W_{ik} - \frac{2}{n}\sum_{i\neq j}W_{ij}^2 - \frac{1}{n}\sum_{i\neq j}\sum_{\substack{ k\neq l \\ \{k,l\}\neq \{i,j\}}}W_{ij}W_{kl}\cr
=\;\; & \frac{n-2}{n}\sum_{i\neq j}W_{ij}^2 +\sum_{i,j,k, \text{ dist}}W_{ij}W_{ik}- \frac{1}{n}\sum_{i\neq j}\sum_{\substack{ k\neq l \\ \{k,l\}\neq \{i,j\}}}W_{ij}W_{kl}. 
\end{align*}
We plug it into \eqref{chisquare-alt-1} to get 
\begin{align} \label{chisquare-alt-2}
(n-1)\alpha_0(1-\alpha_0)Z_n^*& \equiv \sum_{i=1}^n\left( d_i-\bar{d}\right)^2 - n(n-1)\alpha_0(1-\alpha_0) \cr
&= Y_1 + 2Y_2 +Y_3 +Y_4 - Y_5,
\end{align}
where
\begin{align*}
Y_1 &= \sum_{i=1}^n \left( e_i'\tilde{\Omega}^*\bm{1}_n-\frac{1}{n}\bm{1}_n'\tilde{\Omega}^*\bm{1}_n\right)^2,\cr
Y_2 &= \sum_{i=1}^n \sum_{j\neq i} \left( e_i'\tilde{\Omega}^*\bm{1}_n-\frac{1}{n}\bm{1}_n'\tilde{\Omega}^*\bm{1}_n\right)W_{ij},\cr
Y_3 &= \left(\frac{n-2}{n}\sum_{i\neq j}W_{ij}^2\right)-n(n-1)\alpha_0(1-\alpha_0),\cr
Y_4 &= \sum_{i,j,k, \text{ dist}}W_{ij}W_{ik},\cr
Y_5 &= \frac{1}{n}\sum_{i\neq j}\sum_{\substack{ k\neq l \\ \{k,l\}\neq \{i,j\}}}W_{ij}W_{kl}. 
\end{align*}

We now compute the mean of $Z^*_n$. It is easy to see that
\beq \label{chisquare-alt-5}
\E[Z_n^*] = \frac{Y_1+\E[Y_3]}{(n-1)\alpha_0(1-\alpha_0)}.
\eeq
For $Y_1$, note that $\tilde{\Omega}^*=\tilde{\Omega}-\diag(\tilde{\Omega})$. Since $\Pi\bm{1}_K=\bm{1}_n$, we can re-write 
\[
\tilde{\Omega}= \Omega - \alpha_0 \Pi \bm{1}_K \bm{1}_K' \Pi' = \Pi\left(P-\alpha_0\bm{1}_K\bm{1}_K'\right)\Pi'=\Pi M\Pi'. 
\]
As a result, $\tilde{\Omega}\bm{1}_n=n \Pi Mh$, and $\bm{1}_n'\tilde{\Omega}\bm{1}_n=0$. 
We plug them into the expression of $Y_1$ and note that $(a+b)^2\geq \frac{a^2}{2}-b^2$, for any $a, b\in\mathbb{R}$. It follows that
\begin{align*}
Y_1 &= \|\tilde{\Omega}^*\bm{1}_n\|^2 - \frac{1}{n}(\bm{1}_n'\tilde{\Omega}^*\bm{1}_n)^2\cr
&=\|\tilde{\Omega}\bm{1}_n-\diag(\tilde{\Omega})\bm{1}_n\|^2 -\frac{1}{n} (\bm{1}_n'\diag(\tilde{\Omega})\bm{1}_n)^2 \cr
&\geq \frac{1}{2}\|\tilde{\Omega}\bm{1}_n\|^2 - \|\diag(\tilde{\Omega})\bm{1}_n\|^2 -\frac{1}{n} (\bm{1}_n'\diag(\tilde{\Omega})\bm{1}_n)^2\cr
&=\frac{n^2}{2}\|\Pi M h\|^2 - \sum_{i=1}^n \tilde{\Omega}^2_{ii}-\frac{1}{n}\left(\sum_{i=1}^n\tilde{\Omega}_{ii}\right)^2. 
\end{align*} 
Note that $\max_{i}|\tilde{\Omega}_{ii}|\leq \max_{k,l}|M_{kl}|= C\|M\|$. Moreover, since $G=n^{-1}\Pi'\Pi$ and $\lambda_{\min}(G)\geq c$, we have $\|\Pi M h\|^2=n(h'MGMh)\geq Cn\|Mh\|^2$, and $\|\Pi Mh\|^2\leq \|\Pi\|^2\|Mh\|^2\leq Cn\|Mh\|^2$. It follows that
\beq \label{chisquare-alt-Y1}
Y_1= \frac{n^2}{2}\|\Pi M h\|^2 - O(n\|M\|^2) \asymp n^3\|Mh\|^2. 
\eeq
For $Y_3$, we have
\[
\E[Y_3] = \frac{n-2}{n} \sum_{i\neq j}\Omega_{ij}(1-\Omega_{ij})-n(n-1)\alpha_0(1-\alpha_0). 
\]
Write $\Omega_{ij}(1-\Omega_{ij})=\alpha_0(1-\alpha_0)+(1-2\alpha_0)(\Omega_{ij}-\alpha_0)-(\Omega_{ij}-\alpha_0)^2$. Recalling that $\Omega_{ij}-\alpha_0=\tilde{\Omega}_{ij}$, we plug these results into $\E[Y_3]$ to get
\begin{align*} 
  \E[Y_3] &= \frac{n-2}{n}\sum_{i\neq j}\left[\alpha_0(1-\alpha_0) + (1-2\alpha_0)\tilde{\Omega}_{ij} - \tilde{\Omega}_{ij}^2\right] - n(n-1)\alpha_0(1-\alpha_0) \cr
  &=-2(n-1)\alpha_0(1-\alpha_0) +  \frac{n-2}{n}\left[ (1-2\alpha_0)\left( \bm{1}_n'\tilde{\Omega}\bm{1}_n -\sum_{i}\tilde{\Omega}_{ii}\right) - \sum_{i\neq j}\tilde{\Omega}_{ij}^2  \right]\cr
  &=-2(n-1)\alpha_0(1-\alpha_0) -  \frac{n-2}{n}\left[ (1-2\alpha_0)\sum_{i}\tilde{\Omega}_{ii}+ \sum_{i\neq j}\tilde{\Omega}_{ij}^2  \right].
\end{align*}
Then, $|\E[Y_3]| \leq Cn\alpha_0+Cn\|M\|+Cn^2\|M\|^2$. Recall that by assumption, $\|M\|\leq C\|Mh\|$, $n\alpha_0~\to~\infty$ and  $\delta_n=n^{-3/2}\alpha_0^{-1}\|Mh\|^2\to\infty$. It follows that 
\begin{align*}
    \frac{|\E[Y_3]|}{n^3\|Mh\|^2} &\leq \frac{C}{\sqrt{n}\delta_n}+\frac{C}{n^{3/4}\sqrt{n\alpha_0\delta_n}}+\frac{C}{n} \to 0.
\end{align*}
It yields that 
\beq  \label{chisquare-alt-Y3}
\E[Y_3] =   o(n^3\|Mh\|^2). 
\eeq
We plug \eqref{chisquare-alt-Y1}-\eqref{chisquare-alt-Y3} into \eqref{chisquare-alt-5} to get 
\beq   \label{chisquare-alt-3}
\E[Z_n^*] = \frac{(n/2)\|\Pi Mh\|^2 - o(n^3\|Mh\|^2)}{(n-1)\alpha_0(1-\alpha_0)} \asymp n^2\alpha_0^{-1}\|Mh\|^2.  
\eeq

We then compute the variance of $Z_n^*$, it is easy to see that
\[
\V(Z^*_n)\leq \frac{C\V(Y_2)+C\V(Y_3)+C\V(Y_4)+C\V(Y_5)}{(n-1)^2\alpha_0^2(1-\alpha_0)^2}. 
\]
By direct calculations, we know that
\begin{align*}
\V(Y_3)&\leq C\sum_{i< j}\E[W^4_{ij}]\leq C\sum_{i\neq j}\Omega_{ij}\leq Cn^2\alpha_0,\cr
\V(Y_4) &\leq C\sum_{i, j, k \text{ dist}}\E[W_{ij}^2]\E[W_{ik}^2]\leq Cn^3\alpha_0^2,\cr
\V(Y_5) &\leq \frac{C}{n^2}\sum_{\substack{ i\neq j, k\neq l \\ \{i,j\}\neq \{k,l\}}} \E[W_{ij}^2]\E[W_{k,l}^2]\leq Cn^2\alpha_0^2.
\end{align*}
In the previous steps, we have seen that $\tilde{\Omega}^*=\tilde{\Omega}-\diag(\tilde{\Omega})$, $\bm{1}_n'\tilde{\Omega}\bm{1}_n=0$, $\|\tilde{\Omega}\bm{1}_n\|^2=n^3h'MGMh$, $\Omega_{ij}\leq C\alpha_0$, and $|\tilde{\Omega}_{ii}|\leq C\|M\|$. It follows that
\begin{align*}
\V(Y_2) &\leq C \sum_{i\neq j} \left( e_i'\tilde{\Omega}^*\bm{1}_n-\frac{1}{n}\bm{1}_n'\tilde{\Omega}^*\bm{1}_n\right)^2\times \Omega_{ij}(1-\Omega_{ij})\cr
&=  C \sum_{i\neq j} \left[ e_i'\tilde{\Omega}\bm{1}_n + \tilde{\Omega}_{ii} -\frac{1}{n}\left(\bm{1}_n'\diag(\tilde{\Omega})\bm{1}_n\right) \right]^2\times \Omega_{ij}(1-\Omega_{ij})\cr
&\leq C\left[ n\|\tilde{\Omega}\bm{1}_n\|^2 + n\sum_{i}\tilde{\Omega}_{ii}^2 + \left(\bm{1}_n'\diag(\tilde{\Omega})\bm{1}_n\right)^2 \right] \times C\alpha_0\cr
&\leq Cn^4\alpha_0 \|Mh\|^2 + Cn^2\alpha_0\|\diag(M)\|^2\cr
&\leq Cn^4\alpha_0\|Mh\|^2. 
\end{align*}
We combine the above results and note that for $n$ big enough, $n\alpha_0\geq c$. It gives  
\begin{align} \label{chisquare-alt-4}
\V(Z^*_n) & \leq \frac{C}{n^2\alpha_0^2}\left( n^4\alpha_0\|Mh\|^2 + n^3\alpha_0^2\right)\cr
&\leq Cn^2\alpha_0^{-1}\|Mh\|^2 + Cn. 
\end{align}
In conclusion, the mean and variance of $Z_n^*$ are characterized by \eqref{chisquare-alt-3} and \eqref{chisquare-alt-4}, respectively.

\smallskip
\noindent
{\bf The mean and variance of $\psi_n^{DC}$.}  
We now show the claims of this theorem. First, consider the mean of $\psi_n^{DC}$. Recalling \eqref{proof-chi2alt-0} and letting $\Delta_n=(1+U_n)^{-1}U_n$, we have
\begin{align} \label{proof-chi2alt-1}
\sqrt{n}\, \mathbb{E}[\psi^{DC}_n] &\geq \mathbb{E}[Z_n^*] -  \mathbb{E}\Bigl|\frac{U_n}{1+U_n}Z_n^*\Bigr| - n\, \mathbb{E}\Bigl|\frac{U_n}{1+U_n}\Bigr|\cr
&\geq \mathbb{E}[Z_n^*] -  \sqrt{\mathbb{E}[\Delta_n^2]}\sqrt{\mathbb{E}[(Z_n^*)^2]} - n\sqrt{\mathbb{E}[\Delta^2_n]}. 
\end{align}
The mean and variance of $Z_n^*$ have been analyzed above. We now study $\Delta_n$, which is a function of $\hat{\alpha}_n$ and $\alpha_0$. Note that
\[
\max_{i,j}\Omega_{ij}\leq\max_{k,\ell}P_{k\ell}\leq \bm{1}_K'P\bm{1}_K \leq C h'Ph=C\alpha_0,
\]
where $\bm{1}_K'P\bm{1}_K \leq C h'Ph$ is because $\min_kh_k\geq c$. Since $\hat{\alpha}_n=\frac{1}{n(n-1)}\bm{1}_n' A\bm{1}_n$ and $\alpha_0=h'Ph=n^{-2}\bm{1}_n'\Omega\bm{1}_n$, we have
\begin{align} \label{AltProof-hatalpha}
\left|\E[\hat{\alpha}_n]-\alpha_0\right|&=\frac{1}{n(n-1)}\left| \bm{1}_n'\Omega\bm{1}_n - \bm{1}_n'\diag(\Omega)\bm{1}_n -n(n-1)\alpha_0\right|\cr
&= \frac{1}{n(n-1)} \left| n^2\alpha_0 - \bm{1}_n'\diag(\Omega)\bm{1}_n-n(n-1)\alpha_0\right| \;\; \leq Cn^{-1}\alpha_0, \cr
\V(\hat{\alpha}_n) & = \frac{4}{n^2(n-1)^2}\sum_{i<j}\Omega_{ij}(1-\Omega_{ij}) \;\; \leq Cn^{-2}\alpha_0. 
\end{align}
Furthermore, we write $\hat{\alpha}_n-\mathbb{E}[\hat{\alpha}_n]=\frac{2}{n(n-1)}\sum_{i<j}W_{ij}$, where $\{W_{ij}\}_{i<j}$ is a collection of independent, bounded, zero-mean variables. We apply Bernstein's inequality and use \eqref{AltProof-hatalpha} to get
\beq \label{AltProof-add}
\mathbb{P}\bigl(\bigl|\hat{\alpha}_n -\mathbb{E}[\alpha_n]\bigr|> t\bigr)\leq \exp\biggl(-  \frac{t^2/2}{Cn^{-2}\alpha_0+Cn^{-2}t}  \biggr), \qquad\mbox{for all }t>0. 
\eeq
Consider the event $E=\{|\hat{\alpha}_n-\alpha_0|<\delta\cdot \alpha_0\}$, for a sufficiently small constant $\delta>0$ to be determined. Using the above inequality, $\mathbb{P}(E^c)\leq \exp(-C\delta\cdot n^2\alpha_0)$ for big enough $n$. On the event $E$, we can derive a bound for $|\Delta_n|$. Recalling that $U_n=\frac{\hat{\alpha}_n(1-\hat{\alpha}_n)}{\alpha_0(1-\alpha_0)}$, we have
\[
\Delta_n= \frac{U_n}{1+U_n}=\frac{(\hat{\alpha}_n-\alpha_0)(1-\hat{\alpha}_n-\alpha_0)}{\hat{\alpha}_n(1-\hat{\alpha}_n)}. 
\]
Since $\alpha_0\leq 1-c$ for a constant $c\in (0,1)$, 
when $\delta$ is chosen properly small, $|\Delta_n|\leq C\alpha_0^{-1}|\hat{\alpha}_n-\alpha_0|$ on the event $E$, where the constant $C>0$ here does not depend on $\delta$. On the event $E^c$, according to the footnote on Page 3, $|\Delta_n|\leq Cn^2$. It follows that
\begin{align} \label{proof-chi2alt-2}
\mathbb{E}[\Delta_n^2] & \leq Cn^4\cdot\mathbb{P}(E^c) + C\alpha_0^{-2}\mathbb{E}[(\hat{\alpha}_n-\alpha_0)^2]\cr
&\leq Cn^4\cdot\mathbb{P}(E^c) + C\alpha_0^{-2}\bigl[ (\E[\hat{\alpha}_n]-\alpha_0)^2+\V(\hat{\alpha}_n)\bigr]\cr
&\leq Cn^4\exp(-C\delta n^2\alpha_0) + C\alpha_0^{-2}(n^{-2}\alpha_0^2+n^{-2}\alpha_0)\cr
&\leq Cn^{-2}\alpha_0^{-1}. 
\end{align}
We plug \eqref{proof-chi2alt-2} into \eqref{proof-chi2alt-1} and then utilize \eqref{chisquare-alt-3}-\eqref{chisquare-alt-4}.  \textcolor{black}{Recalling that we have defined $\delta_n=n^{3/2}\alpha_0^{-1}\|Mh\|^2$}, it yields that
\begin{align*} 
\mathbb{E}[\psi_n^{DC}]& \geq \frac{C}{\sqrt{n}}\biggl(n^2\alpha_0^{-1}\|Mh\|^2-\sqrt{Cn^{-2}\alpha_0^{-1}}\times\cr
&\sqrt{(n^2\alpha_0^{-1}\|Mh\|^2)^2+\bigl(n^2\alpha_0^{-1}\|Mh\|^2+n\bigr)} - n\sqrt{Cn^{-2}\alpha_0^{-1}}\biggr)\cr
&\textcolor{black}{= \frac{C}{\sqrt{n}}\biggl(\sqrt{n}\delta_n-\sqrt{Cn^{-2}\alpha_0^{-1}}\sqrt{n\delta_n^2+\sqrt{n}\delta_n+n} - n\sqrt{Cn^{-2}\alpha_0^{-1}}\biggr)}\cr
&\textcolor{black}{\geq C\delta_n\left(1-\frac{C}{\sqrt{n^2\alpha_0}}-\frac{C}{\sqrt{n^{5/2}\alpha_0\delta_n}}-\frac{C}{\sqrt{n^2\alpha_0\delta_n^2}}\right)-\frac{C}{\sqrt{n\alpha_0}}}\cr
&\textcolor{black}{\geq C\delta_n\Bigl[1-O\bigl(n^{-5/4}\alpha_0^{-1/2}\delta_n^{-1/2}+n^{-1}\alpha_0^{-1/2}\delta_n^{-1}\bigr)\Bigr]-O\bigl(n^{-1/2}\alpha_0^{-1/2}\bigr).}
\end{align*}
Now, assume that $\delta_n\geq C$. Then, there exists a constant $c_1>0$ such that
\begin{align}\label{proof-chi2alt-3}
    \mathbb{E}[\psi_n^{DC}] \geq c_1\delta_n-O\bigl(n^{-1/2}\alpha_0^{-1/2}\bigr).
\end{align}
This gives the first claim.

Next, consider the variance of $\psi_n^{DC}$. Note that $(1+U_n)^{-1}=1-\Delta_n$ and $(1+U_n)^{-1}U_n=\Delta_n$. It follows from \eqref{proof-chi2alt-0} that $\sqrt{n}\, \psi_n^{DC} =Z^*_n - \Delta_nZ^*_n - n\Delta_n$. Therefore,
\begin{align} \label{proof-chi2alt-4}
\V(\psi_n^{DC})&\leq Cn^{-1}[\V(Z^*_n)+\V(\Delta_nZ^*_n)+n^2\V(\Delta_n)]\cr
&\leq Cn^{-1}\bigl(\V(Z^*_n) + \mathbb{E}[\Delta_n^2(Z^*_n)^2]+n^2\mathbb{E}[\Delta_n^2]\bigr) \cr
&\leq Cn^{-1}\Bigl(n^2\alpha_0^{-1}\|Mh\|^2+n + \mathbb{E}[\Delta_n^2(Z^*_n)^2]+\alpha_0^{-1}\Bigr),
\end{align}
where we have used \eqref{chisquare-alt-4} and \eqref{proof-chi2alt-2} in the last inequality. 

We calculate $\mathbb{E}[\Delta_n^2(Z^*_n)^2]$. For a large enough constant $B_0>0$, we define an event 
\[
E_1=\bigl\{ |\hat{\alpha}_n-\mathbb{E}[\hat{\alpha}_n]|\leq B_0n^{-1}\sqrt{\alpha_0\log(n)}\bigr\}.
\]
By \eqref{AltProof-add}, $\mathbb{P}(E_1^c)\leq \exp(-B\log(n))$, where the constant $B>0$ is a monotone increasing function of $B_0$. With a properly large $B_0$, we can make $\exp(-B\log(n))=o(n^8\alpha_0^{-2})$. Now, on the event $E_1$, we have $|\Delta_n|\leq C\alpha_0^{-1}|\hat{\alpha}_n-\alpha_0|\leq Cn^{-1}\alpha_0^{-1/2}\sqrt{\log(n)}$. On the event $E^c$, we note that $|\Delta_n|\leq Cn^{2}$ and $|Z_n^*|\leq Cn^2\alpha_0^{-1}$ hold uniformly. It follows that 
\begin{align}  \label{proof-chi2alt-5}
\mathbb{E}[\Delta_n^2(Z^*_n)^2]& = \mathbb{E}[\Delta_n^2(Z^*_n)^2\cdot I_{E_1^c}] +  \mathbb{E}[\Delta_n^2(Z^*_n)^2\cdot I_{E_1}]  \cr
& \leq Cn^8\alpha_0^{-2}\cdot \exp(-B\log(n)) + Cn^{-2}\alpha_0^{-1}\log(n)\mathbb{E}\bigl[(Z_n^*)^2\cdot I_{E^c}\bigr]\cr
&\leq o(1) + Cn^{-2}\alpha_0^{-1}\log(n) \bigl[ (\mathbb{E}[Z_n^*])^2+\V(Z_n^*)\bigr]\cr
&\leq o(1) +\frac{C\log(n)}{n^2\alpha_0} \Bigl[ (n^2\alpha_0^{-1}\|Mh\|^2)^2+n^2\alpha_0^{-1}\|Mh\|^2+n\Bigr], 
\end{align}
where in the last inequality we have used \eqref{chisquare-alt-3}-\eqref{chisquare-alt-4}. We plug \eqref{proof-chi2alt-5} into \eqref{proof-chi2alt-4} to get
\begin{align} \label{proof-chi2alt-6}
\V(\psi_n^{DC})&\textcolor{black}{\leq C\Bigl(1+n^{-1}\alpha_0^{-1}+n^{-1/2}\delta_n+\frac{\log(n)}{n^3\alpha_0}(n\delta_n^2+\sqrt{n}\delta_n+n)\Bigr)}\cr
&\leq C\Bigl[1+ n^{-1/2}\delta_n+n^{-2}\alpha_0^{-1}\delta_n^2\log(n)\Bigr]. 
\end{align}
This gives the second claim. \qed


\section{Proof of Theorem~\ref{thm:SQ-alt}}\label{appendix:proof_SQ_alt}
Write $\alpha_1=\E[\hat{\alpha}_n]$, $\overline{\Omega}=\Omega-\alpha_1\mathds{1}_n\mathds{1}_n'$ and $\Delta_n = \alpha_1-\hat{\alpha}_n$. It follows that 
\begin{align*}
    Q_n &= \sum_{i,j,k,l \text{ dist.}}(A_{ij}-\hat{\alpha}_n)(A_{jk}-\hat{\alpha}_n)(A_{kl}-\hat{\alpha}_n)(A_{li}-\hat{\alpha}_n)\\
    &= \sum_{i,j,k,l \text{ dist.}}(W_{ij}+\overline{\Omega}_{ij}+\Delta_n)(W_{jk}+\overline{\Omega}_{jk}+\Delta_n)(W_{kl}+\overline{\Omega}_{kl}+\Delta_n)(W_{li}+\overline{\Omega}_{li}+\Delta_n). 
\end{align*}
Expanding the sum gives $3^4=81$ terms. Combining equal-valued terms, we have the following decomposition:
\begin{align}
Q_n &= X_1+4X_2+4X_3+4X_4+2X_5+8X_6+4X_7+4X_8+2X_9+4X_{10}+8X_{11} \cr
&+ 4X_{12}+8X_{13}+4X_{14}+4X_{15}+X_{16}+4X_{17}+4X_{18}+2X_{19}+4X_{20}+X_{21},
\end{align}
where the expressions of $X_1$-$X_{21}$ are presented in Column 4 of Table~\ref{tab:mean&SD}. In this table, we also list other information of each term, such as the degree in $W$ ($N_W$), in $\overline{\Omega}$ ($N_{\overline{\Omega}}$) and in $\Delta_n$ ($N_\Delta$). We plan to study the mean and variance of each of $X_1$-$X_{21}$ and then combine them to show the claims.

\begin{table}[tb!]   
\caption{The $21$ different types of the $81$ post-expansion sums of $Q_n$. The order of the mean and variance of each term will be derived in the proofs.}  
\centering
\scalebox{0.85}{
\begin{tabular}{lcclcr}
Type  & $\#$ &  ($N_W, N_{\widetilde{\Omega}}, N_\Delta)$ &     Representative  &  Mean & Variance \\
\hline 
$X_{1} $ & 1 & $(4,0,0)$ & $\sum_{i,j,k,l \text{ dist.}}W_{ij}W_{jk}W_{kl}W_{li}$                           & $0$ & $O(n^4\alpha_0^4)$ \\
$X_{2} $ & 4 & $(3,1,0)$ & $\sum_{i,j,k,l\text{ dist.}}W_{ij}W_{jk}W_{kl}\oOmega_{li}$                      & $0$ & $O(n^4\alpha_0^3\|M\|^2)$\\
$X_{3} $ & 4 & $(3,0,1)$ & $\sum_{i,j,k,l\text{ dist.}}W_{ij}W_{jk}W_{kl}\Delta_n$                          & $0$ & $O(n^2\alpha_0^4)$\\
$X_{4} $ & 4 & $(2,2,0)$ & $\sum_{i,j,k,l\text{ dist.}}W_{ij}W_{jk}\oOmega_{kl}\oOmega_{li}$                & $0$ & $O(n^4\alpha_0^2\|M\|^4)$\\
$X_{5} $ & 2 & $(2,2,0)$ & $\sum_{i,j,k,l\text{ dist.}}W_{ij}\oOmega_{jk}W_{kl}\oOmega_{li}$                & $0$ & $O(n^4\alpha_0^2\|M\|^4)$\\
$X_{6} $ & 8 & $(2,1,1)$ & $\sum_{i,j,k,l\text{ dist.}}W_{ij}W_{jk}\oOmega_{kl}\Delta_n$                    & $0$ & $O(n^3\alpha_0^3\|M\|^2)$\\
$X_{7} $ & 4 & $(2,1,1)$ & $\sum_{i,j,k,l\text{ dist.}}W_{ij}\oOmega_{jk}W_{kl}\Delta_n$                    & $0$ & $O(n^2\alpha_0^3\|M\|^2)$\\
$X_{8} $ & 4 & $(2,0,2)$ & $\sum_{i,j,k,l\text{ dist.}}W_{ij}W_{jk}\Delta_n^2$                              & $O(n^{1/2}\alpha_0^2)$ & $O(n\alpha_0^4)$\\
$X_{9} $ & 2 & $(2,0,2)$ & $\sum_{i,j,k,l\text{ dist.}}W_{ij}W_{kl}\Delta_n^2$                              & $O(\alpha_0^2)$ & $O(\alpha_0^4)$\\
$X_{10}$ & 4 & $(1,3,0)$ & $\sum_{i,j,k,l\text{ dist.}}W_{ij}\oOmega_{jk}\oOmega_{kl}\oOmega_{li}$          & $0$ & $O(n^6\alpha_0\|M\|^6)$\\
$X_{11}$ & 8 & $(1,2,1)$ & $\sum_{i,j,k,l\text{ dist.}}W_{ij}\oOmega_{jk}\oOmega_{kl}\Delta_n$              & $O(n^2\alpha_0\|M\|^2)$ & $O(n^4\alpha_0^2\|M\|^4)$\\
$X_{12}$ & 4 & $(1,2,1)$ & $\sum_{i,j,k,l\text{ dist.}}W_{ij}\oOmega_{jk}\Delta_n\oOmega_{li}$              & $O(n^2\alpha_0\|M\|^2)$ & $O(n^4\alpha_0^2\|M\|^4)$\\
$X_{13}$ & 8 & $(1,1,2)$ & $\sum_{i,j,k,l\text{ dist.}}W_{ij}\oOmega_{jk}\Delta_n^2$                        & $O(n^2\alpha_0^{3/2}\|M\|)$ & $O(n^4\alpha_0^3\|M\|^2)$\\
$X_{14}$ & 4 & $(1,1,2)$ & $\sum_{i,j,k,l\text{ dist.}}W_{ij}\oOmega_{kl}\Delta_n^2$                        & $O(n^2\alpha_0^{3/2}\|M\|)$ & $O(n^4\alpha_0^3\|M\|^2)$ \\
$X_{15}$ & 4 & $(1,0,3)$ & $\sum_{i,j,k,l\text{ dist.}}W_{ij}\Delta_n^3$                                    & $O(\alpha_0^2)$ & $O(\alpha_0^4)$\\
$X_{16}$ & 1 & $(0,4,0)$ & $\sum_{i,j,k,l\text{ dist.}}\oOmega_{ij}\oOmega_{jk}\oOmega_{kl}\oOmega_{li}$    & $n^4\|M\|^4$ & 0\\
$X_{17}$ & 4 & $(0,3,1)$ & $\sum_{i,j,k,l\text{ dist.}}\oOmega_{ij}\oOmega_{jk}\oOmega_{kl}\Delta_n$        & $0$ & $O(n^6\alpha_0\|M\|^6)$\\
$X_{18}$ & 4 & $(0,2,2)$ & $\sum_{i,j,k,l\text{ dist.}}\oOmega_{ij}\oOmega_{jk}\Delta_n^2$                  & $O(n^2\alpha_0\|M\|^2)$ & $O(n^4\alpha_0^2\|M\|^4)$\\
$X_{19}$ & 2 & $(0,2,2)$ & $\sum_{i,j,k,l\text{ dist.}}\oOmega_{ij}\oOmega_{kl}\Delta_n^2$                  & $O(n^2\alpha_0\|M\|^2)$ & $O(n^4\alpha_0^2\|M\|^4)$\\
$X_{20}$ & 4 & $(0,1,3)$ & $\sum_{i,j,k,l\text{ dist.}}\oOmega_{ij}\Delta_n^3$                              & $0$ & $0$\\
$X_{21}$ & 1 & $(0,0,4)$ & $\sum_{i,j,k,l\text{ dist.}}\Delta_n^4$                                          & $O(\alpha_0^2)$ & $O(\alpha_0^4)$\\
\hline 
\end{tabular} 
} 
 \label{tab:mean&SD}
\end{table} 

In preparation, we derive some useful results. First, we study $|\overline{\Omega}_{ij}|$. 
Write $M=P-\alpha_0\bm{1}_K\bm{1}_K'$. Then,
\begin{align}\label{SQproof-alt-0}
    |\alpha_1-\alpha_0| &=|\E[\hat{\alpha}_n]-\alpha_0|= \left|\frac{1}{n(n-1)}\sum_{i\neq j}\pi_i'M\pi_j\right| \nonumber \\
    &= \left|\frac{1}{n(n-1)}\sum_{i,j}\pi_i'M\pi_j-\frac{1}{n(n-1)}\sum_{i}\pi_i'M\pi_i\right| \nonumber \\
    &\leq \frac{n}{n-1}\left|h'Mh\right|+\frac{\|M\|}{n-1}\leq \frac{C\|M\|}{n}, 
\end{align}
where we have used in the last line that $h'Mh=h'Ph-\alpha_0h'\bm{1}_K\bm{1}_K'h=0$.

Note that $\overline{\Omega}_{ij}=\pi_i'P\pi_j-\alpha_1=\pi_i'M\pi_j+\alpha_0-\alpha_1$. It follows that
\beq \label{SQproof-alt-1}
     |\overline{\Omega}_{ij}| \leq |\pi_i'M\pi_j|+|\alpha_0-\alpha_1| \leq C\|M\|. 
\eeq
Next, we study $\Delta_n$. By definition, 
\[ 
    \Delta_n = \E[\hat{\alpha}_n]-\hat{\alpha}_n = -\frac{1}{n(n-1)}\sum_{i\neq j}(A_{ij}-\Omega_{ij}) = -\frac{1}{n(n-1)}\sum_{i\neq j}W_{ij}. 
\]
Using properties of Bernoulli variables, we have $\E[W_{ij}^2]=\Omega_{ij}(1-\Omega_{ij})\leq \Omega_{ij}$ and $|\E[W_{ij}^m]|\leq C\Omega_{ij}$, for any fixed $m\geq 3$ (the constant $C$ may depend on $m$). Note that
\begin{equation*}
    \Omega_{ij} = \pi_i'P\pi_j \leq \bm{1}_K'P\bm{1}_K\leq C\alpha_0,
\end{equation*}
where we have used that $\min_{k}h_k>C/K$, which is a consequence of \eqref{cond1}. Additionally, 
\[
    \sum_{i,j}\Omega_{ij}=\sum_{i,j}\pi_i'P\pi_j=n^2 h'Ph=n^2\alpha_0.
\]
It follows that
\begin{align}   \label{SQproof-alt-2}
    \E[\Delta_n^2] &= \frac{4}{n^2(n-1)^2}\sum_{i<j}\E(W^2_{ij})\leq Cn^{-4}\sum_{i\neq j}\Omega_{ij}\leq Cn^{-2}\alpha_0, \cr
    |\E[\Delta_n^3]| &= \frac{8}{n^3(n-1)^3}\biggl|\E\biggl[\sum_{i<j,k<l,u<v}W_{ij}W_{kl}W_{uv}\biggr]\biggr|=\frac{8}{n^3(n-1)^3}\biggl|\E\biggl[\sum_{i<j}W_{ij}^3\biggr]\biggr|\cr
    &\leq Cn^{-6}\sum_{i<j}\Omega_{ij} \; \leq\; Cn^{-4}\alpha_0,\cr
    \E[\Delta_n^4] &=\frac{16}{n^4(n-1)^4}\Biggl(\sum_{i<j}\E[W_{ij}^4]+3\sum_{\substack{i<j,k<l\\(i,j)\neq(k,l)}}\E[W_{ij}^2]\E[W_{kl}^2]\Biggr)\cr
    &\leq Cn^{-8} \Biggl[\sum_{i<j}\Omega_{ij}+\Bigl(\sum_{i<j}\Omega_{ij}\Bigr)\Bigl(\sum_{k<\ell}\Omega_{k\ell}\Bigr)\Biggr]\leq Cn^{-4}\alpha_0^2,\cr
    \E[\Delta_n^8] &\leq Cn^{-16}\Biggl(\sum_{i<j,k<l, m<s,q<t}\E[W^2_{ij}]\E[W^2_{kl}]\E[W^2_{ms}]\E[W^2_{qt}]  \Biggr)\cr
    & \leq Cn^{-16}\Bigl(\sum_{i<j}\Omega_{ij}\Bigr)^4\leq Cn^{-8}\alpha_0^4. 
\end{align}
We shall frequently use \eqref{SQproof-alt-1} and \eqref{SQproof-alt-2} in the proof below.


\bigskip
\noindent
{\bf Mean and variance of $Q_n$}. 
We study the mean and variance of each of $X_1$-$X_{21}$, and combine them to get the mean and variance of $Q_n$. 

Consider $X_1=\sum_{i,j,k,l \text{ dist.}} W_{ij} W_{jk} W_{k l} W_{l i}$. It is easy to see that
\beq \label{SQproof-X1-mean}
\E[X_1]=0. 
\eeq
Furthermore, let $CC(I_n)$ be collection of equivalent classes of 4-tuples $(i,j,k,l)$ (see the proof of \eqref{eq:ideal} for details). By elementary probability,  
\begin{align*}
\V(X_1) &= \V\Biggl(8\sum_{CC(I_n)}W_{ij} W_{jk} W_{kl} W_{l i}\Biggr) \cr
&= 64 \sum_{CC(I_n)}\E[W^2_{ij}]\E[W^2_{jk}]\E[W^2_{kl}] \E[W^2_{l i}] \cr
&\leq C\sum_{i,j,k,l}\Omega_{ij}\Omega_{jk}\Omega_{kl}\Omega_{li}\;\; \leq C\text{Tr}(\Omega^4).
\end{align*}
Note that $\Omega=\Pi P\Pi'$ and $\Pi {\bf 1}_n={\bf 1}_K$. Also, we have defined $G=n^{-1}\Pi'\Pi$ in Section~\ref{subsec:main-alt}. It follows that
\[
\text{Tr}(\Omega^4) = n^4\text{Tr}(PGPGPGPG) = n^4\text{Tr}\left((G^{1/2}PG^{1/2})^4\right) \leq Kn^4\left\|G^{1/2}PG^{1/2}\right\|^4.
\]
From the definition of $G$, we have $G_{kl}=n^{-1}\sum_{i,j}\pi_i(k)\pi_j(l)\leq 1$ for all $1\leq k,l\leq K$. Hence $\|G\|\leq K^2$. In addition, recall that $\alpha_0=h'Ph$. By our assumption \eqref{cond1}, all the entries of $h$ are lower bounded by a constant $C>0$. It follows that $\alpha_0\geq C{\bf 1}_K'P{\bf 1}_K$. We immediately have 
\[
\text{Tr}(\Omega^4) \leq K^9n^4\|P\|^4 \leq K^9n^4(\bm{1}_K'P\bm{1}_K)^4\leq Cn^4\alpha_0^4,
\]
where we have used that $\|P\|\leq \bm{1}_K'P\bm{1}_K$ since $P$ is a nonnegative matrix. Combining the above gives 
\beq \label{SQproof-X1-var}
\V(X_1)\leq Cn^4\alpha_0^4. 
\eeq

Next, consider $X_2=\sum_{i,j,k,l\text{ dist.}}W_{ij}W_{jk}W_{kl}\oOmega_{li}$. It is easy to see that
\beq \label{SQproof-X2-mean}
\E[X_2]=0. 
\eeq
Furthermore,
\[
\V\left(X_2\right) = \V\left(2\sum_{\substack{i,j,k,l\text{ dist.}\\i<l}}W_{ij}W_{jk}W_{kl}\oOmega_{li}\right)\leq C\sum_{\substack{i,j,k,l\text{ dist.}\\i<l}} \E[W_{ij}^2]\E[W_{jk}^2]\E[W_{kl}^2]\oOmega_{li},
\]
where we have used that summands in the expression above are pairwise independent. It follows that 
\beq \label{SQproof-X2-var}
\V(X_2)\leq Cn^4\alpha_0^3\|M\|^2. 
\eeq

Next, consider $X_3=\sum_{i,j,k,l\text{ dist.}}W_{ij}W_{jk}W_{kl}\Delta_n$. Recall that
\begin{equation*}
    \Delta_n = \alpha_1-\hat{\alpha}_n = -\frac{2}{n(n-1)}\sum_{i<j}W_{ij}.
\end{equation*}
It follows that
\beq \label{SQproof-X3-mean}
\E[X_3]= -\frac{2}{n(n-1)}\E\left[\sum_{i,j,k,l \text{ dist.}}\sum_{s<t}W_{ij} W_{jk} W_{k l}W_{st}\right]=0. 
\eeq
Furthermore,
\begin{align*}
    \V(X_3) &= \frac{1}{n^2(n-1)^2}\V\left(\sum_{\substack{i,j,k,l \text{ dist.}\\ s\neq t}}W_{ij} W_{jk} W_{kl}W_{st}\right)\\
    &\leq \frac{C}{n^4}\E\left[\sum_{\substack{i,j,k,l \text{ dist.}\\ a,b,c,d \text{ dist.}\\ s\neq t, u\neq v}}W_{ij} W_{jk} W_{kl}W_{st}W_{ab} W_{bc} W_{cd}W_{uv}\right]\\
    &\leq \frac{C}{n^4}\left(\sum_{i,j,k,l,s,t \text{ dist.}}\E[W_{ij}^2 W_{jk}^2 W_{kl}^2W_{st}^2]+\sum_{i,j,k,l,t \text{ dist.}}\E[W_{ij}^2 W_{jk}^2 W_{kl}^2W_{lt}^2] \right.\\
    &\left. \sum_{i,j,k,l,t \text{ dist.}}\E[W_{ij}^2 W_{jk}^2 W_{kl}^2W_{kt}^2]+\sum_{i,j,k,l \text{ dist.}}\E[W_{ij}^2 W_{jk}^2 W_{kl}^2W_{lj}^2]+\right.\\
    &\left. \sum_{i,j,k,l \text{ dist.}}\E[W_{ij}^2 W_{jk}^2 W_{kl}^2W_{li}^2]+\sum_{i,j,k,l \text{ dist.}}\E[W_{ij}^4 W_{jk}^2 W_{kl}^2]+\right.\\
    &\left. +\sum_{i,j,k,l \text{ dist.}}\E[W_{ij}^2 W_{jk}^4 W_{kl}^2]\right).
\end{align*}
It follows that 
\beq \label{SQproof-X3-var}
\V(X_3)\leq \frac{C}{n^4}(n^6\alpha_0^4+n^5\alpha_0^4+n^4\alpha_0^4+n^4\alpha_0^3)\leq Cn^2\alpha_0^4.
\eeq

Next, consider $X_4=\sum_{i,j,k,l\text{ dist.}}W_{ij}W_{jk}\oOmega_{kl}\oOmega_{li}$. It is straightforward to see that 
\beq \label{SQproof-X4-mean}
\E[X_4] = 0. 
\eeq
Furthermore,
\begin{align*}\label{SQproof-X4-var}
    \V(X_4) &= \E\left[\sum_{\substack{i,j,k,l \text{ dist.}\\u,v,s,t \text{ dist.}}} W_{ij}W_{jk}W_{uv}W_{vs}\oOmega_{k l} \oOmega_{l i}\oOmega_{st} \oOmega_{tu}\right]\\
    &\leq C\|M\|^4\sum_{i,j,k,l \text{ dist.}}\E[W_{ij}^2]\E[W_{jk}^2],
\end{align*}
from which we obtain that
\beq\label{SQproof-X4-var}
\V(X_4) \leq Cn^4\alpha_0^2\|M\|^4.
\eeq

Next, consider $X_5=\sum_{i,j,k,l\text{ dist.}}W_{ij}\oOmega_{jk}W_{kl}\oOmega_{li}$. It is straightforward to see that 
\beq \label{SQproof-X5-mean}
\E[X_5] = 0. 
\eeq
Furthermore,
\begin{align*}\label{SQproof-X5-var}
    \V(X_5) &= \E\left[\sum_{\substack{i,j,k,l \text{ dist.}\\u,v,s,t \text{ dist.}}} W_{ij}W_{kl}W_{uv}W_{st}\oOmega_{jk} \oOmega_{li}\oOmega_{vs} \oOmega_{tu}\right]\\
    &\leq C\|M\|^4\sum_{i,j,k,l \text{ dist.}}\E[W_{ij}^2]\E[W_{kl}^2],
\end{align*}
from which we obtain that
\beq\label{SQproof-X5-var}
\V(X_5) \leq Cn^4\alpha_0^2\|M\|^4.
\eeq

Next, consider $X_6=\sum_{i,j,k,l\text{ dist.}}W_{ij}W_{jk}\oOmega_{kl}\Delta_n$. Using the definition of $\Delta_n$, we have
\begin{equation*}
    X_6 = -\frac{1}{n(n-1)}\sum_{i,j,k,l\text{ dist.}}\sum_{s\neq t}W_{ij}W_{jk}\oOmega_{kl}W_{st}.
\end{equation*}
It follows that
\beq \label{SQproof-X6-mean}
\E[X_6] = 0. 
\eeq
Furthermore,
\begin{align*}
    \V(X_6) &= \frac{1}{n^2(n-1)^2}\sum_{\substack{i,j,k,l\text{ dist.}\\a,b,c,d\text{ dist.}\\s\neq t, u\neq v}}\oOmega_{kl}\oOmega_{cd}\E[W_{ij}W_{jk}W_{st}W_{ab}W_{bc}W_{uv}]\\
    &\leq \frac{C\|M\|^2}{n^2}\sum_{\substack{i,j,k\text{ dist.}\\a,b,c\text{ dist.}\\s\neq t, u\neq v}}\E[W_{ij}W_{jk}W_{st}W_{ab}W_{bc}W_{uv}]\\
    &\leq \frac{C\|M\|^2}{n^2}\left(4\sum_{i,j,k,s,t\text{ dist.}}\E[W_{ij}^2W_{jk}^2W_{st}^2]+4\sum_{i,j,k,t\text{ dist.}}\E[W_{ij}^2W_{jk}^2W_{kt}^2]+\right.\\
    &\left. 2\sum_{i,j,k,t\text{ dist.}}\E[W_{ij}^2W_{jk}^2W_{jt}^2]+\sum_{i,j,k\text{ dist.}}\E[W_{ij}^2W_{jk}^2W_{ki}^2]+4\sum_{i,j,k\text{ dist.}}\E[W_{ij}^2W_{jk}^4] +\right.\\
    &\left. 4\sum_{i,j,k\text{ dist.}}\E[W_{ij}^3W_{jk}^3]\right).
\end{align*}
As a result, we obtain
\beq\label{SQproof-X6-var}
\V(X_6) \leq \frac{C\|M\|^2}{n^2}(n^5\alpha_0^3+n^4\alpha_0^3+n^3\alpha_0^3+n^3\alpha_0^2)\leq Cn^3\alpha_0^3\|M\|^2.
\eeq

Next, consider $X_7=\sum_{i,j,k,l\text{ dist.}}W_{ij}\oOmega_{jk}W_{kl}\Delta_n$. Similarly to $X_6$, it is easy to see that
\beq \label{SQproof-X7-mean}
\E[X_7] = 0. 
\eeq
Furthermore,
\begin{align*}
    \V(X_7) &= \frac{1}{n^2(n-1)^2}\sum_{\substack{i,j,k,l\text{ dist.}\\a,b,c,d\text{ dist.}\\s\neq t, u\neq v}}\oOmega_{jk}\oOmega_{bc}\E[W_{ij}W_{kl}W_{st}W_{ab}W_{cd}W_{uv}]\\
    &\leq \frac{C\|M\|^2}{n^4}\left(\sum_{i,j,k,l,s,t\text{ dist.}}\E[W_{ij}^2W_{kl}^2W_{st}^2]+\sum_{i,j,k,l,t\text{ dist.}}\E[W_{ij}^2W_{kl}^2W_{lt}^2]+\right.\\
    &\left. \sum_{i,j,k,l\text{ dist.}}\E[W_{ij}^2W_{jk}^2W_{kl}^2]+\sum_{i,j,k,l\text{ dist.}}\E[W_{ij}^2W_{kl}^4]+\sum_{i,j,k,l\text{ dist.}}\E[W_{ij}^3W_{kl}^3]\right).
\end{align*}
As a result, we obtain
\beq\label{SQproof-X7-var}
\V(X_7) \leq \frac{C\|M\|^2}{n^4}(n^6\alpha_0^3+n^5\alpha_0^3+n^4\alpha_0^3+n^4\alpha_0^2)\leq Cn^2\alpha_0^3\|M\|^2.
\eeq

Next, consider $X_8=\sum_{i,j,k,l\text{ dist.}}W_{ij}W_{jk}\Delta_n^2$. We have
\begin{align*}
    |\E[X_8]| &= (n-3)\left|\E\left[\Delta_n^2\sum_{i,j,k\text{ dist.}}W_{ij}W_{jk}\right]\right| \leq n\E[\Delta_n^4]^{1/2}\E\left[\left(\sum_{i,j,k\text{ dist.}}W_{ij}W_{jk}\right)^2\right]^{1/2}.
\end{align*}
It follows that
\beq \label{SQproof-X8-mean}
|\E[X_8]| \leq Cn^{-1}\alpha_0n^{3/2}\alpha_0\leq Cn^{1/2}\alpha_0^2. 
\eeq
Furthermore,
\begin{align*}
    \V(X_8) &\leq Cn^2\E\left[\Delta_n^4\left(\sum_{i,j,k\text{ dist.}}W_{ij}W_{jk}\right)^2\right]\\
    &\leq Cn^2\E[\Delta_n^8]^{1/2}\E\left[\left(\sum_{i,j,k\text{ dist.}}W_{ij}W_{jk}\right)^4\right]^{1/2}.
\end{align*}
The summands above can be grouped into 6 categories, where each category corresponds to a specific upper bound in terms of $n$ and $\alpha_0$. We obtain
\beq\label{SQproof-X8-var}
\V(X_8) \leq Cn^{-2}\alpha_0^2(n^6\alpha_0^4+n^5\alpha_0^4+n^4\alpha_0^4+n^4\alpha_0^3+n^3\alpha_0^3+n^3\alpha_0^2)^{1/2}\leq Cn\alpha_0^4.
\eeq

Next, consider $X_9=\sum_{i,j,k,l\text{ dist.}}W_{ij}W_{kl}\Delta_n^2$. We have
\begin{align*}
    |\E[X_9]| &= \left|\E\left[\Delta_n^2\sum_{i,j,k,l\text{ dist.}}W_{ij}W_{kl}\right]\right| \leq \E[\Delta_n^4]^{1/2}\E\left[\left(\sum_{i,j,k,l\text{ dist.}}W_{ij}W_{kl}\right)^2\right]^{1/2}.
\end{align*}
It follows that
\beq \label{SQproof-X9-mean}
|\E[X_9]| \leq Cn^{-2}\alpha_0n^2\alpha_0\leq C\alpha_0^2. 
\eeq
Furthermore,
\begin{align*}
    \V(X_9) &\leq C\E\left[\Delta_n^4\left(\sum_{i,j,k,l\text{ dist.}}W_{ij}W_{kl}\right)^2\right]\\
    &\leq C\E[\Delta_n^8]^{1/2}\E\left[\left(\sum_{i,j,k,l\text{ dist.}}W_{ij}W_{kl}\right)^4\right]^{1/2}.
\end{align*}
As for $X_8$, the summands above can be grouped into 6 categories, where each category corresponds to a specific upper bound in terms of $n$ and $\alpha_0$. We obtain
\beq\label{SQproof-X9-var}
\V(X_9) \leq Cn^{-4}\alpha_0^2(n^8\alpha_0^4+n^7\alpha_0^4+n^6\alpha_0^4+n^5\alpha_0^4+n^4\alpha_0^4+n^4\alpha_0^2)^{1/2}\leq C\alpha_0^4.
\eeq

Next, consider $X_{10}=\sum_{i,j,k,l\text{ dist.}}W_{ij}\oOmega_{jk}\oOmega_{kl}\oOmega_{li}$. It is straightforward to see that
\beq \label{SQproof-X10-mean}
|\E[X_{10}]|=0. 
\eeq
Furthermore,
\begin{align*}
    \V(X_{10}) &=\sum_{\substack{i,j,k,l\text{ dist.}\\a,b,c,d \text{ dist.}}}\overline{\Omega}_{jk}\overline{\Omega}_{kl}\overline{\Omega}_{li}\overline{\Omega}_{bc}\overline{\Omega}_{cd}\overline{\Omega}_{da}\E[W_{ij}W_{ab}]\\
    &\leq C\alpha_0\sum_{\substack{i,j,k,l\text{ dist.}\\c\neq d, c,d\notin\{i,j\}}}|\overline{\Omega}_{jk}\overline{\Omega}_{kl}\overline{\Omega}_{li}\overline{\Omega}_{jc}\overline{\Omega}_{cd}\overline{\Omega}_{di}|.
\end{align*}
As a result,
\beq\label{SQproof-X10-var}
\V(X_{10}) \leq C\alpha_0n^6\|M\|^6.
\eeq

Next, consider $X_{11}=\sum_{i,j,k,l\text{ dist.}}W_{ij}\oOmega_{jk}\oOmega_{kl}\Delta_n$. Using the definition of $\Delta_n$, we obtain
\begin{align*}
    |\E[X_{11}]| &= \left|\frac{1}{n(n-1)}\sum_{\substack{i,j,k,l \text{ dist.}\\u\neq v}} \overline{\Omega}_{jk}\overline{\Omega}_{kl}\E[W_{ij}W_{uv}]\right|\leq C\|M\|^2\sum_{i\neq j, u\neq v} |\E[W_{ij}W_{uv}]|.
\end{align*}
As a result,
\beq \label{SQproof-X11-mean}
|\E[X_{11}]| \leq Cn^2\alpha_0\|M\|^2. 
\eeq
Furthermore,
\begin{align*}
    \V(X_{11}) &\leq C\E\left[\left(\frac{1}{n(n-1)}\sum_{\substack{i,j,k,l\text{ dist.}\\u\neq v}}\overline{\Omega}_{jk}\overline{\Omega}_{kl}W_{ij}W_{u v}\right)^2\right]\\
    &\leq\frac{C}{n^4}\sum_{\substack{i,j,k,l\text{ dist.}\\a,b,c,d\text{ dist.}\\u\neq v,r\neq s}}|\overline{\Omega}_{jk}\overline{\Omega}_{kl}\overline{\Omega}_{bc}\overline{\Omega}_{cd}||\E[W_{ij}W_{u v}W_{ab}W_{r s}]|\\
    &\leq C\|M\|^4\sum_{\substack{i\neq j, a\neq b\\u\neq v,r\neq s}}|\E[W_{ij}W_{u v}W_{ab}W_{r s}]|\\
    &\leq C\|M\|^4\left(\sum_{i,j,a,b \text{ dist.}}\E[W_{ij}^2W_{ab}^2]+\sum_{i,j,b \text{ dist.}}\E[W_{ij}^2W_{jb}^2]+\sum_{i,j \text{ dist.}}\E[W_{ij}^4]\right).
\end{align*}
As a result,
\beq\label{SQproof-X11-var}
\V(X_{11}) \leq C\|M\|^4(n^4\alpha_0^2+n^3\alpha_0^2+n^2\alpha_0) \leq Cn^4\alpha_0^2\|M\|^4.
\eeq

Next, consider $X_{12}=\sum_{i,j,k,l\text{ dist.}}W_{ij}\oOmega_{jk}\Delta_n\oOmega_{li}$. Computations in this case are exactly equivalent to those for $X_{11}$, so we obtain:
\beq \label{SQproof-X12-mean}
|\E[X_{12}]| \leq Cn^2\alpha_0\|M\|^2. 
\eeq
and 
\beq\label{SQproof-X12-var}
\V(X_{12}) \leq C\|M\|^4(n^4\alpha_0^2+n^3\alpha_0^2+n^2\alpha_0) \leq Cn^4\alpha_0^2\|M\|^4.
\eeq

Next, consider $X_{13}=\sum_{i,j,k,l\text{ dist.}}W_{ij}\oOmega_{jk}\Delta_n^2$. We have for the mean:
\begin{align*}
    |\E[X_{13}]| &\leq \sum_{i,j,k,l\text{ dist.}}|\oOmega_{jk}|\E[W_{ij}\Delta_n^2] \leq \sum_{i,j,k,l\text{ dist.}}|\oOmega_{jk}|\E[W_{ij}^2]^{1/2}E[\Delta_n^4]^{1/2}\\
    &\leq Cn^4\|M\|\alpha_0^{1/2}E[\Delta_n^4]^{1/2}.
\end{align*}
It follows that
\beq \label{SQproof-X13-mean}
|\E[X_{13}]| \leq Cn^2\alpha_0^{3/2}\|M\|. 
\eeq
Furthermore,
\begin{align*}
    \V(X_{13}) &\leq \E\left[\sum_{\substack{i,j,k,l \text{ dist.}\\ a,b,c,d \text{ dist.}}} W_{ij}W_{ab}\oOmega_{jk}\oOmega_{bc}\Delta_n^4 \right]\leq Cn^4\|M\|^2\sum_{i\neq j, a\neq b}\E[W_{ij}W_{ab}\Delta_n^4]\\
    &\leq Cn^4\|M\|^2\sum_{i\neq j, a\neq b}\E[W_{ij}^2W_{ab}^2]^{1/2}\E[\Delta_n^8]^{1/2} \leq C\alpha_0^2\|M\|^2\sum_{i\neq j, a\neq b}\E[W_{ij}^2W_{ab}^2]^{1/2}\\
    &\leq C\alpha_0^2\|M\|^2\left(\sum_{i,j,a,b \text{ dist.}}\E[W_{ij}^2]^{1/2}\E[W_{ab}^2]^{1/2}+\sum_{i,j,b \text{ dist.}}\E[W_{ij}^2]^{1/2}\E[W_{jb}^2]^{1/2}+\right.\\
    &\left. \sum_{i,j \text{ dist.}}\E[W_{ij}^4]^{1/2}\right).
\end{align*}
As a result,
\beq\label{SQproof-X13-var}
\V(X_{13}) \leq C\alpha_0^2\|M\|^2(n^4\alpha_0+n^3\alpha_0+n^2\alpha_0^{1/2})\leq Cn^4\alpha_0^3\|M\|^2.
\eeq

Next, consider $X_{14}=\sum_{i,j,k,l\text{ dist.}}W_{ij}\oOmega_{kl}\Delta_n^2$. Computations in this case are exactly equivalent to those for $X_{13}$, so we obtain:
\beq \label{SQproof-X14-mean}
|\E[X_{14}]| \leq Cn^2\alpha_0^{3/2}\|M\|. 
\eeq
and
\beq\label{SQproof-X14-var}
\V(X_{14}) \leq C\alpha_0^2\|M\|^2(n^4\alpha_0+n^3\alpha_0+n^2\alpha_0^{1/2})\leq Cn^4\alpha_0^3\|M\|^2.
\eeq

Next, consider $X_{15}=\sum_{i,j,k,l\text{ dist.}}W_{ij}\Delta_n^3$. Using the definition of $\Delta_n$, note that
\begin{equation*}
    X_{15} = (n-2)(n-3)\Delta_n^3\sum_{i\neq j}W_{ij}= -n(n-1)(n-2)(n-3)\Delta_n^4.
\end{equation*}
It follows that
\beq \label{SQproof-X15-mean}
|\E[X_{15}]| \leq n^4\E[\Delta_n^4] \leq C\alpha_0^2. 
\eeq
and
\beq\label{SQproof-X15-var}
\V(X_{15}) \leq n^8\E[\Delta_n^8]\leq C\alpha_0^4.
\eeq

Next, consider $X_{16}=\sum_{i,j,k,l\text{ dist.}}\oOmega_{ij}\oOmega_{jk}\oOmega_{kl}\oOmega_{li}$. This is a non-stochastic term, whose variance is zero. We the focus on deriving a lower bound for $\E[X_{16}]=X_{16}$. Note that 
\begin{align} \label{SQproof-X2-temp}
    X_{16} &= \sum_{i,j,k,l} \overline{\Omega}_{ij}\overline{\Omega}_{jk}\overline{\Omega}_{kl}\overline{\Omega}_{li}-\sum_{i,j,k,l \text{ not dist.}} \overline{\Omega}_{ij}\overline{\Omega}_{jk}\overline{\Omega}_{kl}\overline{\Omega}_{li}\cr
    &= \text{Tr}(\overline{\Omega}^4)- \sum_{i,j,k,l \text{ not dist.}} \overline{\Omega}_{ij}\overline{\Omega}_{jk}\overline{\Omega}_{kl}\overline{\Omega}_{li}\cr
&=\text{Tr}(\overline{\Omega}^4) - O\bigl(n^3\|M\|^4\bigr), 
\end{align}
where the last equality comes from \eqref{SQproof-alt-1} and the observation that $(i,j,k,l)$ has at most 3 distinct values in this sum. In the derivation of \eqref{SQproof-alt-1}, we have seen that $\overline{\Omega}_{ij}=\pi_i'P\pi_j-\alpha_1=\pi_i'\overline{M}\pi_j$, where $\overline{M}=P-\alpha_1{\bf 1}_K{\bf 1}_K'=M+(\alpha_0-\alpha_1){\bf 1}_K{\bf 1}_K'$. This implies that 
\[
\tOmega=\Pi \overline{M}\Pi'. 
\]
Recall that $G=n^{-1}\Pi'\Pi$. We have
\begin{align*}
    \Tr(\tOmega^4) &= \Tr((\Pi \overline{M}\Pi')^4) = n^4\Tr((G^{1/2}\overline{M}G^{1/2})^4)\cr
    &=n^4\|(G^{1/2}\overline{M}G^{1/2})^2\|_F^2\cr
    & \asymp n^4\|(G^{1/2}\overline{M}G^{1/2})^2\|^2\cr
    &\asymp n^4 \|G^{1/2}\overline{M}G^{1/2}\|^4. 
\end{align*}
Note that $\|G^{1/2}\overline{M}G^{1/2}\|\leq \|\overline{M}\|\|G\|$. Additionally,  $\|\overline{M}\|\leq \|G^{-1}\|\|G^{1/2}\overline{M}G^{1/2}\|$. By the definition of $G$ and our assumption \eqref{cond1}, $\|G\|\leq C$ and $\|G^{-1}\|\leq C$. It follows that $\|G^{1/2}\overline{M}G^{1/2}\|\asymp\|\overline{M}\|$. We thus have
\[
\Tr(\tOmega^4)\asymp n^4\|\overline{M}\|^4 = n^4\|M+(\alpha_0-\alpha_1){\bf 1}_K{\bf 1}_K'\|^4. 
\]
Recall now from \eqref{SQproof-alt-0} that $|\alpha_0-\alpha_1|=O(n^{-1}\|M\|)$.  Hence, by Weyl's inequality
\begin{align*}
    \left|\|\overline{M}\|- \|M\|\right| \leq K|\alpha_0-\alpha_1|\leq \frac{CK\|M\|}{n},
\end{align*}
which implies that $\|\overline{M}\|\asymp\|M\|$, so $\Tr(\tOmega^4)\asymp n^4\|M\|^4$. Plugging it into \eqref{SQproof-X2-temp} gives
\beq \label{SQproof-X2-mean}
X_{16}=\E[X_{16}]\asymp n^4\|M\|^4. 
\eeq 

Next, consider $X_{17}=\sum_{i,j,k,l\text{ dist.}}\oOmega_{ij}\oOmega_{jk}\oOmega_{kl}\Delta_n$. It is straightforward to see that
\beq \label{SQproof-X17-mean}
\E[X_{17}] =0. 
\eeq
Furthermore,
\beq\label{SQproof-X17-var}
\V(X_{17}) \leq \left(\sum_{i,j,k,l \text{ dist.}} \overline{\Omega}_{ij}\overline{\Omega}_{jk}\overline{\Omega}_{kl}\right)^2\E[\Delta_n^2] \leq C\alpha_0n^6\|M\|^6.
\eeq

Next, consider $X_{18}=\sum_{i,j,k,l\text{ dist.}}\oOmega_{ij}\oOmega_{jk}\Delta_n^2$. We first note that $X_{18} = (n-3)\Delta_n^2\sum_{i,j,k \text{ dist.}}\overline{\Omega}_{ij}\overline{\Omega}_{jk}$. Hence,
\beq \label{SQproof-X18-mean}
|\E[X_{18}]| \leq \frac{C\alpha_0}{n}\left|\sum_{i,j,k \text{ dist.}}\overline{\Omega}_{ij}\overline{\Omega}_{jk}\right|\leq C\alpha_0n^2\|M\|^2. 
\eeq
Furthermore,
\beq\label{SQproof-X18-var}
\V(X_{18}) \leq n^2\left(\sum_{i,j,k \text{ dist.}}\tilde{\Omega}_{ij}\tilde{\Omega}_{jk}\right)^2\E[\Delta_n^4]\leq C\alpha_0^2n^4\|M\|^4.
\eeq

Next, consider $X_{19}=\sum_{i,j,k,l\text{ dist.}}\oOmega_{ij}\oOmega_{kl}\Delta_n^2$. We have
\beq \label{SQproof-X19-mean}
|\E[X_{19}]| \leq \frac{C\alpha_0}{n^2}\left|\sum_{i,j,k,l \text{ dist.}}\overline{\Omega}_{ij}\overline{\Omega}_{kl}\right|\leq C\alpha_0n^2\|M\|^2. 
\eeq
Furthermore,
\beq\label{SQproof-X19-var}
\V(X_{19}) \leq \left(\sum_{i,j,k,l \text{ dist.}}\tilde{\Omega}_{ij}\tilde{\Omega}_{kl}\right)^2\E[\Delta_n^4]\leq C\alpha_0^2n^4\|M\|^4.
\eeq

Next, consider $X_{20}=\sum_{i,j,k,l\text{ dist.}}\oOmega_{ij}\Delta_n^3$. Notice that
\begin{align*}
    X_{20} &= \Delta_n^3(n-2)(n-3)\sum_{i\neq j}\overline{\Omega}_{ij} = \Delta_n^3(n-2)(n-3)\left(\sum_{i\neq j}\Omega_{ij}-n(n-1)\alpha_1\right)=0.
\end{align*}
It follows that
\beq \label{SQproof-X20-mean}
\E[X_{20}]=0, 
\eeq
and
\beq\label{SQproof-X20-var}
\V(X_{20}) = 0.
\eeq

Next, consider $X_{21}=\sum_{i,j,k,l\text{ dist.}}\Delta_n^4$. Note that $X_{21}=n(n-1)(n-2)(n-3)\Delta_n^4$. As a result,
\beq \label{SQproof-X21-mean}
\E[X_{21}] \leq C\alpha_0^2, 
\eeq
and
\beq\label{SQproof-X21-var}
\V(X_{21}) \leq C\alpha_0^4.
\eeq

\smallskip
\noindent
{\bf Mean and variance of $Q_n/(2\sqrt{2}n^2\alpha_0^2)$.}
We use the results stored in Table~\ref{tab:mean&SD} in order to provide a lower bound for $\E[Q_n/(2\sqrt{2}n^2\alpha_0^2)]$ and an upper bound for $\V(Q_n/(2\sqrt{2}n^2\alpha_0^2))$. Recall that we defined
\begin{equation*}
    \tau_n = \left(\frac{n\|M\|^2}{\alpha_0}\right)^2.
\end{equation*}
We obtain that
\begin{align}\label{SQproof-alt-3}
    \E\left[\frac{Q_n}{2\sqrt{2}n^2\alpha_0^2}\right] &\asymp n^4\|M\|^4+O(n^{1/2}\alpha_0^2+n^2\alpha_0\|M\|^2+n^2\alpha_0^{3/2}\|M\|) \nonumber\\
    &\asymp \tau_n\left(1+O\left(\frac{1}{n^{3/2}\tau_n}+\frac{1}{n\tau_n^{1/2}}+\frac{1}{n^{1/2}\tau_n^{3/4}}\right)\right).
\end{align}
Similarly, we observe that
\begin{align}\label{SQproof-alt-4}
    \mathrm{Var}\left(\frac{Q_n}{2\sqrt{2}n^2\alpha_0^2}\right) &= O\left(\frac{n^4\alpha_0^4+n^4\alpha_0^3\|M\|^2+n^4\alpha_0^2\|M\|^4+n^6\alpha_0\|M\|^6}{n^4\alpha_0^4}\right) \nonumber\\
    &= O\left(1+\frac{\tau_n^{1/2}}{n}+\frac{\tau_n}{n^2}+\frac{\tau_n^{3/2}}{n}\right).
\end{align}
Assuming that $\tau_n\geq C$, then we can write
\begin{equation}\label{SQproof-alt-4bis}
    \E\left[\frac{Q_n}{2\sqrt{2}n^2\alpha_0^2}\right] \asymp \tau_n, \qquad \text{and} \qquad \mathrm{Var}\left(\frac{Q_n}{2\sqrt{2}n^2\alpha_0^2}\right)=O\left(1+\frac{\tau_n^{3/2}}{n}\right).
\end{equation}

\smallskip
\noindent
{\bf Mean and variance of $\psi_n^{SQ}$. }
Recall that
\begin{equation*}
    \psi_n^{SQ} = \frac{Q_n}{2\sqrt{2}n^2\hat{\alpha}_n^2}.
\end{equation*}
In the sequel, we let $Z_n^*=Q_n/(2\sqrt{2}n^2\alpha_0^2)$ for ease of notation. First, we compute a lower bound on the mean of $\psi_n^{SQ}$. Note that
\begin{align*}
    \E[\psi_n^{SQ}] &= \E\left[\left(\frac{\alpha_0}{\hat{\alpha}_n}\right)^2Z_n^*\right] \geq \E[Z_n^*]+2\E\left[\left(\frac{\alpha_0-\hat{\alpha}_n}{\hat{\alpha}_n}\right)Z_n^*\right]+\E\left[\left(\frac{\alpha_0-\hat{\alpha}_n}{\hat{\alpha}_n}\right)^2Z_n^*\right]\\
    &\geq \E[Z_n^*]-C\sqrt{\E\left[(Z_n^*)^2\right]}\left\{\sqrt{\E\left[\left(\frac{\alpha_0-\hat{\alpha}_n}{\hat{\alpha}_n}\right)^2\right]}+\sqrt{\E\left[\left(\frac{\alpha_0-\hat{\alpha}_n}{\hat{\alpha}_n}\right)^4\right]}\right\}.
\end{align*}
Under the event $E$ defined in Appendix~\ref{appendix:chi2-alt}, it holds that $|\hat{\alpha}_n-\alpha_0|<\delta\alpha_0$, so we can derive the following upper bound:
\begin{equation*}
    \left|\frac{\alpha_0-\hat{\alpha}_n}{\hat{\alpha}_n}\right| \leq \frac{|\alpha_0-\hat{\alpha}_n|}{(1-\delta)\alpha_0}.
\end{equation*}
Under $E^c$, it holds that
\begin{equation*}
    \left|\frac{\alpha_0-\hat{\alpha}_n}{\hat{\alpha}_n}\right| \leq Cn^2.
\end{equation*}
We thus have
\begin{align*}
    \E\left[\left(\frac{\alpha_0-\hat{\alpha}_n}{\hat{\alpha}_n}\right)^2\right] &\leq Cn^4\mathbb{P}(E^c)+C\alpha_0^{-2}\E[(\alpha_0-\hat{\alpha}_n)^2]\\
    &\leq Cn^4\mathbb{P}(E^c)+C\alpha_0^{-2}(\alpha_0-\E[\hat{\alpha}_n])^2+C\alpha_0^{-2}\V(\hat{\alpha}_n)\\
    &\leq Cn^4\mathbb{P}(E^c)+\frac{C}{(n-1)^2}+\frac{C}{n^2\alpha_0}\leq \frac{C}{n^2\alpha_0} = o(1).
\end{align*}
Similarly,
\begin{align*}
    \E\left[\left(\frac{\alpha_0-\hat{\alpha}_n}{\hat{\alpha}_n}\right)^4\right] &\leq Cn^8\mathbb{P}(E^c)+C\alpha_0^{-4}\E[(\alpha_0-\hat{\alpha}_n)^4]\\
    &\leq Cn^8\mathbb{P}(E^c)+C\alpha_0^{-4}\E[(\hat{\alpha}_n-\E[\hat{\alpha}_n])^4]+C\alpha_0^{-4}(\E[\hat{\alpha}_n]-\alpha_0)^4\\
    &\leq Cn^8\mathbb{P}(E^c)+\frac{C}{n^4}\\ &+\frac{C\alpha_0^{-4}}{n^8}\E\left[\sum_{\substack{i<j,k<l\\u<v,r<t}}(A_{ij}-\Omega_{ij})(A_{kl}-\Omega_{kl})(A_{uv}-\Omega_{uv})(A_{rs}-\Omega_{rs})\right]\\
    &\leq Cn^8\mathbb{P}(E^c)+\frac{C}{n^4}+ \frac{C\alpha_0^{-4}}{n^8}(n^4\alpha_0^2+n^2\alpha_0) \leq \frac{C}{\alpha_0^2 n^4} = o(1).
\end{align*}
It follows that, for $n$ big enough,
\begin{align*}
    \E[\psi_n^{SQ}] &\geq \E[Z_n^*]-o\left(\sqrt{\E\left[(Z_n^*)^2\right]}\right)= \E[Z_n^*]-o\left(\sqrt{\V(Z_n^*)+\E[Z_n^*]^2}\right)\\
    &= \E[Z_n^*]-o\left(\sqrt{1+n^{-1}\tau_n^{1/2}+n^{-2}\tau_n+n^{-1}\tau_n^{3/2}+\E[Z_n^*]^2}\right)\\
    &\geq \E[Z_n^*](1-o(1))-o\left(1+n^{-1/2}\tau_n^{1/4}+n^{-1}\tau_n^{1/2}+n^{-1/2}\tau_n^{3/4}\right).
\end{align*}
Assuming that $\tau_n\geq C$, we know from \eqref{SQproof-alt-4bis} that there exists a constant $c_2>0$ such that
\begin{equation}\label{SQproof-alt-5}
    \E[\psi_n^{SQ}] \geq c_2\tau_n-o\left(1+n^{-1/2}\tau_n^{3/4}\right).
\end{equation}

Next, we compute an upper bound on the variance of $\psi_n^{SQ}$. We have
\begin{align*}
    \V(\psi_n^{SQ}) &= \V\left(\left(\frac{\alpha_0}{\hat{\alpha}_n}\right)^2Z_n^*\right)= \V\left(\left(\frac{\alpha_0-\hat{\alpha}_n}{\hat{\alpha}_n}+1\right)^2Z_n^*\right)\\
    &\leq C\V(Z_n^*)+C\E\left[\left(\frac{\alpha_0-\hat{\alpha}_n}{\hat{\alpha}_n}\right)^2(Z_n^*)^2\right]+C\E\left[\left(\frac{\alpha_0-\hat{\alpha}_n}{\hat{\alpha}_n}\right)^4(Z_n^*)^2\right].
\end{align*}
Recall the event $E_1$ defined in Appendix~\ref{appendix:chi2-alt}. We had that $\mathbb{P}(E_1^c)\leq \exp(-B\log(n))$, where $B$ is a constant chosen large enough. Then, on the event $E_1$, we have that $|(\alpha_0-\hat{\alpha}_n)/\hat{\alpha}_n|\leq Cn^{-1}\alpha_0^{-1/2}\sqrt{\log(n)}$. On the event $E_1^c$, it holds uniformly that $|(\alpha_0-\hat{\alpha}_n)/\hat{\alpha}_n|\leq Cn^2$ and $|Z_n^*|\leq n^2\alpha_0^{-2}$. It follows that
\begin{align*}
    \E\left[\left(\frac{\alpha_0-\hat{\alpha}_n}{\hat{\alpha}_n}\right)^2(Z_n^*)^2\right] &\leq Cn^{8}\alpha_0^{-4}\mathbb{P}(E_1^c)+Cn^{-2}\alpha_0^{-1}\log(n)\E[(Z_n^*)^2],\\
    \E\left[\left(\frac{\alpha_0-\hat{\alpha}_n}{\hat{\alpha}_n}\right)^4(Z_n^*)^2\right] &\leq Cn^{12}\alpha_0^{-4}\mathbb{P}(E_1^c)+Cn^{-4}\alpha_0^{-2}\log(n)^2\E[(Z_n^*)^2].
\end{align*}
So we obtain that
\begin{align}\label{SQproof-alt-6}
    \V(\psi_n^{SQ}) &\leq  C\V(Z_n^*)+Cn^{-2}\alpha_0^{-1}\log(n)\E[(Z_n^*)^2]+o(1) \nonumber\\
    &\leq C\V(Z_n^*)+Cn^{-2}\alpha_0^{-1}\log(n)\E[(Z_n^*)]^2+o(1).
\end{align}
Recall from \eqref{SQproof-alt-4bis} that when $\tau_n\geq C$, $\E[Z_n^*]\asymp \tau_n$ and $\V(Z_n^*)=O(1+n^{-1}\tau_n^{3/2})$. It follows that
\begin{equation}\label{SQproof-alt-6bis}
    \V(\psi_n^{SQ}) = O\left(1+\frac{\tau_n^{3/2}}{n}+\frac{\log(n)\tau_n^2}{n^2\alpha_0}\right).
\end{equation}
\qed

\section{Proof of Corollary~\ref{cor:chi2}}\label{appendix:proof_cor_chi2}

Let $\psi^{DC}_n$ denote the degree test statistic as in the proof of Theorem~\ref{thm:chi2-alt}. Let $\epsilon\in(0,1)$ and $q_\epsilon$ be the $(1-\epsilon)$-quantile of the standard normal distribution. 

Under the alternative hypothesis, we suppose that $\delta_n\to\infty$. It follows from Theorem~\ref{thm:chi2-alt} that 
\begin{equation*}
    \E[\psi^{DC}_n] \geq c_1\delta_n, \quad \text{and} \quad \V(\psi^{DC}_n) =O(1+n^{-1/2}\delta_n+n^{-2}\alpha_0^{-1}\delta_n^2\log(n)).
\end{equation*}
We have, for $n$ big enough,
\begin{align*}
    \mathbb{P}\left(\psi^{DC}_n<q_\epsilon\right) &= \mathbb{P}\left(\E[\psi^{DC}_n]-\psi^{DC}_n>\E[\psi^{DC}_n]-q_\epsilon\right)\leq \frac{C\V(\psi^{DC}_n)}{\E[\psi^{DC}_n]^2} \asymp \frac{1}{\text{SNR}(\psi^{DC}_n)^2},
\end{align*}
where we have seen that $\text{SNR}(\psi^{DC}_n)\to\infty$ if $\delta_n\to\infty$ under the alternative (see the paragraph before the statement of Corollary~\ref{cor:chi2}). It follows that under the alternative, the power of the test
\begin{equation}\label{corchi2-proof-1}
    \mathbb{P}\left(\psi^{DC}_n>q_\epsilon\right) \xrightarrow[n\to\infty]{} 1.
\end{equation}

Furthermore, under the null hypothesis, we know from Corollary~\ref{cor:null} that $\psi^{DC}_n\xrightarrow[]{\mathcal{L}}\mathcal{N}(0,1)$, hence the level of the test tends to $\epsilon$ as $n\to\infty$. \qed

\section{Proof of Corollary~\ref{cor:SQ}}\label{appendix:proof_cor_SQ}

Let $\psi^{SQ}_n$ denote the degree test statistic as in the proof of Theorem~\ref{thm:SQ-alt}. Let $\epsilon\in(0,1)$ and $q_\epsilon$ be the $(1-\epsilon)$-quantile of the standard normal distribution. 

Under the alternative hypothesis, we suppose that $\tau_n\to\infty$. It follows from Theorem~\ref{thm:chi2-alt} that 
\begin{equation*}
    \E[\psi^{SQ}_n] \geq c_2\tau_n, \quad \text{and} \quad \V(\psi^{SQ}_n) =O(1+n^{-1}\tau_n^{3/2}+n^{-2}\alpha_0^{-1}\tau_n^2\log(n)).
\end{equation*}
We have, for $n$ big enough,
\begin{align*}
    \mathbb{P}\left(\psi^{SQ}_n<q_\epsilon\right) &= \mathbb{P}\left(\E[\psi^{SQ}_n]-\psi^{SQ}_n>\E[\psi^{SQ}_n]-q_\epsilon\right)\leq \frac{C\V(\psi^{SQ}_n)}{\E[\psi^{SQ}_n]^2} \asymp \frac{1}{\text{SNR}(\psi^{SQ}_n)^2},
\end{align*}
where we have seen that $\text{SNR}(\psi^{SQ}_n)\to\infty$ if $\tau_n\to\infty$ under the alternative (see the paragraph before the statement of Corollary~\ref{cor:SQ}). It follows that under the alternative, the power of the test
\begin{equation}\label{corSQ-proof-1}
    \mathbb{P}\left(\psi^{SQ}_n>q_\epsilon\right) \xrightarrow[n\to\infty]{} 1.
\end{equation}

Furthermore, under the null hypothesis, we know from Corollary~\ref{cor:null} that $\psi^{SQ}_n\xrightarrow[]{\mathcal{L}}\mathcal{N}(0,1)$, hence the level of the test tends to $\epsilon$ as $n\to\infty$. \qed

\section{Proof of Theorem~\ref{thm:PET-alt}}\label{appendix:proof_PET_alt}

As in the proofs of Theorem~\ref{thm:chi2-alt} and Theorem~\ref{thm:SQ-alt}, we let $\psi^{DC}_n$ denote the degree chi-squared test statistic and $\psi^{SQ}_n$ denote the oSQ statistic. Recall that the PET statistic is
\[
S_n = \left(\psi^{DC}_n\right)^2+\left(\psi^{SQ}_n\right)^2.
\]
Let $A>0, \epsilon>0$ be arbitrary constants. Then,
\begin{align*}
    \mathbb{P}(S_n<A) &\leq \min\left\{ \mathbb{P}\left(\psi^{DC}_n<\sqrt{A}\right),  \mathbb{P}\left(\psi^{SQ}_n<\sqrt{A}\right)\right\}\\
    &\leq \min\left\{ \mathbb{P}\left(\E[\psi^{DC}_n]-\psi^{DC}_n>\E[\psi^{DC}_n]-\sqrt{A}\right),  \mathbb{P}\left(\E[\psi^{SQ}_n]-\psi^{SQ}_n>\E[\psi^{SQ}_n]-\sqrt{A}\right)\right\}.
\end{align*}
In the regime where $\max\{\delta_n,\tau_n\}\to\infty$, for any constant $B>0$, there exists $N>0$ such that for all $n>N$, $\delta_n>B$ or $\tau_n>B$. We will denote by $N(B)$ the smallest such constant. We choose $B\gg \sqrt{A}$ and $N>N(B)$ such that for all $n>N$,
\begin{align*}
    \frac{1}{B^2}+\frac{1}{n^{1/2}B}+\frac{\log(n)}{n^2\alpha_0} &< \frac{\epsilon}{C}\\
    \frac{1}{B^2}+\frac{1}{nB^{1/2}}+\frac{\log(n)}{n^2\alpha_0}&<\frac{\epsilon}{C}
\end{align*}
Now, suppose that we are in the case $\delta_n>B$. Then from Theorem~\ref{thm:chi2-alt}, we know that
\begin{align*}
    \E[\psi_n^{DC}] > c\delta_n>cB \qquad \text{ and } \qquad \V(\psi_n^{DC}) <C\left(1+\frac{\delta_n}{n^{1/2}}+\frac{\log(n)}{n^2\alpha_0}\delta_n^2\right).
\end{align*}
Then,
\begin{equation*}
    \mathbb{P}\left(\E[\psi^{DC}_n]-\psi^{DC}_n>\E[\psi^{DC}_n]-\sqrt{A}\right) \leq \frac{\V(\psi^{DC}_n)}{\E[\psi^{DC}_n]^2} \leq C\left(\frac{1}{\delta_n^2}+\frac{1}{n^{1/2}\delta_n}+\frac{\log(n)}{n^2\alpha_0}\right)<\epsilon,
\end{equation*}
which implies that $\mathbb{P}(S_n<A)<\epsilon$.

Now, suppose that we are in the case $\tau_n>B$. By Theorem~\ref{thm:SQ-alt}, we have
\begin{equation*}
    \mathbb{E}[\psi^{SQ}_n]  \geq  c \tau_n >cB \quad \text{ and } \quad \mathrm{Var}(\psi^{SQ}_n) <C\left(1+\frac{\tau_n^{3/2}}{n}+\frac{\log(n)\tau_n^2}{n^2\alpha_0}\right).
\end{equation*}
Then

\begin{equation*}
    \mathbb{P}\left(\E[\psi^{SQ}_n]-\psi^{SQ}_n>\E[\psi^{SQ}_n]-\sqrt{A}\right) \leq \frac{\V(\psi^{SQ}_n)}{\E[\psi^{SQ}_n]^2} \leq C\left(\frac{1}{\tau_n^2}+\frac{1}{n\tau_n^{1/2}+\frac{\log(n)}{n^2\alpha_0}}\right)<\epsilon,
\end{equation*}
which implies that $\mathbb{P}(S_n<A)<\epsilon$.



It follows that for all $n>N$, it holds that $\mathbb{P}(S_n<A)<\epsilon$. We have just shown that
\begin{equation}\label{PETproof-alt-1}
    S_n \xrightarrow[n\to\infty]{\mathbb{P}} \infty.
\end{equation}

Now, fix $\epsilon\in(0,1)$ and let $q$ denote the $(1-\epsilon)$-quantile of the $\chi^2_2(0)$ distribution. From Corollary~\ref{cor:null}, we know that as $n\to\infty$, the level of the test tends to $\epsilon$. From \eqref{PETproof-alt-1}, we know that under the alternative
\begin{equation*}
    \mathbb{P}(S_n>q) \xrightarrow[n\to\infty]{} 1,
\end{equation*}
so the power of the test tends to $1$ as $n\to\infty$.\qed

\section{Proof of Theorem~\ref{thm:LB}}\label{appendix:LB} 
Denote by $D_{\chi^2}(P_0\|P_1)$ the chi-square divergence between two hypotheses, where $P_0$ and $P_1$ denote the probability measures under two model, respectively. and then study the symmetric alternative and the asymmetric alternative separately. By definition,
\[
    1+D_{\chi^2}(P_0\|P_1)=\int\left(\frac{dP_1}{dP_0}\right)^2dP_0. 
\]
Letting $q_{ij}(\Pi)=\pi_i'P\pi_j$, we can write
\[
    dP_0=\prod_{i<j}\alpha^{A_{ij}}(1-\alpha)^{1-A_{ij}}, \qquad
    dP_1=\E_\Pi\left[\prod_{i<j}q_{ij}(\Pi)^{A_{ij}}(1-q_{ij}(\Pi))^{1-A_{ij}}\right]. 
\]
Let $\tilde{\Pi}$ be an independent copy of $\Pi$. Then it follows that
\begin{align*}
    \left(\frac{dP_1}{dP_0}\right)^2 &=\mathbb{E}_{\Pi,\tilde{\Pi}}\left[\prod_{i<j}\left(\frac{q_{ij}(\Pi)q_{ij}(\tilde{\Pi})}{\alpha^2}\right)^{A_{ij}}\left(\frac{(1-q_{ij}(\Pi))(1-q_{ij}(\tilde{\Pi}))}{(1-\alpha)^2}\right)^{1-A_{ij}}\right].
\end{align*}
We denote
\[
    \Sigma(A,\Pi,\tilde{\Pi}) = \prod_{i<j}\left(\frac{q_{ij}(\Pi)q_{ij}(\tilde{\Pi})}{\alpha^2}\right)^{A_{ij}}\left(\frac{(1-q_{ij}(\Pi))(1-q_{ij}(\tilde{\Pi}))}{(1-\alpha)^2}\right)^{1-A_{ij}},
\]
and further obtain using the Tonelli theorem that
\begin{align*}
     1+D_{\chi^2}(P_0\|P_1) & =\mathbb{E}_0[\mathbb{E}_{\Pi,\tilde{\Pi}}[\Sigma(A,\Pi,\tilde{\Pi})]]= \mathbb{E}_{\Pi,\tilde{\Pi}}[\mathbb{E}_0[\Sigma(A,\Pi,\tilde{\Pi})]]. 
\end{align*}
Recalling that for all $i<j$ the $A_{ij}$'s are mutually independent, we can calculate $\mathbb{E}_0[\Sigma(A,\Pi,\tilde{\Pi})]$ easily. The calculations yield that 
\[
     1+D_{\chi^2}(P_0\|P_1) = \mathbb{E}_{\Pi,\tilde{\Pi}}\biggl[ \prod_{i<j}\Bigl(1 + \frac{\Delta_{ij}\tilde{\Delta}_{ij}}{\alpha(1-\alpha)} \Bigr)\biggr],
\]
where for $i<j$ we have $\Delta_{ij} = \pi_i' P \pi_j - \alpha$ and $\tilde{\Delta}_{ij} = \tilde{\pi}_i' P \tilde{\pi_j} - \alpha$ . Since for all $x$ in $\mathbb{R}$ it holds that $1+x\leq e^x$, we can bound the above by
\begin{align} \label{LB-eq1}
 1+& D_{\chi^2}(P_0\|P_1)
\leq \mathbb{E}_{\Pi,\tilde{\Pi}}\biggl[\exp\Bigl(\sum_{i<j}\frac{\Delta_{ij}\tilde{\Delta}_{ij}}{\alpha(1-\alpha)}\Bigr)\biggr] \cr
&=\ \mathbb{E}_{\Pi,\tilde{\Pi}}\left[\exp\left(\frac{S}{2(1-\alpha)}\right)\right], \qquad\mbox{where}\quad S\equiv \alpha^{-1}\sum_{i\neq j}\Delta_{ij}\tilde{\Delta}_{ij}. 
\end{align}
Recall that we chose $\alpha=h'Ph$ for the null model. Let $y_i=\pi_i-h$ for $i=1,...,n$, hence $\mathbb{E}[y_i]=0$. We obtain, for all $i\neq j$
\begin{align*}
    \Delta_{ij}&=\pi_i'P\pi_j-\alpha = y_i'Py_j+h'Py_i+h'Py_j+h'Ph-\alpha\\
    &= y_i'Py_j+h'Py_i+h'Py_j. 
\end{align*}
Hence, $\mathbb{E}[\Delta_{ij}]=0$. We define the matrix $M=P-\alpha\bm{1}_K\bm{1}_K'$. For all $i\in\llbracket1,n\rrbracket$, $\pi_i'\bm{1}_K=h'\bm{1}_K=1$, which implies that $y_i'\bm{1}_K=0$. It follows that 
\beq \label{LB-add}
    \Delta_{ij} = y_i'My_j+h'My_i+h'My_j. 
\eeq
We plug \eqref{LB-add} into \eqref{LB-eq1} to decomposition $\Delta_{ij}\tilde{\Delta}_{ij}$ into 9 terms:
\begin{align*}
\Delta_{ij}\tilde{\Delta}_{ij}& =(y_i'My_j)(\tilde{y}_i'M\tilde{y}_j) +  \Bigl[  (h'My_i)(h'M\tilde{y}_i)+(h'My_j)(h'M\tilde{y}_j)  \Bigr] \cr
&+ \Bigl[ (y_i'My_j)(h'M\tilde{y}_i') + (y_i'My_j)(h'M\tilde{y}_j') + (h'My_i)(\tilde{y}_i'M\tilde{y}_j) + (h'My_j)(\tilde{y}_i'M\tilde{y}_j) \Bigr]\cr
&+ \Bigl[  (h'My_i)(h'M\tilde{y}_j)+(h'My_j)(h'M\tilde{y}_i)  \Bigr]. 
 \end{align*}
Summing over $(i,j)$ such that $i\neq j$ gives a total of 9 partial sums, which we denote by $S_1$, $S_{21}$, $S_{22}$, $S_{31}$, $S_{32}$, $S_{33}$, $S_{34}$, $S_{41}$ and $S_{42}$, respectively. For example, 
\begin{align} \label{LBproof-defineSums}
    S_1 & = \alpha^{-1}\sum_{i\neq j}(y_i'My_j)(\tilde{y}_i'M\tilde{y}_j),\cr
    S_{21} &= \alpha^{-1}(n-1)\sum_{i}(h'My_i)(h'M\tilde{y}_i),\cr
    S_{31} & = \alpha^{-1}\sum_{i\neq j}(y_i'My_j)(h'M\tilde{y}_i),\cr
    S_{41} & = \alpha^{-1}\sum_{i\neq j} (h'My_i)(h'M\tilde{y}_j). 
\end{align}
It follows that 
\[
S = S_1+\sum_{m=1}^2 S_{2m}+\sum_{m=1}^4 S_{3m}+ \sum_{m=1}^2 S_{4m}. 
\]
Combining \eqref{LB-eq1} and \eqref{LB-eq2} with Jensen's inequality, we have 
\begin{align*}   
1+ D_{\chi^2}(P_0\|P_1)
&\leq \mathbb{E}\left[\exp\left(\frac{S_1+\sum_{m=2}S_{2m}+\sum_{m=1}^4S_{3m}+\sum_{m=1}^2 S_{4m}}{2(1-\alpha)}\right)\right]\cr
&\leq \frac{1}{9} \exp\left(\frac{9|S_1|}{2(1-\alpha)} \right) + \frac{1}{9}\sum_{m=1}^2 \exp\left(\frac{9|S_{2m}|}{2(1-\alpha)} \right)+\cr
&\qquad + \frac{1}{9}\sum_{m=1}^4 \exp\left(\frac{9|S_{3m}|}{2(1-\alpha)} \right) + \frac{1}{9}\sum_{m=1}^2 \exp\left(\frac{9|S_{4m}|}{2(1-\alpha)} \right).  
 \end{align*}
Write $c_{\alpha}=9/[2(1-\alpha)]$. To show the claim, it suffices to show that
\beq  \label{LB-eq2}
\E\bigl[\exp(c_{\alpha}|X|)\bigr] = 1+o(1), \qquad\mbox{for each } X\in\{S_1, S_{21}, S_{22}, S_{31},\ldots,S_{34}, S_{41}, S_{42}\}. 
\eeq 
Below, we show \eqref{LB-eq2} for each of $X$ listed above.


First, consider $X=S_1$. Let $\delta_1,\delta_2,\ldots,\delta_K$ be the $K$ eigenvalues of $M$, arranged in the descending order of magnitude, and let $b_1,b_2,\ldots,b_K$ be the associated eigenvectors. Then, $M=\sum_{k=1}^K \delta_kb_kb_k'$. 
It follows that 
\[
    S_1 =\alpha^{-1}\sum_{k,l}\delta_k\delta_l\left(\sum_i(y_i'b_k)(\tilde{y}_i'b_l)\right)^2-\alpha^{-1}\sum_{k,l}\delta_k\delta_l\sum_i(y_i'b_k)^2(\tilde{y}_i'b_l)^2. 
\]
Note that  $\max_k|\delta_k|=\|M\|$, where $\|M\|$ is the operator norm of $M$. Therefore, 
\[
    |S_1|\leq\alpha^{-1}\|M\|^2\left[\sum_{k,l}\left(\sum_i(y_i'b_k)(\tilde{y}_i'b_l)\right)^2+\sum_{k,l}\sum_i(y_i'b_k)^2(\tilde{y}_i'b_l)^2\right].
\]
In addition, for any $i\in\llbracket1,n\rrbracket$ and $k\in\llbracket1,K\rrbracket$, by the Cauchy-Schwarz inequality, we have $(y_i'b_k)^2\leq \|y_i\|^2_2\|b_k\|^2_2=\|y_i\|^2_2\leq \|y_i\|_1\leq 2 $, given that $\|y_i\|_\infty\leq 1$ and that $\|y_i\|_1\leq \|\pi_i\|_1+\|h\|_1\leq 2$. It follows that
\begin{equation} \label{LB-eq4(0)}
    |S_1|\leq  4n\alpha^{-1}K^2\|M\|^2+R_1,  
\end{equation}
where
\[
R_1\equiv \alpha^{-1}K^2\|M\|^2\max_{k,l}\left(\sum_i(y_i'b_k)(\tilde{y}_i'b_l)\right)^2. 
\]
To bound $R_1$, 
we fix a tuple $(k,l)$ and provide an upper bound for $Y_{kl}:=\sum_i(y_i'b_k)(\tilde{y}_i'b_l)$. Note that $Y_{kl}$ is a sum of independent, zero-mean random variables. In addition, $|(y_i'b_k)(\tilde{y}_i'b_l)|\leq\|y_i\|_2\|\tilde{y}_i\|_2\leq 2$. We can apply Hoeffding's inequality, for any $t>0$:
\[
    \mathbb{P}(|Y_{kl}|>t)\leq 2\exp\left(-\frac{2t^2}{\sum_{i=1}^n(2\|y_i\|_2\|\tilde{y}_i\|_2)^2}\right)=2\exp\left(-\frac{t^2}{8n}\right). 
\]
Hence, denoting $Y_*:=\max_{k,l}|Y_{kl}|$, we have
\begin{equation*} 
    \mathbb{P}\left(Y_*>t\right) = \mathbb{P}\left(\bigcup_{k,l}\{|Y_{kl}|>t\}\right)\leq \sum_{k,l}\mathbb{P}(|Y_{kl}|>t)\leq 2K^2\exp\left(-\frac{t^2}{8n}\right). 
\end{equation*}
It follows that, for any $t>0$, 
\begin{equation}  \label{LB-eq4}
    \mathbb{P}(R_1>t)=\mathbb{P}\left(Y^*>\frac{\sqrt{t\alpha}}{K\|M\|}\right)\leq 2K^2\exp\left(-\frac{\alpha t}{8nK^2\|M\|^2}\right).
\end{equation}
We now use \eqref{LB-eq4(0)} and \eqref{LB-eq4} to bound $\mathbb{E}[\exp(c_\alpha |S_1|)]$. For any non-negative variable $X$, it follows from integration by part that $\mathbb{E}[\exp(X)]=1+\int_0^\infty e'\mathbb{P}(X>t)dt$. It follows that
\begin{align*}
\mathbb{E}\bigl[\exp(c_\alpha |S_1|)\bigr] & \leq e^{4c_{\alpha} n\alpha^{-1}K^2\|M\|^2} \cdot \mathbb{E}\bigl[\exp(c_\alpha R_1)\bigr]\cr
&\leq e^{4c_\alpha n\alpha^{-1}K^2\|M\|^2} \left[1+\int_0^\infty e'\mathbb{P}\left(R_1>c_{\alpha}^{-1}t\right)dt\right]\cr
&\leq e^{4c_{\alpha} n\alpha^{-1}K^2\|M\|^2}\left[ 1+ \int_0^{\infty} 2e^{-\bigl(\frac{\alpha}{8c_\alpha n K^2 \|M\|^2}-1\bigr)t}\textcolor{black}{dt}\right]. 
\end{align*}
In our assumption, $\beta_n\to 0$, which implies that 
\[
n\alpha^{-1}\|M\|^2\to 0.
\]
It follows that  $e^{4c_\alpha n\alpha^{-1}K^2\|M\|^2}=\exp(o(1))=1+o(1)$. 
Also, for $n$ big enough, $\frac{\alpha}{8c_\alpha nK^2\|M\|^2}-1>0$. Furthermore, we note that for any value $z>0$, $\int_0^\infty e^{-zt}dt=z^{-1}$. Combining the above gives
\beq \label{LBproof-S1}
\mathbb{E}\bigl[\exp(c_\alpha |S_1|)\bigr] \leq  e^{4c_{\alpha}K^2 n\alpha^{-1}\|M\|^2}\left( 1+ \frac{16c_\alpha K^2 n\alpha^{-1}\|M\|^2}{1-8c_\alpha K^2 n\alpha^{-1}\|M\|^2}\right)=1+o(1). 
\eeq
This proves \eqref{LB-eq2} for $X=S_1$.

Second, consider $X=S_{21}$ (the analysis of $S_{22}$ is similar and thus omitted). We define a unit-norm vector $u=\|Mh\|^{-1}(Mh)$. Then,  
\[
    S_{21} = \alpha^{-1}(n-1)\|Mh\|^2\sum_{i}(y_i'u)(\tilde{y}_i'u). 
\]
The variables $\{(y_i'u)(\tilde{y}_i'u)\}_{1\leq i\leq n}$ are independent, with $|(y_i'u)(\tilde{y}_i'u)|\leq \|y_i\|\|\tilde{y}_i\|$. We have seen that $\|y_i\|^2\leq 2$ and $\|\tilde{y}_i\|^2\leq 2$. It follows that $|(y_i'u)(\tilde{y}_i'u)|\leq 2$. Applying Hoeffding's inequality, we obtain that, for any $t>0$,
\begin{align*} 
    \mathbb{P}(|S_{21}|>t) &= \mathbb{P}\left(\left|\sum_{i}(y_iu)(\tilde{y}_iu)\right|>\frac{t\alpha}{(n-1)\|Mh\|^2}\right)\cr
    &\leq 2\exp\left(-\frac{t^2\alpha^2}{8n(n-1)^2\|Mh\|^4}\right)\leq 2\exp\left(-\frac{t^2\alpha^2}{8n^3\|Mh\|^4}\right).
\end{align*}
Our assumption $\beta_n\to 0$ implies that 
\[
n^{3}\alpha^{-2}\|Mh\|^4\to 0.
\]
Furthermore, for $z>0$, we have $\int_0^\infty e^{-zt^2+t}\leq \sqrt{2\pi z^{-1}} e^{(4z)^{-1}}$. Combining these gives
\begin{align}   \label{LBproof-S2}
\mathbb{E}\bigl[\exp(c_\alpha |S_{21}|)\bigr] &= 1+\int_0^\infty e'\mathbb{P}\left(|S_{21}|>c_{\alpha}^{-1}t\right)dt \cr
&\leq 1+ \int_0^{\infty} 2e^{-\frac{\alpha^2}{8\textcolor{black}{c_\alpha^2} n^3\|Mh\|^4}t^2 + \, t}dt\cr
&\textcolor{black}{\leq 1+2\sqrt{2\pi}\sqrt{8c_{\alpha}^2n^3\alpha^{-2}\|Mh\|^4}\exp\left(-2c_\alpha^2 n^3\alpha^{-2}\|Mh\|^4\right)}\cr
&=1+o(1). 
\end{align}
This proves \eqref{LB-eq2} for $X=S_{21}$.

Next, consider $S_{31}$ (the analyses of $S_{32}$-$S_{34}$ are similar and omitted). Recall that $M=\sum_{k=1}^K\delta_kb_kb_k'$ is the eigen-decomposition of $M$; additionally, we have defined $u=\|Mh\|^{-1}(Mh)$. It follows that   
\begin{align*}
    S_{31} &= \alpha^{-1}\|Mh\|\sum_{i\neq j}(y_i'My_j)(\tilde{y}_i'u)\cr
    &= \alpha^{-1}\|Mh\|\sum_{i\neq j}\biggl[\sum_{k} \delta_k (y_i'b_k)(y_j'b_k)\biggr](\tilde{y}_i'u)\cr 
    &= \alpha^{-1}\|Mh\|\sum_{k}\delta_k \biggl[ \sum_{i}(y_i'b_k)(\tilde{y}_i'u)\biggr]\biggl[\sum_{j}(y_j'b_k)\biggr]\\
    &\qquad -\alpha^{-1}\|Mh\|\sum_{k}\delta_k\biggl[ \sum_{i}(y_i'b_k)^2(\tilde{y}_i'u)\biggr]. 
\end{align*}
We have seen that $\|b_i\|^2=1$, $\|y_i\|^2\leq 2$, $\|\tilde{y}_i\|\leq 2$, $\|u\|=1$, and $|\delta_k|\leq \|M\|$.  
It follows that
\beq \label{LB-eq5}
    |S_{31}| \leq R_{31} + 2\sqrt{2}\, n\alpha^{-1}K\|M\|\|Mh\|, 
\eeq
where
\[
R_{31}:=\alpha^{-1}\|M\|\|Mh\|K\max_{k}Z_{k},\qquad\mbox{with}\quad    Z_k =\biggl[ \sum_{i}(y_i'b_k)(\tilde{y}_i'u)\biggr]\biggl[\sum_{j}(y_j'b_k)\biggr].
\]
We can derive the tail probability bound for $Z_k$: Since $|y_i'b_k|\leq \|y_i\|\leq\sqrt{2}$ and $|\tilde{y}_i'u|\leq \|\tilde{y}_i\|\leq \sqrt{2}$, the Hoeffding's inequality yields that 
\begin{align*}
    \mathbb{P}(|Z_k|>t)&\leq \mathbb{P}\left( \left|\sum_{i}(y_i'b_k)(\tilde{y}_i'u)\right|>\sqrt{t}\right)+\mathbb{P}\left(\left|\sum_{j}(y_j'b_k)\right|>\sqrt{t}\right)\\
    &\leq 2\exp\left(-\frac{t}{8n}\right)+2\exp\left(-\frac{t}{4n}\right) \leq 4\exp\left(-\frac{t}{\textcolor{black}{8}n}\right).
\end{align*}
We thus have
\begin{equation} \label{LB-eq6}
    \mathbb{P}(|R_{31}|>t)=\mathbb{P}\left(\max_{k}Z_{k}>\frac{t\alpha}{K\|M\|\|Mh\|}\right)\leq 4K\exp\left(-\frac{t\alpha}{\textcolor{black}{8}nK\|M\|\|Mh\|}\right). 
\end{equation}
We apply \eqref{LB-eq5}-\eqref{LB-eq6} to bound $\mathbb{E}[\exp(c_\alpha|S_{31}|)]$. Our assumption $\beta_n\to 0$ ensures that $n\alpha^{-1}\|M\|^2~\to~0$. Note that $\|Mh\|\leq \|M\|\|h\|\leq \|M\|\sqrt{\|h\|_1\|h\|_\infty}\leq \|M\|$. It follows that 
\[
n\alpha^{-1}\|M\|\|Mh\|\to 0.
\]
We then mimic the proof of \eqref{LBproof-S1} to get
\begin{align} \label{LBproof-S3}
\mathbb{E}\bigl[\exp(c_\alpha |S_{31}|)\bigr] & \leq e^{2\sqrt{2}c_\alpha n\alpha^{-1}K\|M\|\|Mh\|}  \left[1+\int_0^\infty e'\mathbb{P}\left(|R_{31}|>c_{\alpha}^{-1}t\right)dt\right]\cr
&\leq e^{2\sqrt{2}c_\alpha n\alpha^{-1}K\|M\|\|Mh\|}  \left[ 1+ \int_0^{\infty} 4K e^{-\bigl( \frac{\alpha}{4c_\alpha nK\|M\|\|Mh\|}-1\bigr)t}\right]\cr
&\leq  e^{2\sqrt{2}c_\alpha n\alpha^{-1}K\|M\|\|Mh\|} \left( 1+ \frac{16c_\alpha K^2 n\alpha^{-1}\|M\|\|Mh\|}{1-4c_\alpha K n\alpha^{-1}\|M\|\|Mh\|}\right)\cr
&=1+o(1). 
\end{align}
This proves \eqref{LB-eq2} for $X=S_{31}$. 

Last, consider $S_{41}$ (the analysis of $S_{42}$ is similar and omitted). Since $u=\|Mh\|^{-1}Mh$, we have
\begin{align*}
     S_{41} & = \alpha^{-1}\|Mh\|^2\sum_{i\neq j} (y_i'u)(\tilde{y}_j'u)\cr
     &=  \alpha^{-1}\|Mh\|^2\biggl[ \sum_{i}(y_i'u)\biggr] \biggl[ \sum_j (\tilde{y}_j'u)\biggr] - \alpha^{-1}\|Mh\|^2\sum_{i}(y_i'u)(\textcolor{black}{\tilde{y}_i}'u). 
\end{align*}
Note that $|(y_i'u)(\textcolor{black}{\tilde{y}_i}'u)|\leq \|y_i\|\|\textcolor{black}{\tilde{y}_i}\|\leq 2$. We immediately have
\beq \label{LB-eq8(0)}
|S_{41}|\leq R_{41}+ 2n\alpha^{-1}\|Mh\|^2,
\eeq
where 

\[
R_{41} =  \alpha^{-1}\|Mh\|^2\biggl[ \sum_{i}(y_i'u)\biggr] \biggl[ \sum_j (\tilde{y}_j'u)\biggr]. 
\]
We apply Hoeffding's inequality to derive the tail probability bound: For all $t>0$,
\begin{align} \label{LB-eq8}
    \mathbb{P}(|R_{41}|>t) &= \mathbb{P}\left( \left|\sum_{i}(y_i'u) \right| >\frac{\sqrt{\alpha t}}{\|Mh\|}   \right)+ \mathbb{P}\left( \left|\sum_{j}(\tilde{y}_j'u) \right| >\frac{\sqrt{\alpha t}}{\|Mh\|}   \right)\cr
    &\leq 4\exp\left(-\frac{\alpha t}{8n\|Mh\|^2}\right). 
\end{align}
We have seen that $\|Mh\|\leq  \|M\|$. Therefore, the assumption of $\beta_n\to 0$ leads to 
\[
n\alpha^{-1}\|Mh\|^2\to 0.
\] 
Using \eqref{LB-eq8(0)} and \eqref{LB-eq8}, we have
\begin{align*}
\mathbb{E}\bigl[\exp(c_\alpha |S_{41}|)\bigr] 
&\leq e^{2c_\alpha n\alpha^{-1}\|Mh\|^2} \left[1+\int_0^\infty e'\mathbb{P}\left(|R_{41}|>c_{\alpha}^{-1}t\right)dt\right]\cr
&\leq e^{2c_\alpha n\alpha^{-1}\|Mh\|^2} \left[ 1+ \int_0^{\infty} \textcolor{black}{4}e^{-\bigl(\frac{\alpha}{8c_\alpha n\|Mh\|^2}-1\bigr)t}\right]. \cr
& \leq  e^{2c_\alpha n\alpha^{-1}\|Mh\|^2}\left( 1+ \frac{\textcolor{black}{32}c_\alpha  n\alpha^{-1}\|Mh\|^2}{1-8c_\alpha n\alpha^{-1}\|Mh\|^2}\right)\cr
&=1+o(1). 
\end{align*}
This proves \eqref{LB-eq2} for $X=S_{41}$. \qed

\section{Proof of Theorem~\ref{thm:minimax}}\label{appendix:proof_minimax}

Note: this proof requires Lemma~\ref{lm:minimax_0}  and Lemma~\ref{lm:minimax_1}, which are provided directly after the proof.\\

We start by studying the case $t_0=0$. We consider a sequence of null hypotheses indexed by $n$, where $\Omega_n=\alpha_n\bm{1}_K\bm{1}_K'\in\mathcal{M}_{0n}$ under $H_0^{(n)}$. For our sequence of alternatives, we consider $\Omega_n = \Pi_nP_n\Pi_n'$ under $H_1^{(n)}$, with
\begin{align*}
    P_n = \alpha_n\left[\gamma_nI_K+(1-\gamma_n)\bm{1}_K\bm{1}_K'\right], \quad \mbox{and} \quad \pi_1,...,\pi_n &\overset{\mbox{iid}}{\sim} F,
\end{align*}
where for all $k\in\{1,...,K\}$,
\begin{equation*}
    \mathbb{P}_{\pi\sim F}(\pi=e_k) = \frac{1}{K}.
\end{equation*}
In the above definition, $\{e_k\}_{k=1}^K$ denotes the canonical basis of $\mathbb{R}^K$. It follows that
\begin{equation*}
    h := \E_{\pi\sim F}[\pi] = \frac{1}{K}\bm{1}_K, \quad \mbox{and} \quad \Sigma := \E_{\pi\sim F}[\pi\pi'] = \frac{1}{K}I_K. 
\end{equation*}
Under this random mixed membership model, it is straightforward to verify that
\begin{align*}
    \alpha_0 &= \alpha_n\left(1-\frac{K-1}{K}\gamma_n\right),\\
    \|P_nh-\alpha_0\bm{1}_K\| &= 0,\\
    \|P_n-\alpha_0\bm{1}_K\bm{1}_K'\| &=  \alpha_n\gamma_n.
\end{align*}
Hence
\begin{equation*}
    \beta_n = \max\bigl\{ n^{3/2}\alpha_0^{-1}\|P_n h-\alpha_{0}{\bf 1}_K\|^2,\quad  n^2\alpha_0^{-2}\|P_n - \alpha_{0}{\bf 1}_K{\bf 1}_K' \|^4  \bigr\} = n^2\alpha_0^{-2}\alpha_n^4\gamma_n^4.
\end{equation*}
By assumption, $\gamma_n\goto0$, hence for $n$ sufficiently large, $\alpha_n <2\alpha_0$, hence
\begin{equation*}
    \beta_n=O(n^2\alpha_n^2\gamma_n^4) = o(1),
\end{equation*}
under the assumption that $n^2\alpha_n^2\gamma_n^4 = o(1)$. By Theorem~\ref{thm:LB}, the $\chi^2$-distance between the two distributions satisfies $D_{\chi^2}(f_0^{(n)}\| f_1^{(n)})=o(1)$. By connection between $L_1$-distance and $\chi^2$-distance, it follows that
\begin{equation*}
    \|f_0^{(n)}-f_1^{(n)}\|_1=o(1).
\end{equation*}
We now slightly modify the alternative hypothesis. Let $\{\Pi^0_n\}_n$ be a sequence of non-random membership matrices such that $(P_n,\Pi^0_n)\in\mathcal{M}_{1n}(0)$. Such a sequence can be built e.g. by considering $\lfloor n/K\rfloor$ pure nodes in each community and all other nodes equally mixed across all communities. In the modified alternative hypothesis $\tilde{H}_1^{(n)}$,
\[
\tilde{\Pi} =
\begin{cases}
\Pi_n, & \quad \mbox{if } (\Pi_n,P_n)\in\mathcal{M}_{1n}(0),\\
 \Pi^0_{n}, & \quad \mbox{otherwise}. 
\end{cases}
\]
Let $\tilde{f}_1^{(n)}$ be the probability measure associated with $\tilde{H}_1^{(n)}$. Under $\tilde{H}_1^{(n)}$, all realizations $\tPi_nP_n\tPi_n'$ are in the class $\mathcal{M}_{1n}(0)$, by definition. By the Neyman-Pearson lemma and elementary inequalities,
\begin{align*}
    Risk_n^*(0) &\geq 1 - \inf_{f_0\in\mathcal{M}_{0n},f_1\in\mathcal{M}_{1n}(0)}\{\|f_0-f_1\|_1\}\\
    &\geq 1-\|f_0^{(n)}-\tilde{f}_1^{(n)}\|_1\\
    &\geq 1-\|f_0^{(n)}-f_1^{(n)}\|_1-\|f_1^{(n)}-\tilde{f}_1^{(n)}\|_1\\
    &\geq 1-o(1)-\|f_1^{(n)}-\tilde{f}_1^{(n)}\|_1.
\end{align*}
It follows from Lemma~\ref{lm:minimax_0} that $\tilde{\Pi}_n=\Pi_n$ with probability $1-o(1)$. As a result,
\begin{equation*}
    \|f_1^{(n)}-\tilde{f}_1^{(n)}\|_1=o(1),
\end{equation*}
from which we obtain that $\lim_{n\to\infty}\{ Risk_n^*(0)\} = 1$.\\

Next, we study the case $0<t_0$. Again, we consider a sequence of null hypotheses indexed by $n$, where $\Omega_n=\alpha_n\bm{1}_K\bm{1}_K'\in\mathcal{M}_{0n}$ under $H_0^{(n)}$. For our sequence of alternatives, we consider $\Omega_n = \Pi_nP_n\Pi_n'$ under $H_1^{(n)}$, with
\begin{align*}
    P_n = \alpha_n\left[\gamma_nI_K+(1-\gamma_n)\bm{1}_K\bm{1}_K'\right], \quad \mbox{and} \quad \pi_1,...,\pi_n &\overset{\mbox{iid}}{\sim} F,
\end{align*}
where
\begin{equation*}
    \mathbb{P}_{\pi\sim F}(\pi=e_1) = \frac{K+1}{2K}, \quad \mbox{and}\quad \mathbb{P}_{\pi\sim F}(\pi=e_1) = \frac{1}{2K} \quad  \forall k\in\{2,...,K\}.
\end{equation*}
It follows that
\begin{equation*}
    h := \E_{\pi\sim F}[\pi] = \frac{1}{2K}(Ke_1+\bm{1}_K), \quad \mbox{and} \quad \Sigma := \E_{\pi\sim F}[\pi\pi'] = \frac{1}{2K}(Ke_1e_1'+I_K). 
\end{equation*}
Under this random mixed membership model, it is straightforward to verify that
\begin{align*}
    \alpha_0 &= \alpha_n\left(1-\frac{3K-3}{4K}\gamma_n\right),\\
    \|P_nh-\alpha_0\bm{1}_K\| &= \alpha_n\gamma_n\sqrt{\frac{(K-1)(K+3)}{16K}},\\
    \|P_n-\alpha_0\bm{1}_K\bm{1}_K'\| &= \max\left\{ \alpha_n\gamma_n, \frac{K-1}{4}\alpha_n\gamma_n  \right\}.
\end{align*}
Recall that
\begin{equation*}
    \beta_n = \max\bigl\{ n^{3/2}\alpha_0^{-1}\|P_n h-\alpha_{0}{\bf 1}_K\|^2,\quad  n^2\alpha_0^{-2}\|P_n - \alpha_{0}{\bf 1}_K{\bf 1}_K' \|^4  \bigr\}.
\end{equation*}
Hence
\begin{equation*}
    \beta_n = \max\left\{ n^{3/2}\alpha_0^{-1}\alpha_n^2\gamma_n^2\frac{(K-1)(K+3)}{16K},\quad  \max \left(1,\left(\frac{K-1}{4}\right)^4\right)n^2\alpha_0^{-2}\alpha_n^4\gamma_n^4  \right\}
\end{equation*}
By assumption, $\gamma_n\goto0$, hence for $n$ sufficiently large, $\alpha_n <2\alpha_0$, hence
\begin{equation*}
    \beta_n=O\left(\max\left\{ n^{3/2}\alpha_n\gamma_n^2,\quad  n^2\alpha_n^2\gamma_n^4  \right\}\right) = O\left(\max\left\{ n^{3/2}\alpha_n\gamma_n^2,\quad  \frac{(n^{3/2}\alpha_n\gamma_n^2)^2}{n}  \right\}\right) = o(1),
\end{equation*}
under the assumption that $n^{3/2}\alpha_n\gamma_n^2=o(1)$. By Theorem~\ref{thm:LB}, the $\chi^2$-distance between the two distributions satisfies $D_{\chi^2}(f_0^{(n)}\| f_1^{(n)})=o(1)$. By connection between $L_1$-distance and $\chi^2$-distance, it follows that
\begin{equation*}
    \|f_0^{(n)}-f_1^{(n)}\|_1=o(1).
\end{equation*}
We now slightly modify the alternative hypothesis. Let $\{\Pi^0_n\}_n$ be a sequence of non-random membership matrices such that $(P_n,\Pi^0_n)\in\mathcal{M}_{1n}(t_0)$. Such a sequence can be built e.g. by considering $\lfloor n(K+1)/2K\rfloor$ pure nodes in community 1, $\lfloor n/2K\rfloor$ nodes in communities $2$ to $K$ and all other nodes with mixed membership vector $(2K)^{-1}(Ke_1+\bm{1}_K)$. In the modified alternative hypothesis $\tilde{H}_1^{(n)}$,
\[
\tilde{\Pi} =
\begin{cases}
\Pi_n, & \quad \mbox{if } (\Pi_n,P_n)\in\mathcal{M}_{1n}(t_0),\\
 \Pi^0_{n}, & \quad \mbox{otherwise}. 
\end{cases}
\]
Let $\tilde{f}_1^{(n)}$ be the probability measure associated with $\tilde{H}_1^{(n)}$. Under $\tilde{H}_1^{(n)}$, all realizations $\tPi_nP_n\tPi_n'$ are in the class $\mathcal{M}_{1n}(t_0)$, by definition. By the Neyman-Pearson lemma and elementary inequalities,
\begin{align*}
    Risk_n^*(0) &\geq 1 - \inf_{f_0\in\mathcal{M}_{0n},f_1\in\mathcal{M}_{1n}(t_0)}\{\|f_0-f_1\|_1\}\\
    &\geq 1-\|f_0^{(n)}-\tilde{f}_1^{(n)}\|_1\\
    &\geq 1-\|f_0^{(n)}-f_1^{(n)}\|_1-\|f_1^{(n)}-\tilde{f}_1^{(n)}\|_1\\
    &\geq 1-o(1)-\|f_1^{(n)}-\tilde{f}_1^{(n)}\|_1.
\end{align*}
It follows from Lemma~\ref{lm:minimax_1} that $\tilde{\Pi}_n=\Pi_n$ with probability $1-o(1)$. As a result,
\begin{equation*}
    \|f_1^{(n)}-\tilde{f}_1^{(n)}\|_1=o(1),
\end{equation*}
from which we obtain that $\lim_{n\to\infty}\{ Risk_n^*(t_0)\} = 1$. \qed

\begin{lemma}[Case $t_0=0$] \label{lm:minimax_0}
Fix $K\geq 2$, a sequence $\{\alpha_n\}_n\in[0,1]^{\mathbb{N}}$, and a sequence $\{\gamma_n\}_n\in(\mathbb{R}_+)^{\mathbb{N}}$. Denote by $\{e_k\}_{k=1}^K$ the canonical basis of $\mathbb{R}^K$. Consider the sequence of alternative probability matrices $\Omega_n = \Pi_nP_n\Pi_n'$, with
\begin{align*}
    P_n = \alpha_n\left[\gamma_nI_K+(1-\gamma_n)\bm{1}_K\bm{1}_K'\right], \quad \mbox{and} \quad \pi_1,...,\pi_n &\overset{\mbox{iid}}{\sim} F,
\end{align*}
where for all $k\in\{1,...,K\}$, $\mathbb{P}_{\pi\sim F}(\pi=e_k) = \frac{1}{K}$. Suppose that $\alpha_n\goto0$, $n\alpha_n\goto\infty$, and $\gamma_n\goto0$. Then, with probability $1-o(1)$, $(P_n,\Pi_n)\in\mathcal{M}_{1n}(0)$.
\end{lemma}

\noindent\textit{Proof}\\
From the proof of Theorem~\ref{thm:minimax} for $t_0=0$, we know that
\begin{align*}
    h := \E_{\pi\sim F}[\pi] = \frac{1}{K}\bm{1}_K, \quad \Sigma := \E_{\pi\sim F}[\pi\pi'] = \frac{1}{K}I_K'\quad \mbox{and} \quad \alpha_0 &= \alpha_n\left(1-\frac{K-1}{K}\gamma_n\right).
\end{align*}
We introduce the following random quantities:
\begin{equation*}
    \tH = \frac{1}{n}\sum_{i=1}^n\pi_i,\quad \tG=\frac{1}{n}\sum_{i=1}^n\pi_i\pi_i', \quad \mbox{and} \quad \tilde{\alpha}_0=\tH P_n\tH'.
\end{equation*}
To show that $(P,\Pi)\in\mathcal{M}_{1n}(0)$, we will check that 
\begin{enumerate}
    \item $OSC(\tH)\leq C$ and $\|\tG^{-1}\|\leq C$,
    \item $\tilde{\alpha}_0\leq c$, $n\tilde{\alpha}_0\geq c^{-1}$, and $\tilde{\alpha}_0\geq\alpha_n/2$,
    \item $2\tilde{\alpha}_0^{-1}\|P_n-\tilde{\alpha}_0{\bf 1}_K{\bf 1}_K'\|\geq \gamma_n$.
\end{enumerate}

First, recognize that $\tH=n^{-1}\sum_{i=1}^n\tpi_i\xrightarrow[]{\mbox{as}}K^{-1}\bm{1}_K$ by the Strong Law of Large Numbers. As a consequence, for $n$ sufficiently large, we have $OSC(\tH)<C$ with probability at least $1-o(1)$. Next, let $y_i=\pi_i-\tH$. We have
\begin{align*}
    n\tG &= \sum_{i=1}^n\pi_i\pi_i' = \sum_{i=1}^n (\tilde{h}\tilde{h}'+\tilde{h}y_i'+y_i\tilde{h}'+y_iy_i')\\
    &= n\Sigma + \sum_{i=1}^n (y_iy_i'-\E[y_iy_i']) + \sum_{i=1}^n (\tilde{h}y_i') + \sum_{i=1}^n (y_i\tilde{h}')\\
    &= n\Sigma+Z_0+Z_1+Z_2.
\end{align*}
Notice that $Z_0$ is a sum of $n$ independent mean-zero random matrices, so we can apply the matrix Hoeffding inequality to bound its operator norm. Since $\|y_iy_i'-\E[y_iy_i']\|\leq C$, we obtain for $t>0$,
\begin{equation*}
    \mathbb{P}\left(\|Z_0\|>t\right)\leq \exp\left(-\frac{Ct^2}{n}\right).
\end{equation*}
If we pick $t=C\sqrt{n\log(n)}$, then we have that $\|Z_0\|<C\sqrt{n\log(n)}$ with probability $1-o(1)$. Similarly, it is straightforward to show that $\|Z_1+Z_2\|\leq C\sqrt{n\log(n)}$ with probability $1-o(1)$. Now, recall that $\lambda_{\min}(\Sigma)=K^{-1}$. As a result,
\begin{equation*}
    \lambda_{\min}(n\tG) = \lambda_{\min}(n\Sigma+Z_0+Z_1+Z_2) > \lambda_{\min}(n\Sigma)-\|Z_0+Z_1+Z_2\| > \frac{n}{K}-C\sqrt{n\log(n)}.
\end{equation*}
It follows that
\begin{equation*}
    \lambda_{\min}(\tG) > \frac{1}{K}-C\sqrt{\frac{\log(n)}{n}},
\end{equation*}
which shows that for $n$ sufficiently large, $\|\tG^{-1}\|<C$ with probability $1-o(1)$.\\

Next, we show that $\tilde{\alpha}_0 < c$ and $n\tilde{\alpha}_0>c^{-1}$ with high probability. Denote $z:=\tilde{h}-h$. We can rewrite
\begin{equation*}
    \tilde{\alpha}_0 = \tH'P\tH = z'Pz+2h'Pz+\alpha_0.
\end{equation*}
Notice that both $\|z'Pz\|\leq C\|Pz\|$ and $\|h'Pz\|\leq C\|Pz\|$. We now provide a high-probability bound on the 2-norm of $Pz$, which can be written as a sum of mean-zero independent random variables
\begin{equation*}
    Pz = \frac{1}{n}\sum_{i=1}^n (P\pi_i-Ph),
\end{equation*}
where for all $i=1,...,n$ it holds that $\|P\pi_i-Ph\|\leq C\|P\|\leq C\alpha_n\gamma_n$. For $t>0$, Hoeffding's inequality yields
\begin{equation*}
    \mathbb{P}\left(\|Pz\|>t\right) \leq C\exp\left(-\frac{Cnt^2}{\alpha_n^2\gamma_n^2}\right).
\end{equation*}
Pick $t=\alpha_n\gamma_n\sqrt{\log(n)/n}$. As a consequence, we obtain that $\|Pz\|<\alpha_n\gamma_n\sqrt{\log(n)/n}$ with probability $1-o(1)$. Hence with probability $1-o(1)$,
\begin{equation*}
    \tilde{\alpha}_0 = \alpha_0+O\left(\alpha_n\gamma_n\sqrt{\frac{\log(n)}{n}}\right) = \alpha_n-\frac{K-1}{K}\alpha_n\gamma_n+o\left(\alpha_n\gamma_n\right) = \alpha_n+O(\alpha_n\gamma_n).
\end{equation*}
It follows that for $n$ sufficiently large, $\tilde{\alpha}_0<c$ and $n\tilde{\alpha}_0>c^{-1}$ with probability $1-o(1)$. We also obtain from this last equation that for $n$ sufficiently large, $\tilde{\alpha}_0\geq \alpha_n/2$ with probability $1-o(1)$.\\

It remains to show that $2\tilde{\alpha}_0^{-1}\|P_n-\tilde{\alpha}_0{\bf 1}_K{\bf 1}_K'\|\geq \gamma_n$. With probability $1-o(1)$, the matrix $(P_n-\tilde{\alpha}_0{\bf 1}_K{\bf 1}_K')$ has eigenvalues 
\begin{align*}
    \lambda_+ &= K(\alpha_n-\tilde{\alpha}_0)-(K-1)\alpha_n\gamma_n = o(\alpha_n\gamma_n),\\
    \lambda_- &= \alpha_n\gamma_n.
\end{align*}
hence for $n$ sufficiently large, we must have $\|P_n-\tilde{\alpha}_0{\bf 1}_K{\bf 1}_K'\|=\alpha_n\gamma_n\geq \tilde{\alpha}_0\gamma_n/2$. It follows that 
\begin{equation*}
    2\tilde{\alpha}_0^{-1}\|P-\tilde{\alpha}_0{\bf 1}_K{\bf 1}_K'\| \geq \gamma_n,
\end{equation*}
which concludes the proof. \qed

\begin{lemma}[Case $0 < t_0$] \label{lm:minimax_1}
Fix $K\geq 2$, a sequence $\{\alpha_n\}_n\in[0,1]^{\mathbb{N}}$, and a sequence $\{\gamma_n\}_n\in(\mathbb{R}_+)^{\mathbb{N}}$. Denote by $\{e_k\}_{k=1}^K$ the canonical basis of $\mathbb{R}^K$. Consider the sequence of alternative probability matrices $\Omega_n = \Pi_nP_n\Pi_n'$, with
\begin{align*}
    P_n = \alpha_n\left[\gamma_nI_K+(1-\gamma_n)\bm{1}_K\bm{1}_K'\right], \quad \mbox{and} \quad \pi_1,...,\pi_n &\overset{\mbox{iid}}{\sim} F,
\end{align*}
where 
\begin{equation*}
    \mathbb{P}_{\pi\sim F}(\pi=e_1) = \frac{K+1}{2K}, \quad \mbox{and}\quad \mathbb{P}_{\pi\sim F}(\pi=e_k) = \frac{1}{2K} \quad  \forall k\in\{2,...,K\}.
\end{equation*}
Suppose that $\alpha_n\goto0$, $n\alpha_n\goto\infty$, $\gamma_n\goto0$, and $0<t_0<\sqrt{(K-1)(K+3)/(16K)}$. Then, with probability $1-o(1)$, $(P_n,\Pi_n)\in\mathcal{M}_{1n}(t_0)$.
\end{lemma}

\noindent\textit{Proof}\\
From the proof of Theorem~\ref{thm:minimax} for $t_0>0$, we know that
\begin{align*}
    h :&= \E_{\pi\sim F}[\pi] = \frac{1}{2K}(Ke_1+\bm{1}_K),\\
    \Sigma :&= \E_{\pi\sim F}[\pi\pi'] = \frac{1}{2K}(Ke_1e_1'+I_K),\\
    \mbox{and} \quad \alpha_0 &= \alpha_n\left(1-\frac{3K-3}{4K}\gamma_n\right).
\end{align*}
We introduce the following random quantities:
\begin{equation*}
    \tH = \frac{1}{n}\sum_{i=1}^n\pi_i,\quad \tG=\frac{1}{n}\sum_{i=1}^n\pi_i\pi_i', \mbox{and} \quad \tilde{\alpha}_0=\tH P_n\tH'.
\end{equation*}
To show that $(P,\Pi)\in\mathcal{M}_{1n}(t_0)$, we will check that 
\begin{enumerate}
    \item $OSC(\tH)\leq C$ and $\|\tG^{-1}\|\leq C$,
    \item $\tilde{\alpha}_0 \leq c$, $n\tilde{\alpha}_0\geq c^{-1}$, and $\tilde{\alpha}_0\geq \alpha_n/2$,
    \item $2\tilde{\alpha}_0^{-1}\|P-\tilde{\alpha}_0{\bf 1}_K{\bf 1}_K'\|\geq \gamma_n$ and $\|P\tilde{h}-\tilde{\alpha}_0\bm{1}_K\|\geq t_0\|P-\tilde{\alpha}_0{\bf 1}_K{\bf 1}_K'\|$.
\end{enumerate}
The first two points can be shown with probability at least $1-o(1)$ in the same way as in the proof of Lemma~\ref{lm:minimax_0}. We will focus on the third point. For $n$ sufficiently large, $\alpha_n\gamma_n$ is the largest eigenvalue of $(P-\tilde{\alpha}_0{\bf 1}_K{\bf 1}_K')$ in magnitude. Hence, we must have, for $n$ sufficiently big
\begin{equation*}
    \|P_n-\tilde{\alpha}_0{\bf 1}_K{\bf 1}_K'\| = \alpha_n\gamma_n \geq \tilde{\alpha}_0\gamma_n/2.
\end{equation*}
Now, introduce the (continuous) function with support $\mathbb{R}^K$:
\begin{equation*}
    g(x):=\left\|\begin{bmatrix}
                x_1(1-x_1)-\sum_{k\neq 1}x_k^2 \\
                x_2(1-x_2)-\sum_{k\neq 2}x_k^2  \\
                \vdots \\
                x_K(1-x_K)-\sum_{k\neq K}x_k^2 
                \end{bmatrix}\right\|.
\end{equation*}
Notice that $\|P\tH-\tilde{\alpha}_0\bm{1}_K\| = \alpha_n\gamma_ng(\tH)$ and $g(h)=\sqrt{(K-1)(K+3)/(16K)}$. As a consequence, for $n$ sufficiently large,
\begin{equation*}
    \frac{\|P_n\tH-\tilde{\alpha}_0\bm{1}_K\|}{\|P_n-\tilde{\alpha}_0{\bf 1}_K{\bf 1}_K'\|} = \frac{\alpha_n\gamma_ng(\tH)}{ \alpha_n\gamma_n} \xrightarrow[]{\mbox{as}}g(h) = \sqrt{\frac{(K-1)(K+3)}{16K}}> t_0.
\end{equation*}
It follows that for $n$ sufficiently large, with probability at least $1-o(1)$, 
\begin{equation*}
    \|P_n\tH-\tilde{\alpha}_0\bm{1}_K\|\geq t_0\|P_n-\tilde{\alpha}_0{\bf 1}_K{\bf 1}_K'\|,
\end{equation*}
which concludes the proof. \qed

\section{Proof of Propositions~\ref{prop:iden}-\ref{prop:INC}}\label{appendix:proofs_INC}

\subsection{Proof of Proposition~\ref{prop:iden}} \label{subsec:proof:iden}

We suppose that there exists an eligible tuple $(\Pi_0,P_0,K_0)$ such that $\Omega=\Pi_0P_0\Pi_0'$. To show the first point of the proposition, define the set:
\begin{equation*}
    S = \left\{k\in\mathbb{N}^* \quad \middle| \quad \exists (\Pi,P)\in\mathbb{R}^{n\times k}\times\mathbb{R}^{k\times k} \text{ eligible such that } \Omega=\Pi P\Pi'\right\}.
\end{equation*}
Note that $S$ is a discrete set lower bounded by $0$ which is non-empty since $K_0\in S$ by assumption. It follows that $S$ has a lower bound, which we denote as $k_\Omega$. It corresponds to the INC defined in Definition~\ref{def:INC}.\\

Now, we proceed to showing that when $K=k_\Omega$, the matrix $P$ is identifiable up to permutation. Suppose that we have two pairs of eligible matrices $(\Pi,P),(\Pi^*,P^*)\in\mathbb{R}^{n\times k_\Omega}\times\mathbb{R}^{k_\Omega\times k_\Omega}$ such that $\Omega=\Pi P\Pi'=\Pi^*P^*(\Pi^*)'$. Because $\Pi,\Pi^*$ are eligible, they contain the identity matrix as a submatrix. We assume without loss of generality that the first $k_\Omega$ rows of $\Pi$ and $\Pi^*$ correspond to $k_\Omega$ pure points, one per community. The submatrices
\begin{align*}
    \tPi := \Pi_{|\{1,...,k_\Omega\},\cdot} \quad \text{ and } \tPi^* := \Pi^*_{|\{1,...,k_\Omega\},\cdot}
\end{align*}
are permutations matrix. We have
\begin{align*}
    \Omega_{|\{1,...,k_\Omega\}\times\{1,...,k_\Omega\}} = \tPi P\tPi' = \tPi^*P^*(\tPi^*)',
\end{align*}
which implies that $P^*=DPD'$, where $D=(\tPi^*)'\tPi$ is a permutation matrix.\\

If, in addition, we have that $\text{rank}(P)=k_\Omega$, then $P$ is invertible. It follows that
\begin{align*}
    \tPi P=\tPi^*P^*(\tPi^*)'\tPi=\tPi^*P^*D=\tPi^*DP \quad \Longrightarrow \quad \tPi=\tPi^*D.
\end{align*}
In addition, since $\Omega=\Pi P\Pi'=\Pi^*P^*(\Pi^*)'$, $\Pi$ and $\Pi^*$ have full column rank, which means that there must exist an invertible matrix $B\in\mathbb{R}^{K\times K}$ such that $\Pi=\Pi^*B$. This implies that $\tPi=\tPi^*B$. As a result that $B=D$, so $\Pi=\Pi^*D$. This shows that if $\text{rank}(P)=k_\Omega$, then $\Pi$ is also identifiable up to permutation.\\

Finally, it holds by definition of $k_\Omega$ that $K_0\geq k_\Omega$. Since $\text{rank}(P_0)=\text{rank}(\Omega)$ and $k_\Omega=\text{dim}(P)\geq \text{rank}(\Omega)$, we obtain that 
\begin{equation*}
    K_0 \geq k_\Omega \geq \text{rank}(P_0).
\end{equation*}
Furthermore, if $P_0$ is non-singular, then $K_0=\text{rank}(P_0)$, hence
\begin{equation*}
    K_0 = k_\Omega = \text{rank}(P_0).
\end{equation*}

\subsection{Proof of Proposition~\ref{prop:INC}} \label{subsec:proof:INC}

By Proposition~\ref{prop:iden}, there exists a pair of eligible $\Pi\in\mathbb{R}^{n\times k_\Omega}$ and $P\in\mathbb{R}^{k_\Omega\times k\Omega}$ such that $\Omega=\Pi P\Pi'$, where $k_\Omega$ is the INC. Hence in the rest of the proof, we take $K=k_\Omega$.\\

Denote by $\Lambda\in\mathbb{R}^{r\times r}$ the matrix of eigenvalues of $\Omega$. It follows that we can write $\Omega=\Xi\Lambda\Xi'$. Furthermore, note that the fact that $r=\text{rank}(\Omega)$ implies that we also have $\text{rank}(P)=r$. We can thus denote by $X\in\mathbb{R}^{k_\Omega\times r}$ the matrix of eigenvectors of $P$, and by $L\in\mathbb{R}^{r\times r}$ the corresponding matrix of non-zero eigenvalues, thus obtaining that $P=XLX'$. As a consequence,
\begin{equation*}
    \Omega = \Xi\Lambda\Xi'=\Pi XL (\Pi X)'.
\end{equation*}
Note that $\Lambda\Xi'$ and $L (\Pi X)'$ must have full row-rank $r$, so the column space of $\Xi$ is equal to the column space of $\Pi X$. There must exist a matrix $B\in\mathbb{R}^{r\times r}$ such that $\Xi=\Pi XB$. Hence there exists a matrix $V\in\mathbb{R}^{k_\Omega\times r}$ such that
\begin{equation}\label{proofINC-1}
    \Xi = \Pi V.
\end{equation}

Since $\Pi$ is a membership matrix, it follows that the rows of $\Xi$ are convex combinations of the $k_\Omega$ rows of $V$. Because $\Pi$ is eligible, the identity matrix is a submatrix of $\Pi$. Without loss of generality, assume that $\Pi_{|\{1,...,k_\Omega\},\cdot}=I_{k_\Omega}$. It follows that $V=\Xi_{\{1,...,k_\Omega\},\cdot}$. This shows that $C(\Xi)$ is a polytope with at most $k_\Omega$ vertices and at least $r$ vertices.\\

In the case that $k_\Omega=r$, the desired result follows immediately. In the case that $k_\Omega<r$, \textit{hic jacet lepus}. Suppose by contradiction that $V$ has only $N$ distinct rows, where $r\leq N<k_\Omega$. This means that we can write $\Xi=\Pi B\tilde{V}$, where $\tilde{V}\in\mathbb{R}^{N\times r}$ is the matrix containing the unique rows of $V$ and $B\in\mathbb{R}^{k_\Omega\times N}$ is a row-replication matrix (which admits the identity matrix $I_N$ as a submatrix). It follows that we can write:
\begin{equation*}
    \Omega = \Pi B\tilde{V}\Lambda \tilde{V}'B'\Pi'.
\end{equation*}
We denote $\tPi:=\Pi B$ and $\tilde{P}=\tilde{V}\Lambda \tilde{V}'$, and proceed to showing that these matrices are eligible. First, it is straightforward to see that for any $i\in\{1,...,n\}$, the $i$-th row of $\tPi$ is positive and verifies $\tilde{\pi}_i'\bm{1}_N=\pi_i'B\bm{1}_N=\pi'\bm{1}_K=1$. In addition, since both $\Pi$ admits $I_{k_\Omega}$ as a submatrix and $B$ admits $I_N$ as a submatrix, it follows that $\tPi$ admits $I_N$ as a submatrix. This shows that $\tPi$ is admissible.

Now, from Equation~\eqref{proofINC-1}, we know that $\Omega=\Pi V\Lambda V'\Pi'$, so $P=V\Lambda V'=B\tilde{V}\Lambda\tilde{V}'B'$. By definition, $B$ admits a left inverse, call it $Q\in\{0,1\}^{N\times k_\Omega}$, so that $QB=I_2$. Then $\tilde{P}=QPQ'$. Since both $Q$ and $P$ are nonnegative, it follows that $\tilde{P}$ is nonnegative, thus eligible.\\

We have shown that we can write $\Omega=\tPi\tilde{P}\tPi'$, where $(\tPi, \tilde{P})\in\mathbb{R}^{n\times N}\times\mathbb{R}^{N\times N}$ eligible and $N<k_\Omega$, QEA. \qed


\end{document}